%% file: isoperimetric_nonlocal.tex
\documentclass[a4paper,reqno]{amsart}

\usepackage{ifxetex,ifluatex}
\newif\ifxetexorluatex
\ifxetex
  \xetexorluatextrue
\else
  \ifluatex
    \xetexorluatextrue
  \else
    \xetexorluatexfalse
  \fi
\fi

\usepackage{lmodern}
\usepackage{xargs} 

\ifxetexorluatex
	\usepackage{polyglossia}
	\setmainlanguage{english}

	\usepackage{fontspec}
	\usepackage{amsmath} 
	\usepackage{amsthm}
	\usepackage{mathtools,thmtools} 
	\usepackage[mathrm=sym,mathbf=sym]{unicode-math}

	\ifluatex

	\DeclareMathAlphabet{\mathcal}{OMS}{cmsy}{m}{n}

	\newcommandx{\forcebold}[2][1=0.2]{{\pdfliteral direct {2 Tr #1 w}#2\pdfliteral direct {0 Tr 0 w}}}
	
	\else
	\DeclareMathAlphabet{\mathcal}{OMS}{cmsy}{m}{n}

	\newcommandx{\forcebold}[2][1=0.2]{{\boldmath #2}}
	\fi

	\renewcommand\boldsymbol\mathbf

\else
	\pdfoutput=1
	\usepackage[english]{babel}

	\usepackage[T1]{fontenc}
	\usepackage[utf8]{inputenc}
	\usepackage[final]{microtype}
	\usepackage{amsmath}
	\usepackage{amsthm}
	\usepackage{amssymb}
	\usepackage{mathtools,thmtools} 

	\usepackage{bm}
	\renewcommand{\mathbf}[1]{\bm{\mathrm{#1}}}
	
	\providecommand{\mathbfscr}[1]{\bm{\mathscr{#1}}}
	\providecommand{\mathbfcal}[1]{\bm{\mathcal{#1}}}
	\newcommandx{\forcebold}[2][1=0.2]{\bm{{#2}}}

	\usepackage{scalerel,stackengine}
	\stackMath
	\newcommand\widecheck[1]{%
	\savestack{\tmpbox}{\stretchto{%
	  \scaleto{%
		\scalerel*[\widthof{\ensuremath{#1}}]{\kern-.6pt\bigwedge\kern-.6pt}%
		{\rule[-\textheight/2]{1ex}{\textheight}}
	  }{\textheight}%
	}{0.5ex}}%
	\stackon[1pt]{#1}{\scalebox{-1}{\tmpbox}}%
	}

	\usepackage{mathrsfs,esint}
\fi

\usepackage[top=3cm, bottom=3cm, inner=2.7cm, outer=2.7cm,headsep=0.75cm,footskip=1.25cm]{geometry}

\usepackage{csquotes,enumitem,multicol}

\usepackage{hyperref}
\hypersetup{unicode,hypertexnames=false,colorlinks=true,linkcolor=blue,citecolor=blue,pdfauthor={M. Pegon},pdftitle={Large mass minimizers for isoperimetric problems with integrable nonlocal potentials}}

\usepackage[capitalize,nameinlink,noabbrev]{cleveref}
\usepackage[backend=biber,sorting=nyt,style=numeric,isbn=false,url=false,doi=false,eprint=false,giveninits=true,maxbibnames=50]{biblatex}
\renewbibmacro{in:}{} 
\AtEveryBibitem{\iffieldundef{month}{}{\clearfield{month}}}
\AtEveryBibitem{\iffieldundef{day}{}{\clearfield{day}}}
\AtEveryBibitem{\clearfield{note}} 
\AtEveryBibitem{\clearlist{location}} 

\renewbibmacro{volume+number+eid}{%
    \textbf{
	\printfield{volume}%
    \setunit{\adddot}%
    \printfield{number}%
    \setunit{\addcomma\space}%
    \printfield{eid}}
}
\AtBeginBibliography{%

}

\DeclareFieldFormat
  [article,inproceedings,incollection,report,unpublished]
  {title}{#1\addcomma}
\DeclareFieldFormat
  [book,inbook,incollection,unpublished]
  {date}{(#1)}
\DeclareFieldFormat
  [book,inbook]
  {title}{\textit{#1}\addcomma}
\DeclareFieldFormat
  [article,inbook,book,incollection,inproceedings,report,unpublished]
  {pages}{#1}
\DeclareDelimFormat[bib]{nametitledelim}{\space:\space}
\DeclareNameAlias{labelname}{given-family} 
\usepackage{xifthen}
\usepackage[dvipsnames]{xcolor}
\usepackage{graphicx}
\usepackage[mode=buildnew]{standalone}
\usepackage{tikz}
\usepackage{appendix}
\usepackage{chngcntr}
\usepackage{etoolbox}
\usepackage[colorinlistoftodos,textsize=small,disable=true]{todonotes}
\newcommandx{\unsure}[2][1=]{\todo[linecolor=red,backgroundcolor=red!25,bordercolor=red,#1]{#2}}
\newcommandx{\change}[2][1=]{\todo[linecolor=blue,backgroundcolor=blue!25,bordercolor=blue,#1]{#2}}
\newcommandx{\info}[2][1=]{\todo[linecolor=OliveGreen,backgroundcolor=OliveGreen!25,bordercolor=OliveGreen,#1]{#2}}
\newcommandx{\improvement}[2][1=]{\todo[linecolor=Plum,backgroundcolor=Plum!25,bordercolor=Plum,#1]{#2}}

\newcommand{\nhphantom}[1]{\sbox0{#1}\hspace{-\the\wd0}}

\newcommand{\numberset}[1]{{\mathbb{#1}}}
\newcommand{\N}{\numberset{N}}

\newcommand{\R}{\numberset{R}}

\newcommand{\dd}{\mathop{}\mathopen{}\mathrm{d}} 
\newcommand{\dx}{\mathop{}\mathopen{}\mathrm{d}x} 
\newcommand{\dy}{\mathop{}\mathopen{}\mathrm{d}y} 
\newcommand{\dz}{\mathop{}\mathopen{}\mathrm{d}z} 
\newcommand{\dt}{\mathop{}\mathopen{}\mathrm{d}t} 
\newcommand{\dr}{\mathop{}\mathopen{}\mathrm{d}r}

\newcommand{\dH}{\mathop{}\mathopen{}d\mathscr{H}} 
 
\newcommand{\ddt}{\mathop{}\mathopen{}\frac{\mathrm{d}}{\mathrm{d}t}}

\newcommand{\bbS}{\mathbb{S}}

\newcommand{\bfK}{\mathbf{K}}

\newcommand{\calD}{{\mathscr{D}}}
\newcommand{\calE}{{\mathcal{E}}}
\newcommand{\calF}{{\mathcal{F}}}

\newcommand{\calH}{{\mathscr{H}}}
\newcommand{\calS}{{\mathcal{S}}}
\newcommand{\calM}{{\mathcal{M}}}
\newcommand{\calN}{{\mathcal{N}}}
\newcommand{\calL}{{\mathscr{L}}}
\newcommand{\calQ}{{\mathcal{Q}}}
\newcommand{\calR}{{\mathcal{R}}}
\newcommand{\calV}{{\mathcal{V}}}

\newcommand{\scrR}{{\mathscr{R}}}

\newcommand{\ind}{{\mathbf{1}}}
\newcommand{\Ren}{{\R^n}} 

\newcommand{\Om}{\Omega}
\newcommand{\om}{\omega}
\newcommand{\eps}{\varepsilon}
\newcommand{\vphi}{\varphi}

\newcommand{\ol}{\overline}
\newcommand{\wt}{\widetilde}

\newcommand{\bv}{\mathrm{BV}}

\newcommand{\loc}{{\mathrm{loc}}}
\newcommand{\aev}{a.e.}

\newcommand{\id}{{\mathrm{id}}}

\newcommand{\compl}{{\mathrm{c}}}

\newcommand{\grad}{\nabla}
\newcommand{\subsq}{\subseteq}

\newcommand{\csubset}{\subset\subset} 

\newcommand{\bmass}[2][]{\lbrack B\rbrack_{\ifthenelse{\isempty{#1}}{}{#1,}#2}}

\newcommand{\iintb}{\iint\limits} 

\newcommand{\overarc}[1]{\overset{\frown}{#1}}

\newcommand{\edit}[1]{#1}

\DeclareMathOperator{\Per}{Per}

\DeclareMathOperator{\Tan}{\mathrm{Tan}}

\DeclareMathOperator{\dvg}{\mathrm{div}}

\DeclareMathOperator{\dist}{\mathrm{d}}
\DeclareMathOperator{\spt}{\mathrm{spt}}

\DeclarePairedDelimiter{\abs}{\lvert}{\rvert}

\DeclarePairedDelimiter{\norm}{\lVert}{\rVert}
\DeclarePairedDelimiter{\seminorm}{\lbrack}{\rbrack}

\numberwithin{equation}{section} 

\declaretheorem[name=Theorem]{mainthm}

\declaretheorem[name=Application,numberlike=mainthm]{mainapp}

\declaretheorem[name=Theorem,within=section]{thm}
\declaretheorem[name=Lemma,numberlike=thm]{lem}
\declaretheorem[name=Proposition,numberlike=thm]{prp}
\declaretheorem[name=Corollary,numberlike=thm]{cor}
\declaretheorem[name=Definition,numberlike=thm,style=definition]{dfn}
\declaretheorem[name=Remark,numberlike=thm,style=remark]{rmk}

\crefname{equation}{}{}
\crefname{enumi}{}{}
\crefname{thm}{Theorem}{Theorems}
\crefname{mainthm}{Theorem}{Theorems}
\crefname{lem}{Lemma}{Lemmas}
\crefname{mainapp}{Application}{Applications}
\crefname{prp}{Proposition}{Propositions}
\crefname{cor}{Corollary}{Corollaries}
\crefname{dfn}{Definition}{Definitions}
\crefname{rmk}{Remark}{Remarks}
\crefname{conj}{Conjecture}{Conjectures}
\crefname{ex}{Example}{Examples}
\crefname{appsec}{Appendix}{Appendices}

\def\Xint#1{\mathchoice
{\XXint\displaystyle\textstyle{#1}}%
{\XXint\textstyle\scriptstyle{#1}}%
{\XXint\scriptstyle\scriptscriptstyle{#1}}%
{\XXint\scriptscriptstyle
\scriptscriptstyle{#1}}%
\!\int}
\def\XXint#1#2#3{{%
\setbox0=\hbox{$#1{#2#3}{\int}$}
\vcenter{\hbox{$#2#3$}}\kern-.5\wd0}}

\def\dashint{\Xint-}

\renewcommand{\leq}{\leqslant}
\renewcommand{\geq}{\geqslant}

\title[Large mass minimizers for isoperimetric problems]{\textbf{Large mass minimizers for
isoperimetric problems with integrable nonlocal potentials}}
\author{Marc Pegon}
\date{}

\begin{document}

\begin{abstract}
This paper is concerned with volume-constrained minimization problems derived from Gamow's liquid
drop model for the atomic nucleus, involving the competition of a perimeter term and repulsive
nonlocal potentials. We consider a large class of potentials, given by general radial nonnegative
kernels which are integrable on $\mathbb{R}^n$, such as Bessel potentials, and study the behavior of
the problem for large masses (i.e., volumes). Contrarily to the small mass case, where the nonlocal
term becomes \textsl{negligible} compared to the perimeter, here the nonlocal term \textsl{explodes}
compared to it.
However, using the integrability of those kernels, we rewrite the problem as the minimization of the
difference between the classical perimeter and a nonlocal perimeter, which converges to a
multiple of the classical perimeter as the mass goes to infinity.
Renormalizing to a fixed volume, we show that, if the first moment of the kernels is smaller than an
explicit threshold, \textsl{the problem admits minimizers of arbitrarily large mass}, which
contrasts with the usual case of Riesz potentials.
In addition, we prove that, any sequence of minimizers converges to the ball as the mass goes to
infinity.
Finally, we study the stability of the ball, and show that our threshold on the first moment of
the kernels is \textsl{sharp} in the sense that large balls go from stable to unstable.
A direct consequence of the instability of large balls above this threshold is that there exist
nontrivial compactly supported kernels for which the problems admit minimizers which are not balls,
that is, \textsl{symmetry breaking occurs}.
\end{abstract}


\maketitle
\vspace{-0.8cm}

\tableofcontents

\vspace{-1cm}

\section{Introduction}
\addtocontents{toc}{\protect\setcounter{tocdepth}{1}}


We study large mass minimizers for a variant of Gamow's liquid drop model for the atomic nucleus,
in which the repulsive term is given by a general nonnegative, integrable, radial kernel.
More precisely, given $G:\Ren\setminus\{0\}\to[0,+\infty)$ a measurable nonnegative radial function
with $G\in L^1(\Ren)$ and $n\geq 2$, we originally consider the minimization problem
\begin{gather}\label{origpb}\tag{$\star$}
\min\left\{P(E) + \iint_{E \times E} G(x-y) \dd x \dd y ~:~ \abs{E} = m\right\},
\end{gather}
where the minimum is taken over all sets of finite perimeter of volume $\abs{E}=m$ -- which we call
the mass -- and $P(E)$ denotes the perimeter of $E$.
Observe that this problem exhibits a competition between two terms and is thus nontrivial: the local
perimeter term constrains the set $E$ to concentrate as much as possible, while the nonlocal term
acts as a repulsive term, forcing $E$ to spread. Indeed, it is known that the perimeter is
\textit{minimized} by balls under volume constraint, while the nonlocal term is \textit{maximized}
by balls \edit{if $G$ is in addition radially nonincreasing\footnote{\edit{In fact, even if $G$ is
not radially decreasing, $G_\lambda$ \enquote{concentrates} near the origin when $\lambda$ is large,
heuristically making the nonlocal term repulsive if $G$ does not vanish identically.}} (by Riesz'
symmetric rearrangement, using e.g. \autocite[Chapter 3.7]{LiebLoss} and the fact that $G$ is equal
to its symmetric rearrangement in that case).}

As we show just below, the integrability assumption on $G$ allows us to reformulate this problem as
\begin{equation}\label{minrp}\tag{$\text{P}_{\gamma,\lambda}$}
\min \Bigg\{ P(E)-\gamma\Per_{G_\lambda}(E)~:~\abs{E}=\abs{B_1}\Bigg\},
\end{equation}
where $B_1$ stands for the open unit ball of $\Ren$ centered at the origin, $\lambda$ and $\gamma$
are positive constants, the kernel $G_\lambda$ is given by
\begin{gather}\label{defglambda}
G_\lambda(\cdot)\coloneqq\lambda^{n+1}G(\lambda\,\cdot),
\end{gather}
and the functional $\Per_G$ is defined by
\[
\Per_G(E)\coloneqq\iint_{E \times E^\compl} G(x-y) \dd x \dd y.
\]
The parameter $\lambda>0$ will represent the mass (to the power $\frac{1}{n}$), and $\gamma>0$ will
be chosen to adjust the first moment of $G$ (that is, its integral against the measure $\abs{x}\dx$)
to a particular universal constant.
As will be justified later on, $\Per_G$ should be considered as a \enquote{nonlocal perimeter},
which behaves in several ways as a standard perimeter term rather than as a volume term.

Before elaborating on the reformulation, let us say a few words on the original problem
\cref{origpb}. This problem has been studied extensively in the literature when $G$ is a Riesz
kernel (whose definition is recalled in \cref{subsec:riesz}) or a general integrable kernel with
compact support. In the Riesz case, it is known that the problem admits the ball as unique minimizer
below a critical mass, up to translations, and it is conjectured that there is no minimizer above a
(possibly different) critical mass. This conjecture has already been proven in a few cases. We will
discuss the \edit{Riesz and compact support cases} further, and provide references, in the next
section.

In this paper, we are interested in kernels decreasing faster than Riesz kernels at infinity, enough
to make them integrable, but not necessarily compactly supported, and we focus exclusively in the
case of large masses. Let us remark that, even if the kernel decreases rapidly at infinity, the
asymptotic behavior of the problem for large masses is very different than that of small masses:
indeed, in the case of Riesz (or Bessel) kernels, \textsl{the nonlocal term is negligible compared
to the perimeter} as the mass vanishes, so that the problem consists of minimizing the perimeter
plus a vanishing perturbation; here, the nonlocal term \textsl{explodes} compared to the perimeter
as the mass goes to infinity, \edit{as can be seen} by writing
\[
\begin{aligned}
P(E)+\iint_{E\times E} G(x-y)\dx\dy
&=\lambda^{n-1}\left(P(F)+\lambda^{n+1}\iint_{F\times F} G\big(\lambda(x-y))\dx\dy\right)\\
&=\lambda^{n-1}\Big(P(F)+\lambda\norm{G}_{L^1(\Ren)}-\Per_{G_\lambda}(F)\Big),
\end{aligned}
\]
where $E\coloneqq \lambda F$ with $\abs{F}=\abs{B_1}$, since $\Per_{G_\lambda}(F)$ is of the same
order as $P(F)$, as we will see.


\subsection*{Reformulation of the problem}

Rewriting the nonlocal repulsive term as
\[
\iint_{E\times E} G(x-y)\dx\dy = m\norm{G}_{L^1(\Ren)}-\iint_{E\times E^\compl} G(x-y)\dx\dy,
\]
we see that \cref{origpb} is in fact strictly equivalent to
\[
\min \Bigg\{ P(E)-\Per_G(E)~:~\abs{E}=m\Bigg\}.
\]
Now, to further normalize our problem, for $k\in \{0,1\}$, we define the quantity
\[
I_G^{k} \coloneqq \int_{\R^n} \abs{x}^k G(x)\dx,
\]
and for every positive natural numbers $p$ and $n$, we denote by $\bfK_{p,n}$ the constant defined
by
\begin{equation}\label{eq:defKpn}
\bfK_{p,n} \coloneqq \dashint_{\bbS^{n-1}} \abs{e\cdot x}^p\dH^{n-1}(x),
\end{equation}
which does not depend on $e\in \bbS^{n-1}$ by symmetry. Here $\bbS^{n-1}$ denotes the unit sphere in
$\Ren$ and $\calH^{n-1}$ the $(n-1)$-dimensional Hausdorff measure.
Now, up to dividing $G$ by the constant $\gamma=(I_G^1\bfK_{1,n})/2$, without renaming $G$, we may
assume that
\begin{equation}\label{eq:defIG}
I_G^1 = \frac{2}{\bfK_{1,n}},
\end{equation}
and may look instead at the problem
\begin{equation}\label{tmpMinPb}
\min \Bigg\{ P(E)-\gamma\Per_G(E)~:~\abs{E}=m\Bigg\},
\end{equation}
where $\gamma>0$. The choice of the constant in \cref{eq:defIG} will be justified in
\cref{subsec:reform}.
Keep in mind that the parameter $\gamma$ plays the role of a constant times $I_G^1$, the first
moment of $G$. As a last simplification step, in order to study the asymptotic behavior when the
mass goes to infinity, it is more convenient to look at the rescaled problem with fixed mass equal
to the volume of the unit ball in $\Ren$.
Given $m>0$, setting $\lambda\coloneqq\left(\frac{m}{\abs{B_1}}\right)^{\frac{1}{n}}$ and
$F:=\lambda^{-1}E$, it is then easy to see that $\abs{F}=\abs{B_1}$, $I_{G_\lambda}^{1}=I_G^1$, and
by a change of variables
\[
P(E)-\gamma\Per_G(E)=\lambda^{n-1}\big(P(F)-\gamma\Per_{G_\lambda}(F)\big),
\]
so that the set $E$ is a minimizer of \cref{tmpMinPb} if and only if $F$ is a minimizer of the
rescaled problem \cref{minrp}. Let us emphasize that even in the reformulation \cref{minrp}, the
nonlocal perturbation still \textsl{does not vanish} as $\lambda$ goes to infinity, contrarily to
the Riesz case for small masses, and our results hold for \textsl{any} $\gamma\in(0,1)$.
In the rest of the paper, we shall always work with this equivalent formulation.

\subsection*{General assumptions and results}

Except in \cref{sec:context}, we shall always assume that $G$ satisfies the two following general
hypotheses:
\begin{enumerate}[label=(H\arabic*),parsep=0pt]
\item\label{Krad} $G$ is radial, that is, there exists a nonnegative function $g: (0,+\infty)\to \R$ such that $G(x)=g(\abs{x})$ for $\calL^n-\aev$ $x\in \Ren$;
\item\label{Kint} $G\in L^1(\Ren)$, and
\[
I_G^1 = \frac{2}{\bfK_{1,n}}.
\]
\end{enumerate}
Starting from \cref{sec:stabball}, dedicated to the study of the stability of the ball, we may add
the extra assumptions
\begin{enumerate}[label=(H\arabic*)]
\setcounter{enumi}{2}
\item\label{Korig} $G(x)=o(\abs{x}^{\alpha-n})$ near the origin, for some $\alpha>0$;
\item\label{Ksuperdec} $G(x)=o(\abs{x}^{-(n+\beta)})$ at infinity for some $\beta>0$ when $n\geq 3$, and
$G(x)=o(\abs{x}^{-3})$ at infinity when $n=2$;
\edit{\item\label{KC1} $G\in C^1(\Ren\setminus\{0\})$.}
\end{enumerate}

\edit{These extra assumptions} are required essentially in order to be able to use directly
computations from \autocite{FFMMM} for the second variation of the nonlocal perimeter. The stronger
assumption \cref{Ksuperdec} in dimension~$2$ is due to the presence of a jacobian determinant when
integrating on the sphere, which appears to be singular only in dimension $2$ (see
\cref{lem:intsphere}). As we will see, Bessel kernels satisfy those general assumptions (see
\cref{subsec:prelimbessel,subsec:bessel}).

We are interested in the asymptotic behavior of the minimization problem \cref{minrp} as $\lambda$
(that is, the mass) goes to infinity, and give answers to several natural questions: does
\cref{minrp} admit a minimizer? If so, what do minimizers look like, are they regular? Can the unit
ball be a minimizer? We decided to state in a concise manner just below three of the main results
obtained in the paper, and an application to Bessel kernels, however these results are not
necessarily arranged in the same way in the paper.

\begin{mainthm}[Consequence of \cref{thm:existmin,thm:regminim}]\label{mainthm:existmin}
Assume $\gamma<1$. Then there exists $\lambda_e=\lambda_e(n,\gamma,G)$
such that, for any $\lambda > \lambda_e$, \cref{minrp} admits a minimizer, and in addition,
minimizers have a $C^{1,\frac{1}{2}}$ reduced boundary and are essentially connected.
\
\info{The dependency of $m_e$ in $G$ is not known. It depends on the speed of convergence of
$\eta(\lambda)$ to zero.}
\end{mainthm}

Let us point out that the regularity of the reduced boundary is actually true for any minimizer, if
one exists, no matter the value of $\gamma$ and $\lambda$. Connectedness of minimizers also holds in
any case, provided that $G$ is \edit{strictly positive} (see \cref{thm:regminim}).
In terms of the original problem \cref{origpb}, \cref{mainthm:existmin} means that if
$I_G^1<\frac{2}{\bfK_{1,n}}$, then there exists a critical mass $m_e=m_e(n,G)$, such that, above
this mass, the problem admits a minimizer.
To our knowledge, this is the first time existence of minimizers of arbitrarily large mass is
obtained in a Gamow-type problem on the whole space for non-compactly supported kernels without the
presence of an external attractive background potential, as is often the case (see e.g.
\autocite{ABCT2017,ABCT2018,GO2018}). \edit{Let us also mention \autocite{ABTZ2021}, where it is shown
that, if the perimeter is weighted by a power-law density growing sufficiently fast at infinity,
then minimizers always exist, and are balls in the large mass regime.}

The main obstacle for proving existence with the direct method in the calculus of variations is 
the possibility for a minimizing sequence to have some mass escape at infinity. We solve this
problem by showing that, for large values of $\lambda$, a minimizing sequence may be constrained
inside a ball via a truncation lemma. This relies heavily on the fact that the nonlocal perimeter
behaves to some extent as the classical perimeter and \edit{converges} to it as $\lambda$ goes to
infinity.
Note that the general kernels we consider (and in particular Bessel kernels) do not behave as nicely
as Riesz kernels under scaling, which are homogeneous.
The $C^{1,\frac{1}{2}}$-regularity of the reduced boundary of minimizers follows by results in
\autocite{Rig} on quasi-minimizers for the perimeter.

For any $\gamma$ and $\lambda$, let us define the functional to be minimized
\begin{equation}\label{deffuncF}
\calF_{\gamma,G_\lambda}(E) := P(E)-\gamma\Per_{G_\lambda}(E).
\end{equation}
When $\gamma\in[0,1)$, we are able to compute the $\Gamma$-limit in $L^1$ of these functionals as
$\lambda$ goes to infinity, and we show independently that, up to translations, any sequence of
minimizers converges to the unit ball. More precisely, we prove that minimizers are included in the
set~difference between two balls whose radii converge to $1$ as $\lambda$ goes to infinity, which
implies in particular Hausdorff convergence of the boundaries of minimizers to the unit sphere, up
to translations.

\begin{mainthm}[See \cref{thm:existmin}, \cref{thm:gammacv} for more precise
statements]\label{thm:cvminim}
Assume $\gamma < 1$.
For any minimizer $E$ of \cref{minrp} with $\lambda>\lambda_e$, up to a translation, we have
\begin{equation}\label{eq:dblincl}
\ol{B}_{1-\eta_{\gamma,G}(\lambda)} \subsq E \subsq B_{1+\eta_{\gamma,G}(\lambda)},
\end{equation}
where $\eta_{\gamma,G}$ is a function depending only on $n$, $\gamma$ and $G$ which vanishes at
infinity.
In addition, the family of functionals $\calF_{\gamma,G_\lambda}$ (with the added constraint of
being the indicator function of a set of finite perimeter with volume $\abs{B_1}$)
$\Gamma$-converges in $L^1$ to $\left(1-\gamma\right)P$ as $\lambda$ goes to infinity.
\end{mainthm}

\edit{Of course, the $\Gamma$-convergence of the functional to a positive multiple of the perimeter
implies that any converging sequence of minimizers of \cref{minrp} with $\lambda\to\infty$ converges
to the unit ball, but \cref{eq:dblincl} is stronger and in fact a direct consequence of the proof of
existence above $\lambda_e$.}

Then we recall a well-suited notion of stability for functionals on sets under volume constraint
(see \cref{dfn:stability}), and show that the threshold $\gamma=1$ is a stability threshold of
the unit ball for \cref{minrp} for large values of $\lambda$.

\begin{mainthm}[Consequence of \cref{thm:unstable,thm:stable}]\label{thm:stabthreshold}
Assume that $G$ satisfies all the hypotheses \crefrange{Krad}{KC1}. Then the following holds:
\begin{enumerate}[label=(\roman*),wide, labelwidth=0.5pt, labelindent=0pt]
\item if $\gamma<1$, then there exists $\lambda_s=\lambda_s(n,\gamma,G)>0$ such that, for any
$\lambda>\lambda_s$, $B_1$ is a (critical) stable set for $\calF_{\gamma,G_{\lambda}}$;
\item if $\gamma>1$, then there exists $\lambda_u=\lambda_u(n,\gamma,G)$ such that, for any
$\lambda>\lambda_u$, $B_1$ is a (critical) unstable set for $\calF_{\gamma,G_{\lambda}}$: in
particular, it cannot be a minimizer, i.e., symmetry-breaking occurs.
\end{enumerate}
\end{mainthm}

In terms of the original problem \cref{origpb}, this means that the threshold $\frac{2}{\bfK_{1,n}}$
for $I_G^1$ is a threshold for which large balls go from stable to unstable.

The proofs for the stability and instability of the ball rely essentially on the two following
ingredients:
\begin{enumerate}[label=(\roman*),parsep=0pt,itemsep=0.5pt]
\item the decomposition in spherical harmonics of the Jacobi operator associated with the second
variation of the perimeter and of the nonlocal term (given by the so-called Funk-Hekke formula for
the latter);
\item results analogous to the one by \citeauthor{BBM2} in \autocite{BBM2} for Sobolev spaces on
spheres, that is, computation of the limit and of a sharp \enquote{asymptotic} upper bound for the
quantity
\[
\iint_{\bbS^{n-1}\times \bbS^{n-1}} \frac{\abs{f(x)-f(y)}^2}{\abs{x-y}^2}
\eta_\eps(x-y)\dH^{n-1}_x\dH_y^{n-1},
\]
where $(\eta_\eps)_{\eps>0}$ is a \textsl{$(n-1)$-dimensional approximation of identity}, and $f$
belongs to $H^1(\bbS^{n-1})$.
\end{enumerate}

A particularly interesting consequence of \cref{thm:stabthreshold} is that there exist kernels for
which \cref{minrp} admits nontrivial minimizers, that is, minimizers which are not balls. Indeed,
working from the formulation \cref{origpb}, \citeauthor{Rig} proved in \autocite{Rig} that
\cref{minrp} \textsl{always} admits a minimizer whenever $G$ is compactly supported. Hence,
taking $\gamma>1$, minimizers still exist but cannot be the unit ball when $\lambda$ is large
enough, since it is unstable.

\cref{mainthm:existmin,thm:cvminim,thm:stabthreshold} directly apply when $G$ is a so-called Bessel
kernel: for every $\alpha,\kappa>0$, we denote by $\mathcal{B}_{\kappa,\alpha}$ the Bessel kernel of
order $\alpha$ defined as the fundamental solution of the operator
$(I-\kappa\Delta)^{\frac{\alpha}{2}}$, that is,
\[
(I-\kappa\Delta)^{\frac{\alpha}{2}}\mathcal{B}_{\kappa,\alpha} = \mathbf{\delta}_0
\qquad\text{ in }\calD'(\Ren),
\]
where $\mathbf{\delta}_0$ is the Dirac distribution at the origin. We then have:

\begin{mainapp}[\cref{cor:assumbessel,lem:expKpn}]\label{mainapp:bessel}
For every $\kappa,\alpha\in (0,+\infty)$, we consider the problem \cref{origpb} with
$G=\mathcal{B}_{\kappa,\alpha}$.
Let us define \[
\mathbf{\kappa}_{\alpha} \coloneqq 
\pi\left(\frac{(n+1)\Gamma\left(\frac{\alpha}{2}\right)}
{2\Gamma\left(\frac{1+\alpha}{2}\right)}\right)^2.
\]
Then we have:
\begin{enumerate}[label=(\roman*),wide, labelwidth=0.5pt, labelindent=0pt]
\item if $\kappa<\mathbf{\kappa}_\alpha$, there exist $m_e=m_e(\alpha,\kappa)>0$ and
$m_s=m_s(\alpha,\kappa)>0$ such that, for every $m>m_e$, \cref{origpb} admits a minimizer, and for
every $m>m_s$, the ball of volume $m$ centered at the origin, denoted by $\bmass{m}$, is a
\textbf{stable} critical point for the functional of \cref{origpb}. In addition, rescaling
minimizers so that they are of volume $\abs{B_1}$, and translating them, they converge to the unit
ball as $m$ goes to infinity;
\item if $\kappa>\mathbf{\kappa}_\alpha$, there exists $m_u=m_u(\alpha,\kappa)$ such that
for every $m>m_u$, $\bmass{m}$ is an \textbf{unstable} critical point of the functional of
\cref{origpb}. In particular, $\bmass{m}$ cannot be a minimizer.
\end{enumerate}
\end{mainapp}

In view of \crefrange{mainthm:existmin}{thm:stabthreshold}, we conjecture that for
$\gamma<1$, there should be a critical value $\lambda_B$ such that, for $\lambda>\lambda_B$,
the unique minimizer of \cref{minrp} is the unit ball, up to translations. This conjecture will be
the subject of a future work.

\subsection*{Acknowledgments}

The author is very grateful to his PhD advisor \textsc{V. Millot} for interesting him with this
problem, and to \textsc{B.  Merlet} for his suggestions to improve the clarity of the paper.
M. Pegon is supported by the Labex CEMPI (ANR-11-LABX-0007-01).

\subsection*{Outline of the paper}
This paper is organized as follows. In \cref{sec:context} we discuss a few variants of Gamow's
liquid drop model which have already been studied in the literature, and we motivate the choice of
our assumptions \cref{Krad,Kint}. We also recall some well-known results on isoperimetric
inequalities. In \cref{sec:prelim} we establish basic prerequisites on nonlocal perimeters and
on Bessel kernels, which justify \cref{mainapp:bessel}.
\cref{sec:exist} is devoted to the proofs of \cref{mainthm:existmin,thm:cvminim}. First, we prove
existence of minimizers for $\gamma<1$ and $\lambda$ large enough, as well as convergence to the
ball as $\lambda$ goes to infinity in \cref{thm:existmin}. Then, we compute the $\Gamma$-limit of
the functionals $\calF_{\gamma,G_\lambda}$ as $\lambda$ goes to infinity, and conclude this section
by establishing $C^{1,\frac{1}{2}}$-regularity (applying directly results from \autocite{Rig}) and
connectedness of minimizers.
In \cref{sec:stabball}, we focus on the stability of the unit ball for large $\lambda$, and show that
$\gamma=1$ is a threshold for which the unit ball goes from stable to unstable, i.e.,
\cref{thm:stabthreshold}. To conclude on the stability issue, we need to study asymptotics for some
nonlocal seminorms on the sphere, which is done in \cref{app:mollifsphere}: here we compute the
limit of these seminorms as the kernels concentrate to the Dirac distribution, and obtain a uniform
upper bound which is asymptotically sharp.

\subsection*{Notation}

\paragraph{\textsl{Operations on sets.}}
For any set $E\subsq\Ren$, we define $E^\compl\coloneqq\Ren\setminus E$, and we write $\abs{E}$ for its
volume (that is, its Lebesgue measure) whenever $E$ is measurable. We write $E\sqcup F$ for the
union of two sets which are disjoint. Given two sets $E$ and $F$, we denote by $E\triangle
F\coloneqq (E\setminus F)\sqcup(F\setminus E)$ their symmetric difference. We say that two sets $E$
and $F$ in $\Ren$ are equivalent if $\abs{E\triangle F}=0$.

\vspace{10pt}

\paragraph{\textsl{Hausdorff measures.}}
We denote by $\calH^k$ the $k$-dimensional Hausdorff measure in $\Ren$, and by $\dim_\calH(E)$ the
Hausdorff dimension of a set $E\subsq\Ren$. When integrating w.r.t. the measure $\calH^k$ in a
variable $x$, we use the notation $\dH_x^k$ instead of the more standard but less compact
$\dH^k(x)$.

\vspace{10pt}

\paragraph{\textsl{Balls and spheres.}}
We denote by $B_r(x)$ the open ball in $\Ren$ of radius $r$ centered at $x$. For simplicity we write
$B_r$ when $x$ is the origin.
The volume of $B_1$ is $\om_n\coloneqq \abs{B_1} =
\frac{\pi^{\frac{n}{2}}}{\Gamma\left(1+\frac{n}{2}\right)}$, and the area of the unit sphere
$\bbS^{n-1}$ is $\calH^{n-1}(\bbS^{n-1})=n\om_n$. More generally we denote by $\bbS^k$ the
$k$-dimensional unit sphere, and for simplicity we write $\abs{\bbS^k}=\calH^k(\bbS^k)$ for its
surface area.
For any $m>0$ and $x\in\Ren$ we let $\bmass[x]{m}$ be the open ball of volume $m$ centered at $x$,
or simply $\bmass{m}$ if $x=0$.

\vspace{10pt}

\paragraph{\textsl{Sets of finite perimeter.}}

For any nontrivial open set $\Om\subsq\Ren$, we denote by $\bv(\Om)$ the space of functions with
bounded variation in $\Om$, and for any $f\in \bv(\Om)$ we let $\abs{Df}$ be its total variation
measure, and set $\seminorm{f}_{\bv(\Om)}\coloneqq\int_{\Om}\, \abs{Df}$.
For a set of finite perimeter $E$ in $\Om$, we let $\ind_E\in\bv(\Om)$ be its characteristic
function (i.e., $\ind_E(x)=1$ if $x\in E$ and $0$ otherwise), and define its perimeter in $\Om$ by
$P(E;\Om)\coloneqq \int_{\Om}\,\abs{D\,\ind_E}$. If $\Om=\Ren$ we simply write $P(E)\coloneqq P(E;\Ren)$.
We denote by $\mu_E\coloneqq D\,\ind_E$ the Gauss--Green measure associated with the set of finite
perimeter $E$ and $\nu_E(x)$ the outer unit normal of $\partial^* E$ at $x$, where $\partial^*E$
stands for the reduced boundary of~$E$.
We refer e.g. to \autocite[Chapter~5]{EG2015} or \autocite{Maggi} for further details on functions
of bounded variations and sets of finite perimeter.

\addtocontents{toc}{\protect\setcounter{tocdepth}{2}}

\section{Motivation and context}\label{sec:context}

\subsection{No repulsion: the classical isoperimetric problem}

First let us say a few words about the simplest case for \cref{origpb}, that is, when $G\equiv 0$.
In that case, \cref{origpb} is the classical isoperimetric problem which consists in minimizing the
perimeter under a volume constraint. It is known that the unique minimizer is the ball, up to
translations (see e.g. \autocite{De2006a}), which gives the classical isoperimetric inequality
\[
P(E) \geq P(\bmass{m}),
\]
for any set of finite perimeter $E$ with volume $m$, and can be rewritten
\begin{equation}\label{eq:classicisop}
P(E) \geq n\om_n^{\frac{1}{n}} \abs{E}^{1-\frac{1}{n}}.
\end{equation}
Knowing that balls are solutions to the classical isoperimetric problem, it is then natural to
consider the related question: if the perimeter of a set $E$ of volume $m$ is close to
$P(\bmass{m})$, is $E$ close to the ball $\bmass{m}$ in some sense, and if so, is it possible to
quantify it? An answer to this question has been given in \autocite{FMP} (see also \autocite{FV} for
a refinement), in the form of a so-called quantitative version of the isoperimetric inequality,
which we recall just below.
Given a set with finite perimeter $E$ such that $\abs{E}=m$, we define the \textsl{isoperimetric
deficit} of $E$ by
\[
D(E) \coloneqq \frac{P(E)-P(\bmass{m})}{P(\bmass{m})},
\]
and its \textsl{Fraenkel asymmetry} by
\[
\alpha(E) \coloneqq \min \left\{ \frac{\abs{E\triangle \bmass[x]{m}}}{m}~:~ x\in\Ren\right\}.
\]
The sharp quantitative isoperimetric inequality proven in \autocite{FMP} then states that there
exists a constant $C=C(n)$ such that
\begin{equation}\label{eq:quantisop}
\alpha(E) \leq C\sqrt{D(E)},
\end{equation}
and that the $\frac{1}{2}$ exponent over $D(E)$ is sharp. In addition to their intrinsic interest,
isoperimetric inequalities are a very useful tool to study related isoperimetric problems, and we
shall often rely on them in the rest of the paper.

\subsection{Slow decay at infinity: Riesz potentials}\label{subsec:riesz}

Problems such as \cref{origpb} are essentially inspired by a simple model for the atomic nucleus
introduced by George Gamow in the late 1920s, which is now referred to as Gamow's liquid drop model.
This denomination is due to the fact that in this simple model (then refined by Heisenberg, von
Weizsäcker and Bohr in the 1930s), the protons and neutrons inside the atomic nucleus are treated as
an incompressible and uniformly charged fluid. In this model, the atomic nucleus is represented by a
set $\Om\subsq \R^3$ of volume $m$ (which corresponds to its mass), and its energy is given by
\[
P(\Om)+\frac{1}{8\pi}\iint_{\Om\times\Om} \frac{1}{\abs{x-y}}\dx\dy.
\]
The perimeter term represents the energy associated with the attractive short-range nuclear force,
while the Coulombic repulsive term is due to the positively charged protons pushing themselves away
from each other. This model successfully explained the phenomenon of nuclear fission: indeed, there
are two critical masses $0<m_1\leq m_2<\infty$ such that, below $m_1$, the problem admits a
minimizer (no fission), and above $m_2$, there is no minimizer (fission). In fact, there exists
another threshold $0<m_0\leq m_1$ such that, below it, the ball is the unique minimizer (up to
translations). These results were first rigorously proven in \autocite{KM2}.
Many variants and generalizations of this model have been proposed since then, one of the most
natural being to replace the Newton potential $\frac{1}{\abs{x}^{n-2}}$ in dimension $3$ with Riesz
potentials in arbitrary dimension $n\geq 2$, that is
\[
G(x)=\frac{1}{\abs{x}^{n-\alpha}},\quad \alpha\in (0,n).
\]
The Newton case $\alpha=2$ in dimension $n\geq 3$ was treated e.g. in \autocite{Jul}, the Riesz
cases in dimension $2$ for $\alpha\in(0,2)$ in \autocite{KM1}, in arbitrary dimension for $\alpha\in
(0,n-1)$ in \autocite{BC}, and finally the complete Riesz case in arbitrary dimension for every
$\alpha\in (0,n)$ in \autocite{FFMMM}, where the perimeter $P(E)$ can also be replaced by the
$s$-fractional perimeter $P_s(E)$ with $s\in(0,1)$.

Let us sum up some of what is known in the Riesz case in the following theorems.

\begin{thm}[{\autocite{KM1,KM2,Jul,BC,FFMMM}}]
Given $n\geq 2$ and $\alpha\in (0,n)$, there exists $m_0=m_0(n,\alpha)$ such that for any $m<m_0$
the ball $\bmass{m}$ is the unique minimizer of \cref{origpb} for $G(x)=\abs{x}^{-(n-\alpha)}$, up
to translations.
\end{thm}

There are also some nonexistence results.

\begin{thm}[{\autocite{BC,KM1,KM2,LO}}]
Given $n\geq 2$ and $\alpha\in (n-2,n)$, there exists $m_1=m_1(n,\alpha)$ such that for any $m>m_1$,
\cref{origpb} admits no minimizer for $G(x)=\abs{x}^{-(n-\alpha)}$.
\end{thm}

These nonexistence results for large masses are in a sense not surprising. Indeed, on the one hand,
note that without the perimeter term the problem
\[
\min \left\{ \iint_{E\times E} \frac{1}{\abs{x-y}^{n-\alpha}} \dx\dy ~:~ \abs{E}=m\right\}
\]
admits no minimizer, as it is always better to split a set $E$ into infinitely many pieces and
send them farther from each other at infinity, since Riesz kernels are (strictly) radially decreasing.
On the other hand, the relatively slow decay at infinity of the Riesz kernels make them
nonintegrable, which would explain why the repulsive potential takes over the perimeter term in
\cref{origpb} for large masses, resulting in the nonexistence of minimizers.

As for the thresholds $m_0$, $m_1$, and $m_2$, physical evidence indicate that in dimension $n=3$ at
least, they should be equal, but this has yet to be proven.

\subsection{Compactly supported kernels}

An interesting case is when the kernel $G$ has compact support, i.e., when the long-range
interaction fully disappears at some distance.
Recalling our informal discussion on nonexistence of minimizers for Riesz potentials, we see that in
the compact case, sending disjoint pieces of a set $E$ at infinity does not decrease the energy of
the nonlocal term: when the pieces are far enough, they simply do not interact with each other. Thus
we may imagine that it is possible to build a minimizing sequence lying in a fixed ball, and prove
existence of minimizers by the direct method.
In dimension $n=2$, this strategy could probably be implemented quite readily, provided one controls
the number of connected components of a minimizing sequence (the advantage being that sets of finite
perimeter are essentially bounded, that is, included in a ball), but in higher dimension it is much
more complex.

Fortunately, using the link between minimizers of \cref{origpb} and quasi-minimizers of the
perimeter (see \cref{subsec:reg}), that case was successfully treated by \citeauthor{Rig} in
\autocite{Rig}, yielding the following result.

\begin{thm}[{\autocite{Rig}}]\label{mainthm:existmincompact}
If $G$ is compactly supported, then \cref{origpb} always admits minimizers. In addition, for any
minimizer $E$, $\partial^* E$ is a $C^{1,\frac{1}{2}}$-hypersurface, and, up to choosing a good
representative, $E$ has a finite number of connected components $N$, where $N$ can be bounded
depending only on $G$, $n$ and $m$.
\end{thm}

Note that \cref{mainthm:existmincompact} stands true even if $G$ is not radial. 
A consequence of this theorem and of \cref{thm:stabthreshold} is that we can easily
build kernels such that large mass minimizers exist and are nontrivial, as mentioned in the
introduction.

\subsection{Intermediate case: Bessel kernels}\label{subsec:prelimbessel}

Between Riesz kernels, which are slowly decreasing kernels, and compactly supported kernels, it is
natural to wonder what happens in the intermediate case of rapidly decreasing kernels such as
Bessel kernels.
Bessel kernels are usually given by the operators $(I-\Delta)^{-\frac{\alpha}{2}}$ for
$\alpha\in (0,n)$, i.e., the Bessel kernel of order $\alpha$ is the fundamental solution of
\[
(I-\Delta)^{\frac{\alpha}{2}}f = \mathbf{\delta}_0\quad\text{ in }\calD'(\Ren),
\]
where $\mathbf{\delta}_0$ is the Dirac distribution at the origin.
In fact, we consider the ``generalized'' Bessel kernels given by
$(I-\kappa\Delta)^{-\frac{\alpha}{2}}$, where $\alpha,\kappa\in(0,+\infty)$. As far as we know,
there is little literature on \cref{origpb} when $G$ is a Bessel kernel, and especially on
the asymptotic behavior for large masses.
Compared with Riesz kernels (which are associated with the operators
$(-\Delta)^{-\frac{\alpha}{2}}$), Bessel kernels are generally not explicit, in the sense that they
only have an integral representation, and they do not behave as nicely as Riesz kernels under
scaling. Near the origin, Riesz and Bessel kernels of the same order $\alpha$ behave similarly,
however at infinity Bessel kernels decay much faster. Their decay at infinity is
\textsl{exponential} (in particular, they are integrable), making them an intermediate case between
Riesz kernels and compactly supported kernels.

Physically, Bessel kernels are suggested in \autocite{KMN} to model diblock copolymer melts when the
long-range interactions are partially screened by fluctuations in the background nuclear fluid
density.

Note that even though Bessel kernels decay exponentially, the situation is very different from the
compact support case: here, there is always a little interaction between pieces of $E$, no matter
how far they are to one another, thus we cannot use the strategy implemented in \autocite{Rig} to
get compactness of minimizing sequences, even in dimension $n=2$.

For small  masses, the similarity between Riesz and Bessel kernels near the origin suggests that
\cref{origpb} presents the same kind of behavior whether $G$ is a Riesz or a Bessel kernel
of the same order $\alpha$, that is, there exists a critical mass below which, up to translations,
the ball of volume $m$ is the unique minimizer.
In this ``small volume'' case, we believe the approach for the Riesz case in \autocite{FFMMM} can be
adapted without major difficulties, but this is not the subject of this paper. We are more
interested in the case of large volumes.
For Riesz kernels of order $\alpha\in(n-2,n)$, it is known that above a critical mass, \cref{origpb}
admits no minimizers. Here, the better integrability of the Bessel kernels changes the asymptotic
behavior when the mass goes to infinity: if $\kappa$ is small enough, \cref{origpb} admits large
mass minimizers, and up to translations, any sequence of normalized (to unit mass) minimizers
converges to the unit ball as the mass goes to infinity (see \cref{mainapp:bessel}). We end this
introductory discussion on Bessel kernels here, leaving the more technical reminders for
\cref{subsec:bessel}.

\section{Preliminaries}\label{sec:prelim}

For the rest of the paper, we shall always assume that the kernel $G$ satisfies assumptions
\crefrange{Krad}{Kint}.

\subsection{Nonlocal perimeters}\label{subsec:reform}

As mentioned in the introduction, the rest of our study relies on the introduction of the functional
$\Per_G$, sometimes referred to as the \textit{nonlocal $G$-perimeter of $E$} (see e.g.
\autocite{BP,CN}) and the reformulation \cref{minrp} of the original problem.
One of the reasons why $\Per_G$ can be thought of as a perimeter appears if one imagines that the
kernel $G$ is singular at the origin, and decreases quickly away from it. Heuristically in that case
the part in $\Per_G(E)$ that prevails would be when $x$ and $y$ are close to one another, so that
\[
\Per_G(E) \simeq \iintb_{\substack{E\times E^\compl\\\abs{x-y}<\eps}} G(x-y)\dx\dy,
\]
for a small positive $\eps$.
But notice that the set $\left(E\times E^{\compl}\right)\cap\{\abs{x-y}<\eps\}$ is included in
$\left(\partial E+B_\eps\right)^2$. Hence what seems to prevail in $\Per_G(E)$ is the interaction
near the boundary $\partial E$; however, in our case the kernel may not be singular at the origin,
nor radially nonincreasing.

In addition, $\Per_G(E)$ can be controlled by the classical perimeter, using that
\[
\norm{f(h+\cdot)-f}_{L^1(\Ren)}\leq \abs{h}\seminorm{f}_{\bv(\Ren)},\quad\forall
f\in\bv(\Ren),~\forall
h\in\Ren.
\]
Indeed, by a change of variables and Fubini's theorem we find
\[
\begin{aligned}
\Per_G(E) = \iint_{E\times E^\compl} G(x-y) \dx\dy
&= \frac{1}{2} \iint_{\Ren\times\Ren} \abs{\ind_E(x+h)-\ind_E(x)} G(h) \dx\dd h\\
&= \frac{1}{2} \int_{\Ren} \norm{\ind_E(h+\cdot)-\ind_E}_{L^1(\Ren)} G(h)\dd h\\
&\leq \frac{1}{2} \int_{\Ren} \abs{h}P(E) G(h)\dd h
= \frac{I_G^1}{2} P(E).
\end{aligned}
\]
We can actually refine the constant in this inequality using the following proposition, inspired by
\autocite[Lemma 3]{Dav}.
\begin{prp}\label{prp:boundK1n}
Let $f \in \bv(\Ren)$, and let $\rho: (0,+\infty)\to [0,+\infty)$ be a measurable function such that
$\displaystyle\int_{\Ren} \rho(\abs{x})\dx = 1$.
Then
\begin{equation}\label{boundK1n:globres}
\iint_{\Ren\times \Ren} \frac{\abs{f(x)-f(y)}}{\abs{x-y}} \rho(\abs{x-y}) \dx\dy \leq
\bfK_{1,n} \int_{\Ren} \abs{Df},
\end{equation}
where $\bfK_{1,n}$ is defined by \cref{eq:defKpn}.
\end{prp}
\begin{proof}
The proof is similar to \autocite[Proof of Lemma 3]{Dav}, but we detail it here for the reader's
convenience.
By approximation (see e.g. \autocite[Theorem~5.3]{EG2015}), we may assume that $f\in
C^\infty(\Ren)$, and by cut-off, that $f$ is also compactly supported.
Integrating on lines we have
\[
f(x+h)-f(x)=\int_0^1 \grad f(x+th)\cdot h\dt,
\]
thus, making the change of variables $h=x-y$, using Fubini's theorem we find
\begin{equation}\label{boundK1n:eq1}
\begin{aligned}
\iint_{\Ren\times \Ren} \frac{\abs{f(x)-f(y)}}{\abs{x-y}} \rho(\abs{x-y})\dx\dy
=\int_{\Ren} \left(\int_{\Ren} \abs*{\int_0^1 \grad f(x-th)\cdot h\dt}\dx\right)
\frac{\rho(\abs{h})}{\abs{h}}\dd h.
\end{aligned}
\end{equation}
Applying the coarea formula to \cref{boundK1n:eq1} and then Cauchy-Schwarz inequality and Fubini's
theorem gives
\begin{equation}\label{boundK1n:eq2}
\begin{aligned}
&\iint_{B_R\times B_R} \frac{\abs{f(x)-f(y)}}{\abs{x-y}} \rho(\abs{x-y})\dx\dy\\
&\hphantom{\iint}
=\int_0^{\infty} \int_{\bbS^{n-1}} \left(\int_{\Ren} \abs*{\int_0^1 \grad
f(x-tr\sigma)\cdot \sigma\dt}\dx\right) \rho(r)r^{n-1}\dH_\sigma^{n-1} \dr\\
&\hphantom{\iint}
\leq\int_0^{\infty} \int_{\bbS^{n-1}} \left(\int_0^1\int_{\Ren} \abs{\grad
f(x-tr\sigma)\cdot\sigma}\dx\dt\right) \rho(r)r^{n-1}\dH_\sigma^{n-1} \dr\\
&\hphantom{\iint}
=\int_0^{\infty} \int_{\bbS^{n-1}} \int_{\Ren} \abs{\grad
f(x)\cdot\sigma}\dx \rho(r)r^{n-1}\dH_\sigma^{n-1} \dr.
\end{aligned}
\end{equation}
Using Fubini's theorem once again and the equality
\[
\dashint_{\bbS^{n-1}} \abs{\grad f(x)\cdot \sigma} \dH^{n-1}_\sigma = \bfK_{1,n}\abs{\grad f(x)},
\]
from \cref{boundK1n:eq2} we obtain
\[
\begin{aligned}
&\iint_{\Ren\times\Ren} \frac{\abs{f(x)-f(y)}}{\abs{x-y}} \rho(\abs{x-y})\dx\dy\\
&\hphantom{\int}
\leq \abs{\bbS^{n-1}}\bfK_{1,n}\int_0^{\infty} \int_{\Ren} \abs{\grad f(x)}\rho(r)r^{n-1}
\dx\dr
= \bfK_{1,n} \left(\int_{\Ren} \abs{\grad f(y)}\dy\right)\left(\int_{\Ren} \rho(\abs{x})\dx\right),
\end{aligned}
\]
hence \cref{boundK1n:globres}, since
\[
\int_{\Ren} \rho(\abs{x}) \dx=1.
\]
\end{proof}

\begin{rmk}
Note that by \cref{prp:limK1n} further below, the constant $\bfK_{1,n}$ in \cref{prp:boundK1n} is
optimal.
\end{rmk}

Setting $\displaystyle\rho_G(r)\coloneqq \frac{rg(r)}{2}$, we have $\displaystyle\int_\Ren
\rho_G(\abs{x})\dx = \frac{1}{\bfK_{1,n}}$ in view of \cref{Kint}, and rewriting $\Per_G(E)$ as
\begin{equation}\label{eq:rewritePerG}
\begin{aligned}
\Per_G(E)
&= \frac{1}{2}\iint_{\Ren\times \Ren} \abs{\ind_E(x)-\ind_E(y)}G(x-y)\dx\dy\\
&= \iint_{\Ren\times \Ren} \frac{\abs{\ind_E(x)-\ind_E(y)}}{\abs{x-y}}
\rho_G(\abs{x-y})\dx\dy,
\end{aligned}
\end{equation}
a direct application of the previous proposition with $f=\ind_E$ leads to the following
control of the nonlocal perimeter by the local perimeter.

\begin{cor}\label{cor:PerKbound}
For any set of finite perimeter $E$, we have
\[
\Per_G(E) \leq P(E).
\]
\end{cor}

Recall that the classical perimeter is lower semicontinuous with respect to the classical topology
of $L^1(\Ren)$. Here the nonlocal perimeter is in fact continuous w.r.t. the $L^1$ convergence, as
is shown in the following lemma.

\begin{lem}\label{lem:perkC0}
For any sets $E$ and $F$ with finite Lebesgue measure, we have
\[
\abs{\Per_G(E)-\Per_G(F)} \leq I_G^{0}\abs{E\triangle F}.
\]
\end{lem}

\begin{proof}
\edit{The argument is classical, and can be found e.g. in \autocite[Lemma~5.2]{FFMMM} or
\autocite[Proof~of~Lemma~3.1]{KM1}, but we present it here for the reader's convenience.}
Let $E$ and $F$ be sets with finite (possibly different) Lebesgue measure. Using
\[
\Per_G(E) = \frac{1}{2} \iint_{\Ren\times\Ren} \abs{\ind_E(x)-\ind_E(y)}G(x-y)\dx\dy,
\]
we have
\[
\Per_G(E)-\Per_G(F)=\frac{1}{2}\iint_{\Ren\times\Ren}
\left(\abs{\ind_E(x)-\ind_E(y)}-\abs{\ind_F(x)-\ind_F(y)}\right)G(x-y)\dx\dy.
\]
Thus by the triangle inequality,
\[
\begin{aligned}
\abs{\Per_G(E)-\Per_G(F)}
&\leq\frac{1}{2}\iint_{\Ren\times\Ren}
\left(\abs{\ind_E(x)-\ind_F(x)}+\abs{\ind_E(y)-\ind_F(y)}\right)G(x-y)\dx\dy\\
&=\frac{1}{2}\iint_{\Ren\times\Ren}
\left(\ind_{E\triangle F}(x)+\ind_{E\triangle F}(y)\right)G(x-y)\dx\dy
=I_G^{0}\abs{E\triangle F}.
\end{aligned}
\]
\end{proof}

One of the nice properties of the classical perimeter is also its behavior under scaling.
For any set of finite perimeter $E$ and any $t> 0$, we obviously have $P(t
E)=t^{n-1}P(E)$. This is unfortunately not the case for such nonlocal perimeters.
However just as $\ddt [P(tE)]=(n-1)t^{n-2}P(E)$, one can show that $\ddt[\Per_G(tE)]$ is
(essentially) bounded by $Ct^{n-2}P(E)$, for some $C=C(n)$, thanks to assumptions
\cref{Krad} and \cref{Kint}.

\begin{lem}\label{lem:perkC1}
For any set of finite perimeter $E$, the function $t\mapsto \Per_G(tE)$ is locally Lipschitz
continuous in $(0,+\infty)$, and for almost every $t$,
\begin{equation}\label{eq:perkderiv}
\begin{aligned}
\frac{\dd}{\dt} \left[\Per_G(tE)\right]
&= \frac{n}{t}\Per_G(tE)+\frac{1}{t}\int_{tE}\left(\int_{\partial^* (tE)}
G(x-y)(x-y)\cdot\nu_{tE}(y)\dH^{n-1}_y\right)\dx.
\end{aligned}
\end{equation}
In particular
\begin{equation}\label{eq:perkscalingdiff}
\Big|\ddt[\Per_G(tE)]\,\Big| \leq Ct^{n-2}P(E),\quad\text{ for a.e. }
t\in(0,+\infty),
\end{equation}
for some $C=C(n)$.
\end{lem}

\begin{proof}
By approximation, there exists a sequence of functions $(G_k)_{k\in\N}\subsq C^\infty_c(\Ren)$
converging to $G$ in $L^1(\Ren)$ and such that $\abs{G_k}\leq G$.
By scaling, we have
\[
\Per_{G_k}(tE) = t^{2n}\iint_{E\times E^\compl} G_k(t(x-y))\dx\dy,
\]
thus, differentiating under the integral and changing variables, we find
\begin{equation}\label{perkderiv:eq1}
\begin{aligned}
\ddt [\Per_{G_k}(tE)]
&= \frac{2n}{t}\Per_{G_k}(tE)
+t^{2n}\iint_{E\times E^\compl} \grad G_k(t(x-y))\cdot(x-y)\dx\dy\\
&= \frac{2n}{t}\Per_{G_k}(tE)
+\frac{1}{t}\iint_{(tE)\times (tE)^\compl} \grad G_k(x-y)\cdot(x-y)\dx\dy.
\end{aligned}
\end{equation}
Notice that
\[
\dvg_y \big((x-y)G_k(x-y)\big)=-nG_k(x-y)-\grad G_k(x-y)\cdot(x-y),
\]
thus by the divergence theorem,
\[
\begin{aligned}
\int_{(tE)^\compl} \grad G_k(x-y)\cdot(x-y)\dy
&=-n\int_{(tE)^\compl} G_k(x-y)\dy\\
&\hphantom{-n\int_{(tE)^\compl}}
+\int_{\partial^*(tE)} G_k(x-y)(x-y)\cdot\nu_{tE}(y)\dH^{n-1}_y.
\end{aligned}
\]
Plugging this into \cref{perkderiv:eq1} gives
\begin{equation}\label{perkderiv:eq2}
\ddt [\Per_{G_k}(tE)]
= \frac{n}{t}\Per_{G_k}(tE)
+\frac{1}{t}\int_{tE}\left(\int_{\partial^*(tE)} G_k(x-y)(x-y)\cdot\nu_{tE}(y)\dH^{n-1}_y\right)\dx.
\end{equation}
One the one hand, notice that
\[
\abs{\Per_G(tE)-\Per_{G_k}(tE)}\leq t^n\abs{E}\norm{G_k-G}_{L^1(\Ren)},
\qquad\,\forall k\in\N,
\]
thus $t\mapsto \Per_{G_k}(tE)$ and $t\mapsto \frac{n}{t}\Per_{G_k}(tE)$ converge respectively
locally uniformly in $(0,\infty)$ to $t\mapsto \Per_{G}(tE)$ and $t\mapsto \frac{n}{t}\Per_{G}(tE)$
as $k$ goes to infinity. On the other hand, given $\eps>0$, since $I_G^1$ is finite and
$\abs{G_k}\leq G$, we may choose $R$ large enough so that
\[
\int_{B_R^\compl} \abs{x}G(x)\dx\leq\eps \quad\text{ and }\quad \int_{B_R^\compl}
\abs{x}\abs{G_k(x)}\dx\leq\eps,
\]
for every $k\in\N$.
Thus, splitting the integral below into the two parts $\abs{x-y}<R$ and $\abs{x-y}\geq R$, we find
\[
\begin{aligned}
&\int_{tE}\left(\int_{\partial^*(tE)}
\big|(G_k-G)(x-y)(x-y)\cdot\nu_{tE}(y)\big|\dH^{n-1}_y\right)\dx\\
&\hphantom{\int_{tE}\left(\int\right.}
\leq P(tE)\left(\int_{B_R} \abs{x}\abs{(G_k-G)(x)}\dx + 2\eps\right)\\
&\hphantom{\int_{tE}\left(\int\right.}
\leq t^{n-1}P(E)\left(R\norm{G_k-G}_{L^1(\Ren)}+2\eps\right),\qquad\forall t>0,
\,\forall k\in\N,
\end{aligned}
\]
which shows local uniform convergence in $(0,\infty)$ of the last term in \cref{perkderiv:eq2}, by
the arbitrarines of $\eps$.
Whence, given any $\vphi\in C^\infty_c(0,\infty)$, passing to the limit in
\[
\int_0^\infty \Per_{G_k}(tE)\vphi'(t)\dt = -\int_0^\infty \ddt [\Per_{G_k}(tE)] \vphi(t)\dt
\]
yields that $t\mapsto \Per_G(tE)$ is weakly differentiable, and for a.e. $t\in(0,\infty)$, we have
\[
\ddt [\Per_{G}(tE)]
= \frac{n}{t}\Per_{G}(tE)
+\frac{1}{t}\int_{tE}\left(\int_{\partial^*(tE)} G(x-y)(x-y)\cdot\nu_{tE}(y)\dH^{n-1}_y\right)\dx.
\]
Then \cref{eq:perkscalingdiff} follows by using $\Per_G(tE)\leq P(tE)$ and noticing that
\[
\abs*{\int_{tE} \left(\int_{\partial^*(tE)} G(x-y)(x-y)\cdot\nu_E(y)\dH^{n-1}_y\right)\dx}
\leq I_G^1P(tE)=\frac{2}{\bfK_{1,n}}P(tE).
\]
In turn, \cref{eq:perkscalingdiff} gives the local Lipschitz continuity of $t\mapsto\Per_G(tE)$.
\end{proof}

Now, although \cref{eq:perkscalingdiff} would be sufficient for our needs, we can use
\cref{lem:perkC1} to get a clearer idea of how much $\calF_{\gamma,G_\lambda}(E)$ may increase or
decrease by scaling $E$.

\begin{lem}\label{lem:fscaling}
Assume $\gamma<1$, and let $E$ be a set with finite perimeter. Then for any $t>1$, we have
\[
\calF_{\gamma,G_\lambda}(tE)-\calF_{\gamma,G_\lambda}(E)\leq (t^{n}-1)\calF_{\gamma,G_\lambda}(E),
\]
and for any $0<t<1$,
\[
\calF_{\gamma,G_\lambda}(tE)-\calF_{\gamma,G_\lambda}(E)\leq
2\gamma t^{n-1}\big(1-t\big)P(E).
\]
\end{lem}

\begin{proof}
\edit{Let $u(t):=\Per_{G_\lambda}(tE)$, and observe that $G_\lambda$ satisfies the same assumptions
as $G$, thus by \cref{lem:perkC1}, $u$ is locally Lipschitz continuous, and
$u'(t)=\frac{n}{t}u(t)+\frac{1}{t}f(t)$, where}
\[
f(t):=\int_{tE}\left(\int_{\partial^* (tE)}
G_\lambda(x-y)(x-y)\cdot\nu_{tE}(y)\dH^{n-1}_y\right)\dx.
\]
Since we want an upper bound on $\calF_{\gamma,G_\lambda}$, we are going to give a lower bound for
$f(t)$.
Dropping the subscript in $\nu_{tE}$, let us write
\[
\begin{aligned}
f(t)
&\geq \iint\limits_{\substack{tE\times\partial^* (tE)\\(x-y)\cdot\nu(y)<0}}
G_\lambda(x-y)(x-y)\cdot\nu(y)\dx\dH_y^{n-1}\\
&\geq -\int_{\partial^*(tE)} \left(\int_{\Ren\,\cap\,\{z\cdot\nu(y)<0\}}
G_\lambda(z)\abs{z\cdot\nu(y)}\dz\right)\dH_y^{n-1}.
\end{aligned}
\]
Using the coarea formula, it follows
\[
\begin{aligned}
f(t)
&\geq -\int_{\partial^* (tE)}\left(\int_0^\infty r^{n}g(r)
\int_{\bbS^{n-1}\,\cap\,\{\sigma\cdot\nu(y)<0\}}
\abs{\sigma\cdot\nu(y)}\dd\sigma\dr\right)\dH_y^{n-1}\\
&= -\frac{\bfK_{1,n}}{2}\int_{\partial^* (tE)}\left(\abs{\bbS^{n-1}}\int_0^\infty r^{n}g(r)
\dr\right)\dH_y^{n-1}\\
&= -\frac{\bfK_{1,n}I_{G}^1}{2}P(tE)=-P(tE),
\end{aligned}
\]
and similarly
\[
f(t)
\leq \int_{\partial^* (tE)}\left(\int_0^\infty r^{n}g(r)
\int_{\bbS^{n-1}\,\cap\,\{\sigma\cdot\nu(y)>0\}}
\abs{\sigma\cdot\nu(y)}\dd\sigma\dr\right)\dH_y^{n-1}
= \frac{\bfK_{1,n}I_{G}^1}{2}P(tE)=P(tE),
\]
thus
\[
\abs*{\ddt[t^{-n}u(t)]}=\abs*{t^{-n}\left(u'(t)-\frac{n}{t}u(t)\right)}=t^{-(n+1)}\abs{f(t)}
\leq \frac{P(E)}{t^2}.
\]
In both cases $t<1$ and $t>1$, integrating \edit{between $1$ and $t$} gives
\[
u(t)\geq t^nu(1)-t^{n-1}\abs{1-t}P(E).
\]
Hence,
\[
\begin{aligned}
\calF_{\gamma,G_\lambda}(tE)-\calF_{\gamma,G_\lambda}(E)
&=(t^{n-1}-1)P(E)-\gamma(u(t)-u(1))\\
&\leq (t^{n-1}-1)P(E)-\gamma (t^{n}-1)u(1)\edit{+\gamma t^{n-1}\abs{1-t}}P(E).\\
\end{aligned}
\]
If $t>1$ we obtain
\[
\begin{aligned}
\calF_{\gamma,G_\lambda}(tE)-\calF_{\gamma,G_\lambda}(E)
&\leq (t^n-1+t^{n-1}-t^n)P(E)-\gamma(t^n-1)u(1)-\gamma t^{n-1}(1-t)P(E)\\
&= (t^n-1)\calF_{\gamma,G_\lambda}(E)-(1-\gamma)t^{n-1}(t-1)P(E)\\
&\leq (t^n-1)\calF_{\gamma,G_\lambda}(E),
\end{aligned}
\]
and, if $t<1$, using that $u(1)=\Per_{G_\lambda}(E)\leq P(E)$ by \cref{cor:PerKbound},
\[
\begin{aligned}
\calF_{\gamma,G_\lambda}(tE)-\calF_{\gamma,G_\lambda}(E)
&\leq (1-\gamma)(t^{n-1}-1)P(E)+\gamma\big((1-t^n)u(1)+t^{n-1}(1-t)P(E)\\
&\hphantom{\leq(1-}
+(t^{n-1}-1)P(E)\big)\\
&\leq (1-\gamma)(t^{n-1}-1)P(E)+2\gamma t^{n-1}\big(1-t\big)P(E)\\
&\leq 2\gamma t^{n-1}\big(1-t\big)P(E)\\
\end{aligned}
\]
where in both cases we used the fact that $\gamma<1$. This concludes the proof.
\end{proof}

Starting from $\rho_G(r)=\frac{rg(r)}{2}$, let us introduce the rescalings
\begin{equation}\label{eq:defrhog}
\rho_{G,\eps}(r):=\eps^{-n}\rho_G(\eps^{-1}r),
\end{equation}
for every $\eps>0$. Then it is convenient to rewrite $\Per_{G_\lambda}(E)$ as
\begin{equation}\label{eq:defperglambda}
\begin{aligned}
\Per_{G_\lambda}(E)
= \iint_{\Ren \times \Ren}
\frac{\abs{\ind_E(x)-\ind_E(y)}}{\abs{x-y}} \rho_{G,1/\lambda}(\abs{x-y})\dd x \dd y.
\end{aligned}
\end{equation}
Notice that, by \cref{Kint}, 
\begin{equation}\label{eq:approxid1}
\int_{\Ren} \rho_{G,1/\lambda}(\abs{x})\dx = \frac{I_G^1}{2} =\frac{1}{\bfK_{1,n}}.
\end{equation}

The nice thing with this expression is that the family of functions $(\rho_{G,\eps})_{\eps>0}$
constitutes, up to multiplication by the universal constant $\bfK_{1,n}$, a \textsl{$n$-dimensional
approximation of identity}, which we define just below.

\begin{dfn}[Approximation of identity]\label{dfn:mollif}
For any $k\in\N\setminus\{0\}$, we say that a family of measurable functions $(\rho_\eps)_{\eps>0}$
from $(0,+\infty)$ into $[0,+\infty)$ is a \textsl{$k$-dimensional approximation of identity} if,
for all $\eps>0$, we have
\begin{enumerate}[label=(\roman*),parsep=0pt,itemsep=0.5pt]
\begin{multicols}{2}
\item $\displaystyle \abs{\bbS^{k-1}}\int_0^\infty \rho_{\eps}(r)r^{k-1}\dr = 1$;
\item $\displaystyle\lim_{\eps \to 0}\, \int_\delta^\infty \rho_{\eps}(r)r^{k-1} \dr= 0,\quad \forall
\delta > 0$.
\end{multicols}
\end{enumerate}
\end{dfn}

The following proposition will be crucial to shed a light on the behavior of \cref{minrp} when
$\lambda$ goes to infinity.

\begin{prp}[\autocite{Dav}]\label{prp:limK1n}
Let $f \in \bv(\R^n)$, and $(\rho_\eps)_{\eps>0}$ be a $n$-dimensional approximation of identity.
Then we have
\[
\iint_{\Ren\times\Ren} \frac{\abs{f(x)-f(y)}}{\abs{x-y}} \rho_\eps(\abs{x-y}) \dd x\dd y
~\xrightarrow{\eps\to 0}~ \bfK_{1,n} \int_{\Ren} \abs{Df},
\]
where $\bfK_{1,n}$ is defined by \cref{eq:defKpn}.
\end{prp}

Taking $f=\ind_E$ in \cref{prp:limK1n} gives the following immediate corollary, in view
of~\cref{eq:defperglambda,eq:approxid1}.

\begin{cor}\label{cor:cvcalV}
For any set of finite perimeter $E$ in $\Ren$, we have
\[
\Per_{G_\lambda}(E)~\xrightarrow{\lambda\to\infty}~ P(E).
\]
\end{cor}

This implies that the functional $\calF_{\gamma,G_\lambda}$ converges pointwise to $(1-\gamma)P$
when $\lambda$ goes to infinity. Thus we may guess that if $\gamma<1$, \cref{minrp} will reduce to
minimizing the classical perimeter under the volume constraint $\abs{B_1}$ when $\lambda$ is large.

\subsection{Bessel kernels}\label{subsec:bessel}

For any $\kappa$ and $\alpha$ in $(0,+\infty)$, the Bessel kernel $\mathcal{B}_{\kappa,\alpha}$ is
the fundamental solution $f$ of
\[
(I-\kappa\Delta)^{\frac{\alpha}{2}}\,f = \mathbf{\delta}_0\quad\text{ in }\calD'(\Ren).
\]
When $\kappa=1$, we write $\mathcal{B}_\alpha\coloneqq\mathcal{B}_{1,\alpha}$.
The following proposition sums up some basic properties of $\mathcal{B}_\alpha$ (see e.g.
\autocite[Chapter I.2.2]{Graf}, \autocite[Chapter V.3]{Stein} and \autocite[Chapter II.3]{AroSmi}).

\begin{prp}\label{prp:propbessel}
The Bessel kernel of order $\alpha\in (0,+\infty)$ in $\Ren$ is given by
\begin{equation}\label{eq:besselform}
\mathcal{B}_\alpha(x) =
\frac{1}{(4\pi)^{\frac{\alpha}{2}}}\frac{1}{\Gamma\left(\frac{\alpha}{2}\right)}\int_0^\infty
e^{-\frac{\pi\abs{x}^2}{t}}e^{-\frac{t}{4\pi}} t^{\frac{\alpha-n}{2}} \frac{\dt}{t}.
\end{equation}
The kernel $\mathcal{B}_\alpha$ is radial \edit{and $C^\infty$ away from the origin}. In
addition
\[
I_{\mathcal{B}_\alpha}^{0}=1,\quad\text{ and }\quad
I_{\mathcal{B}_\alpha}^{1}=
n\frac{\Gamma\left(\frac{1+\alpha}{2}\right)}{\Gamma\left(\frac{\alpha}{2}\right)}
\frac{\Gamma\left(\frac{1+n}{2}\right)}{\Gamma\left(1+\frac{n}{2}\right)}.
\]
The asymptotic behavior of $\mathcal{B}_\alpha$ is
\[
\mathcal{B}_\alpha(x) \underset{0}{\sim}
\begin{dcases} 
\frac{\Gamma\left(\frac{n-\alpha}{2}\right)}{2^\alpha\pi^{\frac{n}{2}}
\Gamma\left(\frac{\alpha}{2}\right)} \frac{1}{\abs{x}^{n-\alpha}}
&
\text{ if }\/ 0<\alpha<n\\
\frac{-\log(\abs{x})}{2^{n-1}\pi^{\frac{n}{2}}\Gamma\left(\frac{n}{2}\right)}&\text{ if }\/
\alpha=n\\
\frac{\Gamma\left(\frac{\alpha-n}{2}\right)}{2^n\pi^{\frac{n}{2}}}&\text{ if }\/ n<\alpha,
\end{dcases}
\]
and
\[
\mathcal{B}_\alpha(x) \underset{\infty}{\sim}
\frac{1}{2^{\frac{n+\alpha-1}{2}}\pi^{\frac{n-1}{2}}\Gamma\left(\frac{\alpha}{2}\right)}
\abs{x}^{\frac{\alpha-n-1}{2}}e^{-\abs{x}}.
\]
By scaling, for any $\alpha,\kappa>0$, the (generalized) Bessel kernel $\mathcal{B}_{\kappa,\alpha}$
is given by
\begin{equation}\label{propbessel:scaling}
\mathcal{B}_{\kappa,\alpha}(x) = \frac{1}{\kappa^{\frac{n}{2}}}
\mathcal{B}_\alpha\left(\frac{x}{\sqrt{\kappa}}\right),\quad\forall x\in\Ren\setminus\{0\},
\end{equation}
thus
\[
I_{\mathcal{B}_{\kappa,\alpha}}^{0} = 1,\quad\text{ and }\quad
I_{\mathcal{B}_{\kappa,\alpha}}^{1}=\kappa^{\frac{1}{2}}I_{\mathcal{B}_\alpha}^{0,1}.
\]
\end{prp}
\begin{proof}
The integral representation \cref{eq:besselform} and the asymptotics can be found respectively in
\autocite{Graf} and \autocite{AroSmi}, and \edit{the facts that $\mathcal{B}_{\alpha}\in
C^\infty(\Ren\setminus\{0\})$ and that
$I_{\mathcal{B}_\alpha}^{0}=\norm{\mathcal{B}_\alpha}_{L^1(\Ren)}=1$ are well-known}, so we detail
only the computations of $I_{\mathcal{B}_\alpha}^{1}$.
By \cref{eq:besselform}, using Fubini's theorem, we find
\begin{equation}\label{propbessel:eq1}
I_{\mathcal{B}_\alpha}^{1}=\int_\Ren \abs{x}\mathcal{B}_\alpha(x)\dx
= \frac{1}{(4\pi)^{\frac{\alpha}{2}}}\frac{1}{\Gamma\left(\frac{\alpha}{2}\right)}\int_0^\infty
e^{-\frac{t}{4\pi}}t^{\frac{\alpha-n}{2}}\left(\int_{\Ren}
\abs{x}e^{-\frac{\pi\abs{x}^2}{t}}\dx\right)\frac{\dt}{t}.
\end{equation}
Changing variables, we compute
\begin{equation}\label{propbessel:eq2}
\begin{aligned}
\int_{\Ren}
\abs{x}e^{-\frac{\pi\abs{x}^2}{t}}\dx
&=n\om_n\left(\frac{t}{\pi}\right)^{\frac{n+1}{2}}
\int_0^\infty r^n e^{-r^2}\dr\\
&=n\om_n\left(\frac{t}{\pi}\right)^{\frac{n+1}{2}}
\frac{\Gamma\left(\frac{n+1}{2}\right)}{2}
=\frac{nt^{\frac{n+1}{2}}}{\sqrt{4\pi}}\frac{\Gamma\left(\frac{n+1}{2}\right)}
{\Gamma\left(1+\frac{n}{2}\right)}.
\end{aligned}
\end{equation}
Injecting \cref{propbessel:eq2} into \cref{propbessel:eq1} yields
\[
\begin{aligned}
I_{\mathcal{B}_\alpha}^{1}
&= \frac{1}{(4\pi)^{\frac{1+\alpha}{2}}}\frac{1}{\Gamma\left(\frac{\alpha}{2}\right)}
\frac{n\Gamma\left(\frac{n+1}{2}\right)}{\Gamma\left(1+\frac{n}{2}\right)}
\int_0^\infty  e^{-\frac{t}{4\pi}} t^{\frac{\alpha+1}{2}} \frac{\dt}{t}\\
&= \frac{1}{(4\pi)^{\frac{1+\alpha}{2}}}\frac{1}{\Gamma\left(\frac{\alpha}{2}\right)}
\frac{n\Gamma\left(\frac{n+1}{2}\right)}{\Gamma\left(1+\frac{n}{2}\right)}
(4\pi)^{\frac{1+\alpha}{2}}\Gamma\left(\frac{1+\alpha}{2}\right)
=
n\frac{\Gamma\left(\frac{1+\alpha}{2}\right)}{\Gamma\left(\frac{\alpha}{2}\right)}
\frac{\Gamma\left(\frac{1+n}{2}\right)}{\Gamma\left(1+\frac{n}{2}\right)}.
\end{aligned}
\]
\end{proof}

\begin{rmk}\label{rmk:explbessel}
In some cases, the Bessel kernels $\mathcal{B}_\alpha$ take a very simple form. Indeed, they can
also be expressed in terms of the modified Bessel functions of the third kind
$\mathrm{K}_\nu:(0,+\infty)\to(0,+\infty)$, defined for any $\nu\in\R$ by
\begin{equation}\label{eq:3rdbessel}
\mathrm{K}_\nu(r)=\left(\frac{\pi}{2}\right)^{\frac{1}{2}} \frac{r^\nu
e^{-r}}{\Gamma\left(\nu+\frac{1}{2}\right)}\int_0^\infty
e^{-rt}\left(t+\frac{t^2}{2}\right)^{\nu-\frac{1}{2}}\dt\qquad\text{ if }\nu>-\frac{1}{2},
\end{equation}
and the relation $\mathrm{K}_\nu=\mathrm{K}_{-\nu}$ (see \autocite[Chapter~II.3]{AroSmi}). Then by
\autocite[Chapter~II.4]{AroSmi}, $\mathcal{B}_\alpha$ is given by
\begin{equation}\label{eq:besselaltform}
\mathcal{B}_\alpha(x) =
\frac{1}{2^{\frac{n+\alpha-2}{2}}\pi^{\frac{n}{2}}\Gamma\left(\frac{\alpha}{2}\right)}
\frac{\mathrm{K}_{\frac{n-\alpha}{2}}(\abs{x})}{\abs{x}^{\frac{n-\alpha}{2}}},\quad\forall x\in
\Ren\setminus\{0\}.
\end{equation}
From \cref{eq:3rdbessel,eq:besselaltform} it is easy to see that when $\alpha=n-1$,
$\mathcal{B}_\alpha$ takes the explicit form
\[
\mathcal{B}_{n-1}= \frac{1}{(4\pi)^{\frac{n-1}{2}}\Gamma\left(\frac{n-1}{2}\right)}
\frac{e^{-\abs{x}}}{\abs{x}}.
\]
In particular, when $n=3$ and $\alpha=2$,
$\mathcal{B}_\alpha(x)=\frac{1}{4\pi}\frac{e^{-\abs{x}}}{\abs{x}}$.
When $\alpha=n+1$, changing variables in \cref{eq:besselform}, one can compute $\mathcal{B}_\alpha$
explicitly as well. Indeed, in that case,
\[
\mathcal{B}_{n+1}(x)
= 
\frac{1}{(4\pi)^{\frac{n+1}{2}}}\frac{1}{\Gamma\left(\frac{n+1}{2}\right)}
2\sqrt{4\pi}\int_0^\infty e^{-t^2-\frac{\abs{x}^2}{4t^2}}\dt
= \frac{2\pi}{(4\pi)^{\frac{n+1}{2}}}\frac{1}{\Gamma\left(\frac{n+1}{2}\right)}
e^{-\abs{x}}.
\]
\end{rmk}

We have the following straightforward corollary of \cref{prp:propbessel}.

\begin{cor}\label{cor:assumbessel}
For every $\alpha,\kappa\in(0,+\infty)$, the kernel $\mathcal{B}_{\kappa,\alpha}$
satisfies assumptions \crefrange{Krad}{KC1}.
\end{cor}

We can express the constants $\bfK_{p,n}$ in terms of the Gamma function as follows, in order to
compute explicitly $\gamma$ such that $\gamma
I_{\mathcal{B}_{\kappa,\alpha}}^{1}=\frac{2}{\bfK_{1,n}}$ and apply
\cref{mainthm:existmin,thm:cvminim,thm:stabthreshold}.

\begin{lem}\label{lem:expKpn}
For any $n\in\N$ and $p>0$, we have
\[
\bfK_{p,n} = \frac{\Gamma\left(\frac{n}{2}\right)\Gamma\left(\frac{1+p}{2}\right)}
{\sqrt{\pi}\,\Gamma\left(\frac{n+p}{2}\right)},
\]
where $\bfK_{p,n}$ is given by \cref{eq:defKpn}.
In particular
\[
\bfK_{1,n} = \frac{\Gamma\left(\frac{n}{2}\right)}{\sqrt{\pi}\,\Gamma\left(\frac{n+1}{2}\right)}
\quad\text{ and }\quad
\bfK_{2,n} = \frac{1}{n}.
\]
\end{lem}

\begin{proof}
In the definition of $\bfK_{p,n}$, we may assume $e=(0,\dotsc,0,1)$ without loss of generality, so
that $e\cdot x=x_n$.
Recall that for every nonnegative $\calH^{n-1}$-measurable function $f$ on $\bbS^{n-1}$, we have
(see e.g. \autocite[Section D.2]{Gra2014} or \autocite[Corollary A.6]{ABW})
\[
\int_{\bbS^{n-1}} f\dH^{n-1} = \int_{-1}^1 (1-t^2)^{\frac{n-3}{2}} \int_{\bbS^{n-2}}
f(\sqrt{1-t^2}x,t) \dH^{n-2}_x\dt.
\]
This way we compute
\[
\int_{\bbS^{n-1}} \abs{x_n}^p \dH^{n-1} = \abs{\bbS^{n-2}}\int_{-1}^1
\abs{t}^p(1-t^2)^{\frac{n-3}{2}}\dt=\abs{\bbS^{n-2}}
\frac{\Gamma\left(\frac{n-1}{2}\right)\Gamma\left(\frac{p+1}{2}\right)}
{\Gamma\left(\frac{n+p}{2}\right)},
\]
thus 
\[
\bfK_{p,n} = \frac{\abs{\bbS^{n-2}}}{\abs{\bbS^{n-1}}}
\frac{\Gamma\left(\frac{n-1}{2}\right)\Gamma\left(\frac{p+1}{2}\right)}
{\Gamma\left(\frac{n+p}{2}\right)}=\frac{\Gamma\left(\frac{n}{2}\right)
\Gamma\left(\frac{1+p}{2}\right)}{\sqrt{\pi}\Gamma\left(\frac{n+p}{2}\right)}.
\]
\end{proof}

In view of of the computations of this section, \cref{mainapp:bessel} is then a direct application
of \crefrange{mainthm:existmin}{thm:stabthreshold}.

\section{Existence and study of large mass minimizers}\label{sec:exist}

\subsection[Existence and convergence as \texorpdfstring{$m\to\infty$}{m goes to infinity}]{%
\for{toc}{Existence and convergence as \texorpdfstring{$m\to\infty$}{m goes to
infinity}}\except{toc}{Existence and convergence as \texorpdfstring{\forcebold{$m\to\infty$}}{m
goes to infinity}}%
}


In order to prove existence of minimizers for large masses, we want to use the direct method in the
calculus of variations, starting from a minimizing sequence. When $\gamma<1$, we will see that any
minimizing sequence is bounded in $\bv(\Ren)$, but in order to get compactness in $L^1(\Ren)$ and
pass to the limit, we need to show that no mass escapes at infinity. To do so, we will need to
establish a few lemmas.
First we show that for large masses, if the energy $\calF_{\gamma,G_\lambda}$ of some set $E$,
(where $\calF_{\gamma,G_\lambda}$ is defined by \cref{deffuncF}) is smaller than that of $B_1$,
then $E$ is actually close to~$B_1$.

\begin{lem}\label{lem:asym_mass}
Assume that $\gamma<1$. Then for any set $E$ of finite perimeter with volume $\abs{B_1}$ such that
\[
\calF_{\gamma,G_\lambda}(E)\leq \calF_{\gamma,G_\lambda}(B_1),
\]
we have
\[
\abs{E\triangle B_1(x)} \leq C_0\eta(\lambda),
\]
for some $x\in\Ren$, where
\begin{equation}\label{eq:defeta}
\eta(\lambda) \coloneqq \Big((P-\Per_{G_\lambda})(B_1)\Big)^{\frac{1}{2}},
\end{equation}
and $C_0=C_0(n,\gamma)$. \edit{Here $B_1(x)$ is the ball achieving the minimum in the definition of
the Fraenkel asymmetry.}
\end{lem}

\begin{rmk}
Notice that $\eta(\lambda)$ vanishes as $\lambda$ goes to infinity by \cref{cor:cvcalV}.
\end{rmk}

\begin{proof}[Proof of {\cref{lem:asym_mass}}]
The inequality
\[
\calF_{\gamma,G_\lambda}(E)\leq \calF_{\gamma,G_\lambda}(B_1)
\]
rewrites
\[
P(E)-P(B_1) \leq \gamma\big(\Per_{G_\lambda}(E)-\Per_{G_\lambda}(B_1)\big).
\]
By \cref{cor:PerKbound}, this implies
\[
\begin{aligned}
P(E)-P(B_1)
&\leq \gamma\big(P(E)-\Per_{G_\lambda}(B_1)\big)\\
&= \gamma\big(P(E)-P(B_1)\big)+\gamma\big(P(B_1)-\Per_{G_\lambda}(B_1)\big),
\end{aligned}
\]
thus
\begin{equation}\label{asym_mass:eq1}
P(E)-P(B_1) \leq \frac{\gamma}{1-\gamma}\Big((P-\Per_{G_\lambda})(B_1)\Big).
\end{equation}
By the quantitative isoperimetric inequality \cref{eq:quantisop}, we have $\alpha(E) \leq
C\sqrt{P(E)-P(B_1)}$, \edit{where $C=C(n)$, which concludes the proof, by definition of the Fraenkel
	asymmetry.}
\end{proof}

We also need a truncation lemma akin to \autocite[Lemma 29.12]{Maggi} or \autocite[Lemma 4.5]{FFMMM}
to quantify by how much the energy decreases after properly cutting a set which is already close to
a ball.

\begin{lem}[Truncation]\label{lem:truncperim}
Assume $\gamma < 1$. There exist $C_1$, $C_2 \in (0,+\infty)$ depending only on
$n$ and $\gamma$ such that the following holds. If $E$ is a set of finite perimeter satisfying
$\abs{E\setminus B_{r_0}} \leq \eta$, for some positive constants $\eta$ and $r_0$, then there
exists $r \in [r_0, r_0+C_1\eta^{\frac{1}{n}}]$ such that
\begin{equation}\label{eq:truncperimccl}
\calF_{\gamma,G_\lambda}(E\cap B_r) \leq \calF_{\gamma,G_\lambda}(E)-\frac{\abs{E\setminus
B_r}^{1-\frac{1}{n}}}{C_2}.
\end{equation}
\end{lem}

\begin{proof}
\edit{Let $C_2\coloneqq \frac{2}{n\om_n^{1/n}(1-\gamma)}$, and $C_1\coloneqq
4\left[\om_n^{\frac{1}{n}}(1-\gamma)-\frac{1}{nC_2}\right]^{-1}>0$, and let $E$ be a set
of finite perimeter such that $\abs{E\setminus B_{r_0}} \leq \eta$.
We define $u(r) \coloneqq \abs{E\setminus B_r}$. If $u(r_0+C_1\eta^{\frac{1}{n}})=0$, then
\cref{eq:truncperimccl} trivially holds with $r=r_0+C_1\eta^{\frac{1}{n}}$, thus from now on we
assume that $u(r_0+C_1\eta^{\frac{1}{n}})>0$.}
Since $u$ is nonincreasing, we have $u(r) > 0$, for all $r
\in [r_0,r_0+C_1\eta^{\frac{1}{n}}]$. Notice that $u$ is absolutely continuous, and $u'(r) =
-\calH^{n-1}(E\cap\partial B_r)$ for $\calL^1$-almost every $r\in [r_0,r_0+C_1\eta^{\frac{1}{n}}]$.
By contradiction, let us assume that
\begin{equation}\label{eq:truncperim1}
P(E)-\gamma \Per_{G_\lambda}(E) < P(E\cap B_r)-\gamma\Per_{G_\lambda}(E\cap
B_r)+\frac{\abs{E\setminus B_r}^{1-\frac{1}{n}}}{C_2},
\end{equation}
for all $r \in [r_0,r_0+C_1\eta^{\frac{1}{n}}]$. Recall that for almost every $r\in [r_0,
r_0+C_1\eta^{\frac{1}{n}}]$ we have $\calH^{n-1}(\partial^* E\cap \partial B_r)=0$ (see e.g.
\autocite[Proposition 2.16]{Maggi}).
Given such an $r$, note that $P(E)=P(E;B_r)+P(E;\ol{B_r}^\compl)$, and $P(E\cap
B_r)=P(E;B_r)+\calH^{n-1}(E\cap\partial B_r)$ (see e.g. \autocite[Lemma 15.12]{Maggi}). Thus
\begin{equation}\label{eq:truncperim2}
\begin{aligned}
P(E)-P(E\cap B_r)
&= P(E; \ol{B_r}^\compl)-\calH^{n-1}(E\cap\partial B_r)\\
&= P(E\setminus B_r)-2\calH^{n-1}(E\cap\partial  B_r)
= P(E\setminus B_r)+2u'(r),
\end{aligned}
\end{equation}
where we also used $P(E\setminus B_r)=P(E; \ol{B_r}^\compl)+\calH^{n-1}(E\cap\partial B_r)$.
On the other hand, noticing that
\[
\begin{aligned}
\Per_{\edit{G_\lambda}}(E)
&= \iint_{\left(E\setminus B_r\right)\times E^\compl} \edit{G_\lambda}(x-y)\dx\dy+\iint_{\left(E\cap
B_r\right)\times E^\compl} \edit{G_\lambda}(x-y)\dx\dy\\
&= \iint_{\left(E\setminus B_r\right)\times (E\setminus B_r)^\compl} \edit{G_\lambda}(x-y)\dx\dy
-\iint_{\left(E\setminus B_r\right)\times\left(E\cap B_r\right)} \edit{G_\lambda}(x-y)\dx\dy\\
&\phantom{-\iint_{\left(E\setminus B_r\right)}}
+\iint_{\left(E\cap B_r\right)\times E^\compl} \edit{G_\lambda}(x-y)\dx\dy\\
&=\Per_{\edit{G_\lambda}}(E\setminus B_r)
-\iint_{\left(E\setminus B_r\right)\times\left(E\cap B_r\right)} \edit{G_\lambda}(x-y)\dx\dy
+\iint_{\left(E\cap B_r\right)\times E^\compl} \edit{G_\lambda}(x-y)\dx\dy
\end{aligned}
\]
and
\[
\Per_{\edit{G_\lambda}}(E\cap B_r)
=\iint_{\left(E\cap B_r\right)\times E^\compl} \edit{G_\lambda}(x-y)\dx\dy+\iint_{\left(E\cap B_r\right)\times
\left(E\setminus B_r\right)} \edit{G_\lambda}(x-y)\dx\dy,
\]
we find
\begin{equation}\label{eq:truncperim3}
\begin{aligned}
\Per_{\edit{G_\lambda}}(E)-\Per_{\edit{G_\lambda}}(E\cap B_r)
&= \Per_{\edit{G_\lambda}}(E\setminus B_r)-2 \iint_{\left(E\cap B_r\right) \times \left(E\cap B_r^\compl\right)}
\edit{G_\lambda}(x-y) \dx\dy.
\end{aligned}
\end{equation}
Injecting \cref{eq:truncperim2} and \cref{eq:truncperim3} into \cref{eq:truncperim1}, one gets
\begin{equation}\label{eq:truncperim4}
\begin{aligned}
P(E\setminus B_r)-\gamma\Per_{\edit{G_\lambda}}(E\setminus B_r)
&<-2u'(r)-2\gamma\iint_{\left(E\cap B_r\right)\times \left(E\cap B_r^\compl\right)} \edit{G_\lambda}(x-y)
\dx\dy+\frac{u(r)^{1-\frac{1}{n}}}{C_2}\\
&\leq-2u'(r)+\frac{u(r)^{1-\frac{1}{n}}}{C_2},
\end{aligned}
\end{equation}
for almost every $r\in[r_0,r_0+C_1\eta^{\frac{1}{n}}]$.
\edit{By \cref{cor:PerKbound}}, we know that $\Per_{\edit{G_\lambda}}(E\setminus B_r)\leq
P(E\setminus B_r)$, \edit{since $G_\lambda$ satisfies the same assumptions as $G$}, thus
\begin{equation}\label{eq:truncperim5}
\begin{aligned}
\left(1-\gamma\right)P(E\setminus B_r) <
-2u'(r)+\frac{u(r)^{1-\frac{1}{n}}}{C_2}\qquad\text{for a.e. }r\in[r_0,r_0+C_1\eta^{\frac{1}{n}}].
\end{aligned}
\end{equation}
Now by the isoperimetric inequality \cref{eq:quantisop} we have
\begin{equation}\label{eq:truncperim5b}
P(E\setminus  B_r) \geq n\om_n^{\frac{1}{n}}u(r)^{1-\frac{1}{n}},
\end{equation}
hence from \cref{eq:truncperim5}, it follows
\[
\left[n\om_n^{\frac{1}{n}}(1-\gamma)-\frac{1}{C_2}\right]
u(r)^{1-\frac{1}{n}}<-2u'(r).
\]
\edit{By our choice of $C_1$ and $C_2$, this gives
\begin{equation}\label{eq:truncperim6}
\frac{2n}{C_1}u(r)^{\frac{n-1}{n}} < -u'(r).
\end{equation}
Then \cref{eq:truncperim6} can be rewritten
\[
\left(u(r)^{\frac{1}{n}}\right)' = \frac{1}{n}u'(r)u(r)^{\frac{1}{n}-1} < -\frac{2}{C_1}\qquad\text{
	for a.e. }r\in[r_0,r_0+C_1\eta^{\frac{1}{n}}],
\]
so that integrating between $r_0$ and $r_0+C_1\eta^{\frac{1}{n}}$ and using the fact that
$u(r_0)\leq\eta$ yields
\[
u(r_0+C_1\eta^{\frac{1}{n}})^{\frac{1}{n}} \leq u(r_0)^{\frac{1}{n}}-2\eta^{\frac{1}{n}} \leq
-\eta^{\frac{1}{n}}<0,
\]
a contradiction. Hence \cref{eq:truncperimccl} holds.}
\end{proof}

We can also prove the following counterpart of the truncation lemma when taking the union with a
ball.

\begin{lem}[Union with a ball]\label{lem:unionperim}
Assume $\gamma < 1$. There exist $C_1$, $C_2 \in (0,+\infty)$ depending only on
$n$ and $\gamma$ such that the following holds. If $E$ is a set of finite perimeter satisfying
$\abs{B_{r_0}\setminus E} \leq \eta$, for some positive constants $\eta$ and $r_0$, then there
exists $r \in [r_0-C_1\eta^{\frac{1}{n}},r_0]$ such that
\[
\calF_{\gamma,G_\lambda}(E\cup B_r) \leq \calF_{\gamma,G_\lambda}(E)-\frac{\abs{B_r\setminus
E}^{1-\frac{1}{n}}}{C_2}.
\]
\end{lem}

\begin{proof}
The proof is almost identical to the one of the truncation lemma, by considering the quantity
$u(r):=\abs{B_r\setminus E}$, and so is left to the reader.
\end{proof}

We are now able to prove the existence of large mass minimizers.
We would like to kindly thank \textsc{F. Générau} for pointing out that the proof of
\cref{thm:existmin} actually shows that \edit{the boundary of} any minimizer $E$ lies between two
balls close to the unit ball.

\begin{thm}\label{thm:existmin}
Assume $\gamma<1$. Then there exists $\lambda_e=\lambda_e(n,\gamma,G)$ such that, for any $\lambda >
\lambda_e$, \cref{minrp} admits a minimizer. In addition, up to a translation and up to a Lebesgue
negligible set, any minimizer of \cref{minrp} with $\lambda>\lambda_e$ satisfies
\begin{equation}\label{existmin:incl}
\ol{B}_{1-C\eta(\lambda)^{\frac{1}{n}}} \subsq E \subsq B_{1+C\eta(\lambda)^{\frac{1}{n}}},
\end{equation}
where $C=C(n,\gamma)$ and $\eta(\lambda)$ is given by \cref{eq:defeta} and vanishes as $\lambda$
goes to infinity by \cref{cor:cvcalV}.
\end{thm}

\begin{proof} In this proof, $C$ denotes a constant depending only on $n$ and $\gamma$, possibly
changing from line to line. We proceed in three steps.\\
\textit{Step 1.} Let us show that there exists $\lambda_e>0$ depending only on $n$, $\gamma$ and
$G$ such that the following holds.
For any set of finite perimeter $E$ of mass $\abs{B_1}$ satisfying $\calF_{\gamma,G_\lambda}(E)
\leq\calF_{\gamma,G_\lambda}(B_1)$,
up to a translation, there exists a set of finite perimeter $\wt{E}$ of same mass such that
\begin{equation}\label{eq:exist00}
\calF_{\gamma,G_\lambda}(\wt{E}) \leq \calF_{\gamma,G_\lambda}(E)-\abs{E\setminus
B_{1+C\eta(\lambda)^{\frac{1}{n}}}}\quad\text{ and }
\quad \wt{E}\subsq B_{1+C\eta(\lambda)^{\frac{1}{n}}}.
\end{equation}
Let $\lambda_e>0$ to be fixed later.
By \cref{lem:asym_mass}, if $E$ satisfies
$\calF_{\gamma,G_\lambda}(E)\leq\calF_{\gamma,G_\lambda}(B_1)$, then up to translating $E$, we have
\[
\abs{E\triangle B_1} \leq C_0\eta(\lambda),
\]
where $C_0=C_0(n,\gamma)$. Using Lemma~\ref{lem:truncperim} with $C_0 \eta(\lambda)$ in
place of $\eta$, we can find
$r\leq 1+C_1C_0^{\frac{1}{n}}\eta(\lambda)^{\frac{1}{n}}$ such that
\begin{equation}\label{eq:exist1}
\calF_{\gamma,G_\lambda}(E\cap B_r)
\leq\calF_{\gamma,G_\lambda}(E)-\frac{u^{1-\frac{1}{n}}}{C_2}
\end{equation}
where $C_2$ depends only on $n$ and $\gamma$, and where we have set $u \coloneqq \abs{E\setminus B_r}$.
If $u=0$, then $E\subsq B_r\subsq B_{1+C_1C_0^{\frac{1}{n}}\eta^{\frac{1}{n}}}$ (up to a Lebesgue
negligible set) and \cref{eq:exist00} holds with $C=C_1C_0^{\frac{1}{n}}$ and $\wt{E}=E$, so we may
now assume that $u>0$.
To compensate for the loss of some mass, we define the rescaled set $\wt{E}\coloneqq\mu(E\cap B_r)$,
where \edit{$\mu \coloneqq \left(1-\frac{u}{\abs{B_1}}\right)^{-\frac{1}{n}}$}, so that
$\abs{\widetilde{E}}=\abs{B_1}$. Then since $\mu\geq 1$, by \cref{lem:fscaling}, we have
\begin{equation}\label{eq:exist2}
\calF_{\gamma,G_\lambda}(\wt{E})-\calF_{\gamma,G_\lambda}(E\cap B_r)
\leq (\mu^n-1)\calF_{\gamma,G_\lambda}(E\cap B_r)
\leq u\calF_{\gamma,G_\lambda}(E\cap B_r),
\end{equation}
provided $\lambda_e$ is large enough (depending only on $n$ and $G$), since $u\leq C_0\eta(\lambda)$
and thus vanishes as $\lambda$ goes to infinity.
Injecting \cref{eq:exist1} into \cref{eq:exist2} and using that $u>0$, one gets
\begin{equation}\label{eq:exist3}
\begin{aligned}
\calF_{\gamma,G_\lambda}(\wt{E})
&\leq (1+u)\left(\calF_{\gamma,G_\lambda}(E)-\frac{u^{1-\frac{1}{n}}}{C_2}\right)\\
&\leq \calF_{\gamma,G_\lambda}(E)+u\left(\calF_{\gamma,G_\lambda}(B_1)
-\frac{1}{C_2u^{\frac{1}{n}}}\right)
\leq \calF_{\gamma,G_\lambda}(E)+u\left(P(B_1)-\frac{1}{C_2u^{\frac{1}{n}}}\right).
\end{aligned}
\end{equation}
Since $u\leq C_0\eta(\lambda)$ vanishes uniformly as $\lambda$ goes to infinity, this gives
\[
\calF_{\gamma,G_\lambda}(\wt{E})\leq \calF_{\gamma,G_\lambda}(E)-u 
\leq \calF_{\gamma,G_\lambda}(E)-\abs{E\setminus B_{1+C\eta(\lambda)^{\frac{1}{n}}}},
\]
provided $\lambda_e$ is chosen large enough depending only on $n$, $\gamma$ and $G$.
Recall that $\wt{E} \subsq B_{\mu r}$, and notice that,
\[
\mu r=\edit{r\left(1-\frac{u}{\abs{B_1}}\right)^{-\frac{1}{n}}\leq r(1+Cu)}
\leq \big(1+C\eta(\lambda)^{\frac{1}{n}}\big)\left(1+C_0\eta(\lambda)\right)
\leq 1+C\eta(\lambda)^{\frac{1}{n}},
\]
whenever $\lambda_e$ is large enough, thus $\wt{E}\subsq B_{1+C\eta(\lambda)^{\frac{1}{n}}}$, which
concludes this step.\\
\textit{Step 2.}
We prove existence of minimizers.
For $\lambda\geq\lambda_e$, consider a minimizing sequence $(E_k)_{k\in\N}$ for \cref{minrp}.
There are two cases: either $B_1$ is a minimizer of $\calF_{\gamma,G_\lambda}$, and we are done, or
$B_1$ is not a minimizer of $\calF_{\gamma,G_\lambda}$, and up to a subsequence (not relabeled), for
all $k\in\N$, we have $\calF_{\gamma,G_\lambda}(E_k)\leq \calF_{\gamma,G_\lambda}(B_1)$. In the
latter case, by Step 1 we can build another minimizing sequence $(\wt{E}_k)_{k\in\N}$ of sets
included in $B_{1+C\eta(\lambda)^{\frac{1}{n}}}$ such that $\calF_{\gamma,G_\lambda}(\wt{E}_k) \leq
\calF_{\gamma,G_\lambda}(E_k)$, for all $k\in\N$.  Now $(\ind_{\wt{E}_k})_{k\in\N}$ is also bounded
in $\bv(\Ren)$. Indeed $\calF_{\gamma,G_\lambda}(\wt{E}_k)\leq \calF_{\gamma,G_\lambda}(B_1)$
implies
\[
P(\wt{E}_k) \leq (1-\gamma)^{-1}P(B_1),
\]
thus $\seminorm{\ind_{\wt{E}_k}}_{\bv(\Ren)}$ is bounded, and
$\norm{\ind_{\wt{E}_k}}_{L^1(\Ren)}=\abs{B_1}$. By \edit{the compactness of the immersion of}
$\bv(\Ren)$ in $L^1_{\loc}(\Ren)$ and the fact that $\wt{E}_k \subsq
B_{1+C\eta(\lambda)^{\frac{1}{n}}}$ after translation, up to the extraction of a subsequence (still
not relabeled), $\ind_{\wt{E}_k}$ converges to some function $f\in \bv(\Ren)$ in $L^1(\Ren)$ and
almost everywhere. The almost everywhere convergence implies that $f$ is the indicator function of
some set of finite perimeter $E$, and the $L^1$ convergence ensures $\abs{E}=\abs{B_1}$.
Now by lower semicontinuity of the perimeter w.r.t. the $L^1$ convergence, we have $P(E)\leq
\edit{\liminf_{k\to\infty}}\, P(\wt{E}_k)$, and by continuity of the nonlocal perimeter in
$L^1(\Ren)$ shown in \cref{lem:perkC0}, it follows that $\Per_{G_\lambda}(\wt{E}_k)$ converges to
$\Per_{G_\lambda}(E)$ as $k$ goes to infinity.
Hence $\calF_{\gamma,G_\lambda}(E) \leq
\edit{\liminf_{k\to\infty}}\,\calF_{\gamma,G_\lambda}(\wt{E}_k)$, which shows that $E$ is a
minimizer of \cref{minrp}, since $\wt{E}_k$ is a minimizing sequence.\\
\textit{Step 3.}
\edit{By the previous steps, for any $\lambda$ large enough, \cref{minrp} admits a minimizer which
satisfies the right inclusion of \cref{existmin:incl}, up to a translation and a negligible set.
There remains to check that this holds for any such minimizer $E$, as well as the reverse inclusion
$\ol{B}_{1-C\eta(\lambda)^{\frac{1}{n}}}\subsq E$. This follows easily by Step~1: consider
$E$ a minimizer} of \cref{minrp} with $\lambda>\lambda_e$, then by minimality, we have
$\calF_{\gamma,G_\lambda}(E)\leq \calF_{\gamma,G_\lambda}(B_1)$, thus applying Step~1 there exists a
set of finite perimeter $\wt{E}$ with mass $\abs{B_1}$ such that
\[
\calF_{\gamma,G_\lambda}(\wt{E}) \leq \calF_{\gamma,G_\lambda}(E)-\abs{E\setminus
B_{1+C\eta(\lambda)^{\frac{1}{n}}}}.
\]
By minimality of $E$, necessarily $\abs{E\setminus B_{1+C\eta(\lambda)^{\frac{1}{n}}}}=0$,
which shows the upper inclusion of \cref{existmin:incl}. For the other inclusion, one can show,
proceeding almost exactly as in Step~1 but using \cref{lem:unionperim} instead of the truncation
lemma and the case $t<1$ of \cref{lem:fscaling}, that if $E$ is a minimizer, then there exists a set
of finite perimeter $\wt{E}$ with mass $\abs{B_1}$ such that
\[
\calF_{\gamma,G_\lambda}(\wt{E}) \leq
\calF_{\gamma,G_\lambda}(E)-\abs{B_{1-C\eta(\lambda)^{\frac{1}{n}}}\setminus E},
\]
which implies $\abs{B_{1-C\eta(\lambda)^{\frac{1}{n}}}\setminus E}=0$.
\edit{Note that the translations of $E$ to obtain the lower and upper inclusions are the same, as
they are given by the ball achieving the minimum in the Fraenkel asymmetry.}
\end{proof}

An immediate consequence of \cref{thm:existmin} is the convergence of minimizers to the unit ball
and of their boundaries to the unit sphere as $\lambda$ goes to infinity.

\begin{cor}[Convergence]\label{cor:CV}
Assume $\gamma<1$. Then for any $\lambda>\lambda_e$ given by \cref{thm:existmin}, any minimizer $F$
of \cref{minrp} satisfies, up to a translation,
\[
\dist_\calH(E,B_1)+\dist_\calH(\partial E,\partial B_1)\leq C\eta(\lambda)^{\frac{1}{n}},
\]
where $\dist_\calH$ denotes the Hausdorff distance, $C$ is a constant depending only on $n$ and
$\gamma$, and $\eta(\lambda)$ is given by \cref{eq:defeta}.
\end{cor}

\subsection{Indecomposability of minimizers}\label{subsec:indecomp}

The aim of this section is to prove connectedness of minimizers of \cref{minrp} when $\gamma<1$ and
$\lambda$ is large enough, or in any case when the kernel $G$ is not compactly supported.
Since any minimizer is always defined up to a set of vanishing Lebesgue measure, and we do not know
yet if there is a precise (partially) regular representative, we work here with a measure theoretic
notion of connectedness for sets of finite perimeter (see e.g.  \autocite{ACMM}), which is referred
to as \textsl{indecomposability}.

\begin{dfn}
We say that a set of finite perimeter $E$ is decomposable if there exist two sets of finite
perimeter $E_1$ and $E_2$ such that $E=E_1\sqcup E_2$, $\abs{E_1}>0$, $\abs{E_2}>0$ and
$P(E)=P(E_1)+P(E_2)$. Naturally, we say that a set of finite perimeter is indecomposable if it is
not decomposable.
\end{dfn}

As with the usual topological notion of connectedness, it is possible to partition a set of finite
perimeter $E$ into indecomposable sets (see \autocite[Theorem 1]{ACMM}) in a unique way (up to sets
of vanishing Lebesgue measure). We call the sets composing this partition the $\calM$-connected
components of $E$. We have the following result establishing a link between the $\calM$-connected
components of a set of finite perimeter and the topological connected components.

\begin{thm}[{\autocite[Theorem 2]{ACMM}}]\label{thm:equivCC}
If $E$ is an open set of finite perimeter such that $\calH^{n-1}(\partial E)=\calH^{n-1}(\partial^*
E)$, then the $\calM$-connected components of $E$ coincide with their topological connected
components.
\end{thm}

\edit{We can show that if $G$ is \textsl{strictly} positive, then any minimizer of \cref{minrp} (if
one exists, no matter the value of $\gamma$ or $\lambda$) is indecomposable.}

\begin{prp}\label{prp:indecomp}
\edit{If $G$ is \edit{strictly positive}, then any minimizer $E$ of \cref{minrp} is indecomposable.}
\end{prp}

\begin{rmk}
In fact, the restriction that $E$ is included in a ball can be dropped, since any minimizer of
\cref{minrp} (if one exists) is in fact included in a ball, as is pointed out later in the proof of
\cref{thm:regminim}.
\end{rmk}

\begin{proof}[Proof of {\cref{prp:indecomp}}]
We proceed by contradiction. Assume that there exists a minimizer $E\subsq B_R$ of \cref{minrp} and
two sets of finite $E_1$ and $E_2$ such that $E=E_1\sqcup E_2$, $\abs{E_1}>0$, $\abs{E_2}>0$ and
$P(E)=P(E_1)+P(E_2)$. Here we prefer to work with the equivalent formulation \cref{origpb} of the
minimization problem, and show that $E$ cannot be a minimizer of $P+\gamma\calV_{G_\lambda}$ under
the volume constraint $\abs{E}=\abs{B_1}$, where $\calV_{G_\lambda}$ is the functional defined by
\[
\calV_{G_\lambda}(E)\coloneqq \iint_{E\times E} G_\lambda(x-y)\dx\dy,
\]
for every measurable set $E$. We have the decomposition
\[
P(E)+\gamma\calV_{G_\lambda}(E)
=P(E_1)+P(E_2)+\gamma\calV_{G_\lambda}(E_1)+\gamma\calV_{G_\lambda}(E_2)
+2\gamma\iint_{E_1\times E_2} G_\lambda(x-y)\dx\dy.
\]
\edit{Let $h\in\Ren$ be such that $\abs{h}>2R$.
Since $E_1$ and $E_2$ are included in $B_R$, then for any $(x,y)\in E_1\times (E_2+h)$ we have
$\abs{x-y}\geq\abs{h}-2R>0$. In particular, defining $E_h\coloneqq E_1\sqcup (E_2+h)$, we have
$\abs{E_h}=\abs{E_1}+\abs{E_2}=\abs{E}$, so that $E_h$ is a valid competitor for $E$.
For the same reason, $P(E_h)=P(E_1)+P(E_2+h)=P(E_1)+P(E_2)$.
Thus we compute}
\begin{equation}\label{indecomp:eq1}
\begin{aligned}
&P(E_h)+\gamma\calV_{G_\lambda}(E_h)\\
&\phantom{P}
=P(E_1)+P(E_2)+\gamma\calV_{G_\lambda}(E_1)+\gamma\calV_{G_\lambda}(E_2)
+2\gamma\iint_{E_1\times (E_2+h)} G_\lambda(x-y)\dx\dy\\
&\phantom{P}
=P(E)+\gamma\calV_{G_\lambda}(E)
+2\gamma\iint_{E_1\times (E_2+h)} G_\lambda(x-y)\dx\dy-2\gamma\iint_{E_1\times E_2}
G_\lambda(x-y)\dx\dy,
\end{aligned}
\end{equation}
where we used the fact that $P(E)=P(E_1)+P(E_2)$ for the last equality.
On the one hand, by a change of variables, we have
\begin{equation}\label{indecomp:eq3}
\iint_{E_1\times(E_2+h)} G_\lambda(x-y)\dx\dy
\leq \abs{E_1}\int_{\edit{B_{\abs{h}-2R}^\compl}} G_\lambda(y) \dy
\,\xrightarrow{\abs{h}\to\infty}\, 0,
\end{equation}
since $G\in L^1(\Ren)$.
On the other hand, since $G$ is \edit{strictly positive}, we have
\[
\int_A G_\lambda(x)\dx > 0,
\]
for any measurable set $A$ such that $\abs{A}>0$.
In particular
\begin{equation}\label{indecomp:eq2}
\iint_{E_1\times E_2} G_\lambda(x-y)\dx\dy = \int_{E_1} \left(\int_{x-E_2} G_\lambda(y)\dy\right)\dx > 0,
\end{equation}
since $\abs{E_1}>0$ and $\abs{E_2}>0$.
Thus by \cref{indecomp:eq1,indecomp:eq2,indecomp:eq3}, for any $\abs{h}$ large enough we have
\[
P(E_h)+\gamma\calV_{G_\lambda}(E_h)<P(E)+\gamma\calV_{G_\lambda}(E),
\]
which contradicts the minimality of $E$.
\end{proof}

\edit{Now, even if $G$ is not strictly positive (for example, if it is compactly supported)},
we know that when $\gamma<1$ and $\lambda$ is large enough, any minimizer of \cref{minrp} lies
between two balls which are very close to the unit ball, which implies that it is indecomposable.
\begin{prp}\label{prp:indecomp2}
Assume $\gamma<1$. Then any minimizer of \cref{minrp} with $\lambda>\lambda_e$ is indecomposable,
where $\lambda_e$ is given by \cref{thm:existmin}.
\end{prp}

\begin{proof}
The proof is trivial. Let $E$ be a minimizer of \cref{minrp} with $\lambda>\lambda_e$, and assume
that it can be written $E=E_1\sqcup E_2$ with $\abs{E_1}>0$, $\abs{E_2}>0$, and
$P(E)=P(E_1)+P(E_2)$.
Since $P(E)=P(E_1)+P(E_2)$ we may translate $E_2$ without changing the perimeter.
Hence, translating it far and outside the ball $B_{1+C\eta(\lambda)^{\frac{1}{n}}}$ by some vector
$h\in\Ren$, we have
\[
\calF_{\gamma,G_\lambda}\big(E_1\sqcup (h+E_2)\big) \leq \calF_{\gamma,G_\lambda}(E),
\]
so that $E_1\sqcup (h+E_2)$ is still a minimizer but outside any ball
$B_{1+C\eta(\lambda)^{\frac{1}{n}}}(y)$, $y\in\Ren$, which contradicts \cref{thm:existmin}.
\end{proof}

\subsection[\texorpdfstring{$\Gamma$}{Gamma}-convergence to the classical perimeter]{%
\for{toc}{\texorpdfstring{$\Gamma$}{Gamma}-convergence to the classical
perimeter}\except{toc}{\texorpdfstring{$\mathbf{\Gamma}$}{Gamma}-convergence to the classical
perimeter}%
}


Using the results from \cref{subsec:reform}, we can easily compute the $\Gamma$-limit of the
functionals of the rescaled problems \cref{minrp}.
In view of \cref{eq:defperglambda}, for any $\lambda\in(0,+\infty)$, let us define on $L^1(\Ren)$ the
functional
\begin{equation}\label{eq:defFeps}
\calE_{\gamma,G_\lambda}(f) \coloneqq \left\{
\begin{aligned}
&\int_{\Ren} \abs{Df}-\gamma\iint_{\Ren\times\Ren} \frac{\abs{f(x)-f(y)}}{\abs{x-y}}
\rho_{G,1/\lambda}(x-y) \\
&{\hphantom{\int_{\Ren} \abs{Df}}}
\text{if there exists a set of finite perimeter $E$ s.t. $f=\ind_E$ and $\abs{E}=\abs{B_1}$},\\
&{+\infty}\nhphantom{$\displaystyle{+\infty}$}
{\hphantom{\int_{\Ren} \abs{Df}}\text{otherwise}},
\end{aligned}
\right.
\end{equation}
which is well-defined and finite whenever $f=\ind_E$ for some set of perimeter $E$ such that
$\abs{E}=\abs{B_1}$ by \cref{prp:boundK1n}. It is obviously defined so that it \enquote{coincides}
with $\calF_{\gamma,G_\lambda}$ on sets of finite perimeter, in the sense that
$\calE_{\gamma,G_\lambda}(\ind_E)=\calF_{\gamma,G_\lambda}(E)$ for every set of finite perimeter~$E$
with volume $\abs{B_1}$, and so that minimizers of $\calE_{\gamma,G_\lambda}$ are precisely those
functions which are indicator functions of sets of finite perimeter solving \cref{minrp}.

\begin{thm}\label{thm:gammacv}
If $\gamma<1$, then the functionals $\calE_{\gamma,G_\lambda}$ defined by \cref{eq:defFeps}
$\Gamma$-converge in $L^1$ to the functional
\[
\calE_{\gamma,\infty}(f) \coloneqq \left\{
\begin{aligned}
&\left(1-\gamma\right)\int_{\Ren} \abs{Df}\\
&\hphantom{\left(1-\frac{I_G}{}\right)}
\text{if there exists a set of finite perimeter $E$ s.t. $f=\ind_E$ and $\abs{E}=\abs{B_1}$},\\
&{+\infty}\nhphantom{$\displaystyle{+\infty}$}
{\hphantom{\left(1-\frac{I_G}{}\right)}}\text{otherwise},
\end{aligned}
\right.
\]
as $\lambda$ goes to infinity.
\end{thm}

\begin{proof}
We shall check, in that order, that
\[
\Gamma-\limsup\,\calE_{\gamma,G_\lambda}(f) \leq \calE_{\gamma,\infty}(f),
\quad\text{ and }\quad
\calE_{\gamma,\infty}(f) \leq \Gamma-\liminf \calE_{\gamma,G_\lambda}(f),
\]
where
\[
\Gamma-\limsup\,\calE_{\gamma,G_\lambda}(f) \coloneqq \min \left\{ \edit{\limsup_{k\to\infty}}~
\calE_{\gamma,G_{\lambda_k}}(f_k) : f_k \,\xrightarrow{L^1(\Ren)}\, f\right\}.
\]
and
\[
\Gamma-\liminf\,\calE_{\gamma,G_\lambda}(f) \coloneqq \min \left\{ \edit{\liminf_{k\to\infty}}~
\calE_{\gamma,G_{\lambda_k}}(f_k) : f_k \,\xrightarrow{L^1(\Ren)}\, f\right\}.
\]
\textit{Step 1.} Let $f\in L^1(\Ren)$.
If $f$ is not the indicator function of a set of finite perimeter of volume $\abs{B_1}$,
$\calE_{\gamma,\infty}(f) = +\infty$ so the \edit{first} inequality is trivial. Let us assume
$f=\ind_F$ for some set of finite perimeter $F$ such that $\abs{F}=\abs{B_1}$, and consider the
constant sequence $f_k \equiv \ind_F$. Then by \cref{prp:limK1n} we have
\[
\calE_{\gamma,G_{\lambda_k}}(\ind_F) \,\xrightarrow{k\to\infty}\, \calE_{\gamma,\infty}(\ind_F),
\]
thus $\Gamma-\limsup \calE_{\gamma,G_\lambda}(f) \leq \calE_{\gamma,\infty}(f)$.\\
\textit{Step 2.} Given $f\in L^1(\Ren)$, consider a sequence $f_k \in L^1(\Ren)$ such that $f_k
\,\xrightarrow{L^1(\Ren)}\, f$.\\
If $f$ is not given by the indicator function of a measurable set $F$, we claim that there exists
$k_0\in\N$ such that for any $k\geq k_0$, $f_k$ is also not an indicator function.
By contradiction, let us assume that there exists a subsequence $(f_{n_k})_{k\in\N}$ such that for
every $n_k$, $f_{n_k}$ is the indicator function of a set $F_{n_k}$. By $L^1$ convergence, up to a
further subsequence (not relabeled), we may assume that $f_{n_k}=\ind_{F_{n_k}}$ converges almost
everywhere to $f$. But then, for almost every $x\in\Ren$, $f(x)=\lim_{k\to\infty}
\ind_{F_{n_k}}(x)\in \{0,1\}$, so that $f$ is an indicator function, which is a contradiction.
Thus for any $k\geq k_0$, $f_{n_k}$ is not an indicator function, hence
$\calE_{\gamma,G_{\lambda_k}(f_k)}=+\infty$, and we indeed have
\begin{equation}\label{gammacv:eq1}
+\infty=\liminf_{k\to\infty}\,\calE_{\gamma,G_{\lambda_k}}(f_k)\geq
\calE_{\gamma,\infty}(f)=+\infty.
\end{equation}
Now we assume that $f=\ind_F$ for some measurable set $F$. If $\abs{F}\neq\abs{B_1}$, then by $L^1$
convergence there exists $k_0$ such that for all $k\geq k_0$, $\int_{\Ren} f(x)\dx\neq \abs{B_1}$,
thus $\calE_{\gamma,G_{\lambda_k}}(f_k)=+\infty$ and \cref{gammacv:eq1} holds as well. Hence we may
now assume that $\abs{F}=\abs{B_1}$.
By \cref{prp:boundK1n}, we have
\[
\calE_{\gamma,G_{\lambda_k}}(f_k) \geq \left(1-\gamma\right) \int_{\Ren} \abs{Df_k},
\]
which trivially holds even if $f_k \not\in \bv(\Ren)$ or $\int_{\Ren} f_k(x)\dx\neq\abs{B_1}$.
Since the $\bv$ seminorm is lower semicontinuous with respect to the usual $L^1$ topology, and
$\big(1-\gamma\big)>0$, we find
\[
\edit{\liminf_{k\to\infty}}~ \calE_{\gamma,G_{\lambda_k}}(f_k) \geq \left(1-\gamma\right)
\int_{\Ren} \abs{Df} =  \calE_{\gamma,\infty}(f),
\]
where we used the fact that $\int_{\Ren} f(x)\dx=\abs{F}=\abs{B_1}$ for the last equality.
Putting these cases together, we get that $\calE_{\gamma,\infty}(f) \leq \Gamma-\liminf
\calE_{\gamma,G_\lambda}(f)$.
\end{proof}

As usual, the $\Gamma$-convergence implies that any converging sequence of minimizers of
\hyperref[minrp]{\textnormal{(P$_{\gamma,\lambda_k}$)}}, associated with positive numbers
$\lambda_k$ going to infinity, necessarily converges in $L^1$ to a minimizer of the perimeter, hence
to the unit ball; however, we already knew this by \cref{cor:CV}.

\subsection{Regularity of minimizers}\label{subsec:reg}

We address here the question of the regularity of minimizers of \cref{minrp}.
Applying \citeauthor{Rig}'s work in \autocite{Rig}, we show that minimizers are in fact
almost-minimizers of the perimeter in the sense of \citeauthor{Tam1984} (see \autocite{Tam1984}),
defined just below, to obtain partial $C^{1,\frac{1}{2}}$-regularity of the boundary of minimizers.

Since sets of finite perimeter are defined up to a Lebesgue negligible set, we shall specify which
boundary we are referring to. Here the boundary we refer to is the support of the Gauss-Green
measure of $E$, which is given by
\[
\spt \mu_E = \Big\{ x\in\Ren ~:~ 0<\abs{E\cap B_r(x)}<\abs{B_r(x)},~\forall r>0\Big\},
\]
and includes the reduced boundary.
It is known that any set of finite perimeter $E$ admits a representative whose topological boundary
agrees with the support of its Gauss-Green measure, namely $\partial E=\spt\mu_E$ (see e.g.
\autocite[Proposition 12.19]{Maggi}; note that the representative built is Borel but not necessarily
open).

\begin{dfn}[Almost-minimizer of the perimeter]\label{dfn:tamalmostmin}
Let $r_0>0$ and $\om:(0,r_0)\to [0,+\infty)$ be a nondecreasing function vanishing in
$0^+$. We say that a set of finite perimeter $E\subsq\Ren$ is an $\om$-almost-minimizer of the
perimeter if
\[
P(E;B_r(x))\leq P(F;B_r(x))+\om(r)r^{n-1}
\]
for every ball $B_r(x)$ with $r_0<r$, and every set of finite perimeter $F$ such that $E\triangle
F\csubset B_r(x)$.
\end{dfn}

Observe that \citeauthor{Tam1984}'s notion of almost-minimizers includes the notion of
$(\Lambda,r_0)$-perimeter minimizers defined by \citeauthor{Maggi} in \autocite[Chapter~21]{Maggi}
when $\om(r)=C r$ for some $C>0$, and that the approach of Tamanini yields stronger results
($C^{1,\frac{1}{2}}$-regularity instead of $C^{1,\alpha}$-regularity for every
$\alpha\in(0,\frac{1}{2})$).

Notice that by \autocite[Lemma 5.2.1]{Rig} (or as a consequence of \cref{lem:perkC0}), any minimizer
of \cref{minrp} satisfies
\begin{equation}\label{eq:volconstquasimin}
P(E) \leq  P(F)+\lambda\gamma I_{G}^{0}\abs{E\triangle F}
\end{equation}
for any set of finite perimeter $F$ such that $\abs{F}=\abs{B_1}$. Let us emphasize that, due to the
\textsl{volume constraint}, this does \textsl{not} immediately imply that minimizers of \cref{minrp}
are $(\Lambda,r_0)$-perimeter minimizers in the sense of \citeauthor{Maggi}, or even
almost-minimizers of the perimeter in the sense of \cref{dfn:tamalmostmin}.

There are several classical ways to deal with the volume constraint, e.g. by making volume-fixing
variations as in \autocite[§17.21]{Maggi}, or by dropping the constraint in the minimization and
introducing instead a penalization of the form $\Lambda\big|\abs{E}-\abs{B_1}\big|$, with
$\Lambda>0$ large enough, in the functional $\calF_{\gamma,G_\lambda}$. Here, we can take advantage
of \citeauthor{Rig}'s work in \autocite{Rig}, which gives precisely that minimizers of
\cref{minrp} are $\om$-almost-minimizers of the perimeter, where $\om:(0,r_0)\to[0,+\infty)$ is
given by $\om(r)=C r$ with $C=C(n,G,\gamma,\lambda)$ and $r_0=r_0(n,G,\gamma,\lambda)$.
Note that $r_0$ and $\om$ depend only on the parameters of the minimization problem, and not on the
minimizer $E$ itself, which would \textsl{a priori} be the case by making volume-preserving
variations to deal with the volume constraint.

Eventually, let us point out that the assumption made in \autocite{Rig} that $G$ is compactly
supported is in fact only used to get existence of a minimizer; the regularity results rely only
on the integrability of $G$ in $\Ren$ (not even on the fact that it is radial).

Combining results from \autocite{Tam1984,Rig,Rig2}, we get the following partial regularity result
for minimizers of~\cref{minrp}.

\begin{thm}\label{thm:regminim}
Let $E$ be a minimizer\footnote{If one exists, no matter whether $\gamma$ is lower than $1$ or not.}
of \cref{minrp}.
Then $\partial^* E$ is locally a $(n-1)$-dimensional graph of class $C^{1,\frac{1}{2}}$, with
$C^{1,\frac{1}{2}}$-regularity constants depending only on $n$, $\gamma$, $\lambda$ and $I_G^0$.
In addition, defining
\[
E_0 \coloneqq \Big\{\,x\in\Ren~:~ \text{there exists }r>0 \text{ s.t. } \abs{B_r(x)\cap
E}=\abs{B_r(x)}\,\Big\},
\]
$E_0$ is an open set equivalent to $E$ such that $\partial E_0=\spt \mu_E$, whose topological
connected components coincide with the $\calM$-connected components of $E$, and it is included in
some ball $B_R$, where $R$ depends only on $n$, $\gamma$, $\lambda$ and $I_G^{0}$.
If $n<8$, then $\partial E_0=\partial^* E$, making the topological boundary of $E_0$ a
$C^{1,\frac{1}{2}}$-hypersurface, and if $n\geq8$, then $\dim_\calH(\partial E_0\setminus\partial^*
E)\leq n-8$.
Furthermore, if $\gamma<1$ and $\lambda>\lambda_e$ (given by \cref{thm:existmin}), or if $G$ is
strictly positive, then $E_0$ is connected.
\end{thm}

\begin{proof}
The fact that $E_0$ is an \textsl{open} set equivalent to $E$ such that $\partial E_0=\spt\mu_E$ is
due to \autocite[Lemma~2.1.3 \& Proposition~2.2.1]{Rig} (more precisely, see equations (2.1.6) \&
(2.2.5) therein, which characterize $\partial E_0$; see also \autocite[Lemma~3.6]{Rig2}).

By \cref{eq:volconstquasimin}, any minimizer $E$ of \cref{minrp} is a \textsl{volume-constrained}
quasi-minimizer of the perimeter in the sense of \autocite{Rig}, thus by \autocite[Proposition
4.3.1]{Rig}, $E$ is a $\om$-almost-minimizer of the perimeter as in \cref{dfn:tamalmostmin}, where
$\om:(0,r_0)\to[0,+\infty)$ is of the form $\om(r)=Cr$, for some positive constants $C$ and $r_0$
depending only on $n$, $\gamma$, $\lambda$ and $I_G^0$.
By \autocite[§1.9 and §1.12]{Tam1984}, $\partial^* E$ is a $C^{1,\frac{1}{2}}$-hypersurface, with
regularity constants depending only on $n$, $\gamma$, $\lambda$ and $I_G^0$. In addition
$\dim_\calH(\spt \mu_E\setminus\partial^* E)\leq n-8$ for $n\geq 8$, and $\spt
\mu_E\setminus\partial^* E=\emptyset$ for $n<8$.

The fact that $E_0$ is included in a ball $B_R$ such that $R=R(n,\gamma,\lambda,I_G^{0})$ comes from
the density estimate
\begin{equation}\label{regmin:eq1}
\abs{E\cap B_r(x)} \geq c\abs{B_r(x)}, \quad\text{ for }\calL^n\text{-a.e. }x\in E\text{ and every }
0<r<r_1,
\end{equation}
where $c=c(n,\gamma,\lambda,I_G^{0})$, $r_1=r_1(n,\gamma,\lambda,I_G^{0})$, which is a consequence
of the results in \autocite{Rig}: more precisely, a uniform version (in the sense that the density
and radius constants depend only on the minimized functional, not on a specific minimizer $E$) of
\autocite[Lemma 2.1.3]{Rig} is obtained in \autocite[Section 4.1]{Rig} (see in particular paragraph
4.1.3 therein), which readily implies \cref{regmin:eq1}.

Eventually, since $E_0$ is open and $\calH^{n-1}(\partial
E_0)=\calH^{n-1}(\spt\mu_E)=\calH^{n-1}(\partial^* E)$, \cref{thm:equivCC} implies that the
$\calM$-connected components of $E$ coincide with the topological connected components of $E_0$.
The connectedness of $E_0$ in the two cases of the theorem is then a direct consequence of
\cref{prp:indecomp,prp:indecomp2}.
\end{proof}


\section{Stability of the ball}\label{sec:stabball}

As in the previous section, we always assume that the kernel $G$ satisfies assumptions
\cref{Krad,Kint}.


\subsection{First and second variations of perimeters}\label{subsec:variations}

In this subsection we recall formulas for the first and second variations of the classical and
nonlocal perimeters, which can be found e.g. in \autocite[Section 6]{FFMMM}. In all this subsection,
$E$ denotes an open set of finite \edit{perimeter} such that $\partial E$ is a $C^2$ hypersurface.
First we define some terminology.

Given a vector field $X\in C^\infty_c(\Ren;\Ren)$, we define the flow induced by $X$ as the solution
in $t\in\R$ of the ODEs
\[
\left\{
\begin{aligned}
&\partial_t\,\Phi_t(x) = X(\Phi_t(x))\\
&\Phi_0(x) = x,
\end{aligned}
\right.
\]
for every $x\in\Ren$.
It is well-known that $\Phi_t(x)$ is well-defined for every $t\in\R$ and $x\in\Ren$, and that
$(\Phi_t)_{t\in\R}$ is a one-parameter group of smooth diffeomorphisms on $\Ren$, i.e.,
$\Phi_t\circ\Phi_s=\Phi_{s+t}$ for all $s,t\in\R$, and $\Phi_0=\id_{|\Ren}$.
Given $\Phi_t$ a flow induced by $X$, we let $E_t\coloneqq\Phi_t(E)$.

\begin{dfn}
We say that a vector field $X\in C^\infty_c(\Ren;\Ren)$ \textsl{induces a volume-preserving flow on
$E$} if there exists $\delta>0$ such that $\abs{E_t}=\abs{E}$ for all $\abs{t}<\delta$.
\end{dfn}

Given a functional $\calF$  on sets of finite perimeter such that $t\mapsto\calF(E_t) \in
C^2(-\delta,\delta)$ for some $\delta>0$, we define the first and second variations of $\calF$ at
$E$ in the direction $X\in C^\infty_c(\Ren)$ by
\[
\delta \calF(E)[X] \coloneqq \left[\frac{\dd}{\dt} \calF(E_t)\right]_{|t=0},\quad \delta^2 \calF(E)[X] \coloneqq 
\left[\frac{\dd^2}{\dt^2} \calF(E_t)\right]_{|t=0}.
\]
Then we define the notion of volume-constrained stationary sets (that is, critical points) for a
functional.

\begin{dfn}
We say that $E$ is a volume-constrained stationary set for the functional $\calF$ if
for every $X\in C^\infty_c(\Ren;\Ren)$ inducing a volume-preserving flow on $E$, we have
$\delta\calF(E)[X]=0$.
\end{dfn}

We are interested in the variations of the classical perimeter $P$ and of the nonlocal perimeter
$\Per_G$, which we will deduce from the variations of the nonlocal term
\begin{equation}\label{eq:relpercalv}
\calV_G(E) \coloneqq \iint_{E\times E} G(x-y)\dx\dy = \abs{E}I_G^{0}-\Per_G(E).
\end{equation}
For the classical perimeter, it is known that $t\mapsto P(E_t)$ is smooth in $(-\delta,\delta)$
whenever $E$ is a set of finite perimeter, and if $\partial E$ is a $C^2$-hypersurface,
the first variation is
\[
\delta P(E)[X] = \int_{\partial E} H_{\partial E}\, \zeta_X \dH^{n-1},
\]
where $\zeta_X\coloneqq X\cdot\nu_E$, $\nu_E$ denotes the outer unit normal to $E$, and $H_{\partial
E}$ is the scalar mean curvature of $\partial E$.
The second variation is given by
\[
\delta^2 P(E)[X]
=\int_{\partial E} \abs{\grad_\tau\,\zeta_X}^2-c_{\partial E}^2\,\zeta_X^2\dH^{n-1}
+\int_{\partial E} H_{\partial E}\big((\dvg X)\zeta_X-\dvg_\tau(\zeta_X X_\tau)\big)\dH^{n-1},
\]
where $c_{\partial E}^2(x)$ is the sum of the squares of the principal curvatures of $\partial E$ at
$x$, $X_\tau\coloneqq X-\zeta_X\nu_E$, and $\grad_\tau$ and $\dvg_\tau$ denote respectively the tangential
gradient and divergence on $\partial E$.
In addition, if $E$ is a volume-constrained stationary set for the perimeter, and $X$ induces a
volume-preserving flow on $E$, then the second variation of the perimeter takes the simpler form
\begin{equation}\label{eq:var2perimstat}
\delta^2 P(E)[X]
=\int_{\partial E} \abs{\grad_\tau\, \zeta_X}^2-c_{\partial E}^2\,\zeta_X^2\dH^{n-1}.
\end{equation}
Indeed, in that case $H_{\partial E}$ is constant, and the fact that $t\mapsto \abs{E_t}$ is
constant in a neighborhood of $0$ implies
\[
0=\left[\ddt\, \abs{E_t}\right]_{|t=0} = \int_{\partial E} \zeta_X \dH^{n-1}\quad\text{ and }\quad
0=\left[\frac{\dd^2}{\dt^2}\, \abs{E_t}\right]_{|t=0}=\int_{\partial E} (\dvg X)\zeta_X \dH^{n-1}.
\]
As for the \edit{first and second variations} of $\calV_G$, \edit{assuming that \edit{$G\in
C^1(\Ren\setminus\{0\})$}, $G(x)=o(\abs{x}^{\alpha-n})$ at the origin and
$G(x)=o(\abs{x}^{-(n+\beta)})$ at infinity for some $\alpha>0$ and $\beta>0$, since $E$ is an open
set with finite volume such that $\partial E$ is a $C^{2}$-hypersurface, we can apply directly
\autocite[Theorem~6.1]{FFMMM} to get}
\begin{equation}\label{eq:varG}
\begin{aligned}
\delta \calV_G(E)[X] &= \int_{\partial E} H^*_{G,\partial E}\,\zeta_X\dH^{n-1},\\
\delta^2 \calV_G(E)[X] &= -\iint_{\partial E\times\partial E}
G(x-y)\abs{\zeta_X(x)-\zeta_X(y)}^2\dH_x^{n-1}\dH_y^{n-1}+\int_{\partial E} c^2_{G,\partial
E}\,\zeta_X^2\dH^{n-1}\\
&\phantom{=-\iint}
+\int_{\partial E} H^*_{G,\partial E}\big((\dvg X)\zeta_X-\dvg_\tau(\zeta_X X_\tau)\big)\dH^{n-1},
\end{aligned}
\end{equation}
 where
\[
c^2_{G,\partial E}(x) \coloneqq \int_{\partial E} G(x-y)\abs{\nu_E(x)-\nu_E(y)}^2\dH_y^{n-1},\quad\forall
x\in\partial E,
\]
and $H_{G,\partial E}^*$, which plays the role of a nonlocal mean curvature of $\partial E$, is
defined by
\begin{equation}\label{eq:defnonloccurv}
H^*_{G,\partial E}(x) \coloneqq 2\int_E G(x-y)\dy,\quad\forall x\in\partial E.
\end{equation}
Note that all the integrals in \cref{eq:varG} are finite whenever $\partial E$ is a
$C^2$-hypersurface by the assumptions on $G$, since $X_\tau$, $\zeta_X$ and $\nu_E$ are bounded and
$C^1$ functions.
Similarly to the perimeter functional, if $E$ is a volume-constrained stationary set for $\calV_G$,
and $X$ induces a volume-preserving flow on $E$, the fact that $t\mapsto\abs{E_t}$ is constant in a
neighborhood of $0$ implies that the second variation of $\calV_G$ is simply given by
\[
\begin{aligned}
\delta^2 \calV_G(E)[X]
= -\iint_{\partial E\times\partial E}
G(x-y)\abs{\zeta_X(x)-\zeta_X(y)}^2\dH_x^{n-1}\dH_y^{n-1}
+\int_{\partial E} c^2_{G,\partial E}\,\zeta_X^2\dH^{n-1}.
\end{aligned}
\]
Recalling that,
\[
\Per_G(E_t) = I_G^0\abs{E_t}-\calV_G(E_t),\quad\quad\forall t\in(-\delta,\delta),
\]
we see that $E$ is a volume-constrained stationary set for $\calV_G$ if and only if it is such a set
for $\Per_G$. Hence, we get the following expression of the second variation of the nonlocal
perimeter.

\begin{prp}\label{prp:var2perimkstat}
Assume that \edit{$G\in C^1(\Ren\setminus\{0\})$}, $G(x)=o(\abs{x}^{\alpha-n})$ at the origin, and
$G(x)=o(\abs{x}^{-(n+\beta)})$ at infinity for some positive constants $\alpha$ and $\beta$. Then,
if $X$ induces a volume-preserving flow on $E$, we have
\[
\begin{aligned}
\delta^2 \Per_G(E)[X] 
&= -\delta^2 \calV_G(E)[X]\\
&= \iint_{\partial E\times\partial E}
G(x-y)\abs{\zeta_X(x)-\zeta_X(y)}^2\dH_x^{n-1}\dH_y^{n-1}
-\int_{\partial E} c^2_{G,\partial E}\,\zeta_X^2\dH^{n-1}.
\end{aligned}
\]
\end{prp}

We end this section by recalling the definition of stability in that setting.

\begin{dfn}\label{dfn:stability}
We say that $E$ is a volume-constrained stable set for the functional $\calF$ if~$E$ is a
volume-constrained stationary set for $\calF$, and $\delta^2 \calF(E)[X] \geq 0$ for every vector
field $X\in C^\infty_c(\Ren;\Ren)$ inducing a volume-preserving flow on $E$.
\end{dfn}

Obviously, if $E$ is a volume-constrained minimizer for some functional $\calF$, it is a
volume-constrained stable set for $\calF$ in this sense.

\subsection{The stability threshold for the ball}\label{subsec:stabball}

We are interested in the stability of the unit ball for \cref{minrp}.

\subsubsection[Instability when \texorpdfstring{$\gamma>1$}{gamma is larger than 1}]{Instability when \texorpdfstring{\forcebold{$\gamma>1$}}{gamma is larger than 1}}

First, we give the expression of the second variation of $\calF_{\gamma,G_\lambda}$ at $B_1$.

\begin{prp}[Second variation]\label{prp:var2fk}
\edit{Assume $G\in C^1(\Ren\setminus\{0\})$.}
For every $\gamma>0$, $\lambda>0$, the unit ball $B_1$ is a volume-constrained stationary set for
$\calF_{\gamma,G_\lambda}$. Furthermore, if $G(x)=o(\abs{x}^{\alpha-n})$ at the origin, and
$G(x)=o(\abs{x}^{-(n+\beta)})$ at infinity for some positive constants $\alpha$ and $\beta$, then
for any vector field $X\in C^\infty_c(\Ren;\Ren)$ inducing a volume-preserving flow on $B_1$, the
second variation of $\calF_{\gamma,G_\lambda}$ at $B_1$ in the direction $X$ is given by
\begin{equation}\label{eq:var2fkexp}
\begin{aligned}
&\delta^2 \calF_{\gamma,G_\lambda}(B_1)[X]\\
&= \int_{\bbS^{n-1}} \abs{\grad_\tau\,\zeta_X}^2\dH^{n-1}
-\gamma\iint_{\bbS^{n-1}\times\bbS^{n-1}}
\frac{\abs{\zeta_X(x)-\zeta_X(y)}^2}{\abs{x-y}^2}\eta_{G,1/\lambda}(\abs{x-y})\dH_x^{n-1}\dH_y^{n-1}\\
&\phantom{=\int_{\bbS^{n-1}}}
-\left(c^2_{\bbS^{n-1}}-\gamma\, c_{G,\lambda,\bbS^{n-1}}^2 \right)\int_{\bbS^{n-1}}
\zeta_X^2\dH^{n-1}.
\end{aligned}
\end{equation}
where $c^2_{\bbS^{n-1}}=n-1$ is the sum of the squares of the principal curvatures of $\bbS^{n-1}$,
\edit{the quantity $c_{G,\lambda,\bbS^{n-1}}$ is defined by}
\[
c_{G,\lambda,\bbS^{n-1}}^2 \coloneqq \int_{\bbS^{n-1}} \eta_{G,1/\lambda}(\abs{x-y}) \dH_y^{n-1},
\qquad\text{ for any }x\in\bbS^{n-1},
\]
\edit{(by symmetry, it is independent of $x\in\bbS^{n-1}$), and} for all $r\in(0,+\infty)$ and
$\eps>0$, we defined
\begin{equation}\label{eq:defetag}
\eta_G(r) \coloneqq r^2g(r),\quad\quad\quad \eta_{G,\eps}(r)\coloneqq \eps^{-(n-1)}\eta_G(\eps^{-1}
r).
\end{equation}
\end{prp}

\begin{proof}
\edit{Since the unit ball is a minimizer of the perimeter under volume constraint, it is a
volume-constrained stationary set for the perimeter.
As for the nonlocal perimeter, it is not necessarily minimized by the ball, since $G$ is not assumed
to be radially nonincreasing. However, one can notice from \cref{eq:defnonloccurv} that
$H^*_{G,\bbS^{n-1}}$ is constant by symmetry, which directly gives the stationarity of $B_1$ for
$\calV_{G_\lambda}$ and thus for $\Per_{G_\lambda}$, in view of \cref{eq:relpercalv} and the
expression of the first variation of $\calV_{G_\lambda}$ in \cref{eq:varG}.}
Thus $B_1$ is a stationary set for $\calF_{\gamma,G_\lambda}$.  Since $B_1$ has a smooth boundary,
applying \cref{prp:var2perimkstat} to $\Per_{G_\lambda}$, we find
\begin{equation}\label{eq:var2vk1}
\begin{aligned}
\delta^2 \Per_{G_\lambda}(B_1)[X]
&= \iint_{\bbS^{n-1}\times\bbS^{n-1}}
G_\lambda(x-y)\abs{\zeta_X(x)-\zeta_X(y)}^2\dH_x^{n-1}\dH_y^{n-1}\\
&\phantom{=\iint_{\bbS^{n-1}}}
-\iint_{\bbS^{n-1}\times\bbS^{n-1}} G_\lambda(x-y)\abs{x-y}^2
\zeta_X(x)^2\dH_x^{n-1}\dH_y^{n-1}
\end{aligned}
\end{equation}
for every $X\in C^\infty_c(\Ren;\Ren)$ inducing a volume-preserving flow on $B_1$.\\
Now, we can rewrite \cref{eq:var2vk1} in terms of $\eta_{G,1/\lambda}$ by
\[
\begin{aligned}
\delta^2 \Per_{G_\lambda}(B_1)[X]
&= \iint_{\bbS^{n-1}\times\bbS^{n-1}}
\frac{\abs{\zeta_X(x)-\zeta_X(y)}^2}{\abs{x-y}^2}\eta_{G,1/\lambda}(x-y)\dH_x^{n-1}\dH_y^{n-1}\\
&\phantom{=\frac{I_G^1\abs{\bbS^{n-2}}}{\abs{\bbS^{n-1}}}\iint_{\bbS^{n-1}}}
-\iint_{\bbS^{n-1}\times\bbS^{n-1}}
\eta_{G,1/\lambda}(x-y) \zeta_X(x)^2\dH_x^{n-1}\dH_y^{n-1}.
\end{aligned}
\]
Note that
\[
c_{G,\lambda,\bbS^{n-1}}^2 = \int_{\bbS^{n-1}} \eta_{G,1/\lambda}(x-y) \dH_y^{n-1}
\]
does not depend on $x\in\bbS^{n-1}$, since $\eta_{G,1/\lambda}$ is invariant under rotations, thus
by Fubini's theorem we find
\[
\begin{aligned}
\delta^2 \Per_{G_\lambda}(B_1)[X]
&= \iint_{\bbS^{n-1}\times\bbS^{n-1}}
\frac{\abs{\zeta_X(x)-\zeta_X(y)}^2}{\abs{x-y}^2}\eta_{G,1/\lambda}(x-y)\dH_x^{n-1}\dH_y^{n-1}\\
&\phantom{=\frac{I_G^1\abs{\bbS^{n-2}}}{\abs{\bbS^{n-1}}}\iint_{\bbS^{n-1}}}
\vphantom{\iint_{\bbS^{n-1}\times\bbS^{n-1}}}
-c_{G,\lambda,\bbS^{n-1}}^2\int_{\bbS^{n-1}} \zeta_X^2\dH^{n-1},
\end{aligned}
\]
which concludes the proof, since the second variation of the perimeter is given by
\cref{{eq:var2perimstat}}.
\end{proof}

Similarly to the family $(\rho_{G,\eps})_{\eps>0}$ introduced in \cref{eq:defrhog} for rewriting
$\Per_{G_\lambda}$, we see that the family $(\eta_{G,\eps})_{\eps>0}$ is, up to multiplication by a
constant, an approximation of identity as well, as defined in \cref{dfn:mollif} -- however, here it
is a $\mathbf{(n-1)}$-dimensional approximation of identity.

\begin{lem}\label{lem:etaapproxid}
After multiplication by $\bfK_{2,n-1}$, the family $(\eta_{G,\eps})_{\eps>0}$ defined by
\cref{eq:defetag} is a $(n-1)$-dimensional approximation of identity.
\end{lem}

\begin{proof}
First, let us check that
\[
\bfK_{2,n-1}\abs{\bbS^{n-2}} \int_0^{\infty} r^{n-2}\eta_G(r)\dr=1.
\]
Using that
\[
I_G^1=\frac{2}{\bfK_{1,n}}
\]
and the expression of $\bfK_{p,n}$ given by \cref{lem:expKpn}, we find
\[
\bfK_{2,n-1}\abs{\bbS^{n-2}} \int_0^{\infty} \eta_G(r)\edit{r^{n-2}}\dr
=\frac{\bfK_{2,n-1}\abs{\bbS^{n-2}}}{\abs{\bbS^{n-1}}} I_G^1
= \frac{2\bfK_{2,n-1}\abs{\bbS^{n-2}}}{\bfK_{1,n}\abs{\bbS^{n-1}}}
=1.
\]
By a change of variables, we see that
\[
\int_0^\infty r^{n-2}\eta_{G,\eps}(r)\dr = \int_0^\infty r^{n-2}\eta_G(r)\dr,
\]
thus there remains only to check that for every $\delta>0$, we have
\[
\int_\delta^\infty \eta_{G,\eps}(r)r^{n-2} \dr\xrightarrow{\eps\to 0} 0.
\]
By yet another change of variables
\[
\int_\delta^\infty r^{n-2}\eta_{G,\eps}(r)\dr
=\frac{1}{\abs{\bbS^{n-1}}}\int_{B_\delta^\compl} \frac{\eta_{G,\eps}(\abs{x})}{\abs{x}}\dx
=\frac{1}{\abs{\bbS^{n-1}}}\int_{B_{\frac{\delta}{\eps}}^\compl} \abs{x}G(x)\dx,
\]
which vanishes as $\eps$ goes to $0$, since $G$ is integrable on $\Ren$ w.r.t. to the measure
$\abs{x}\dx$ by~\cref{Kint}.
\end{proof}

To see if the unit ball is stable or not when $\lambda$ is large, we wish to pass to the limit in
the second variation of $\calF_{\gamma,G_\lambda}$ at $B_1$, but since we integrate over
$\bbS^{n-1}$ instead of $\Ren$, we cannot use \cref{prp:limK1n}.
In \autocite[Theorem~1.1]{KreMor}, an equivalent to this proposition is given for smooth Riemannian
manifolds, however, the requirements on the family $(\eta_{G,\eps})_{\eps>0}$ seem too strong:
namely the requirement that the functions $\eta_{G,\eps}$ are nonincreasing, i.e., that $r\mapsto
r^2 g(r)$ is nonincreasing. 


\edit{We establish the following proposition that we use to prove that the unit ball is unstable for
$\lambda$ large enough, provided that $\gamma>1$. The proof of this proposition is postponed to
\cref{app:mollifsphere}.}
\edit{\begin{prp}\label{prp:limK1nsph}
Let $(\eta_\eps)_{\eps>0}$ be a $(n-1)$-dimensional approximation ou
identity. In dimension $n=2$, assume that in addition, it satisfies \cref{eq:Tcond}.
Then for any $u\in H^1(\bbS^{n-1})$, we have
\[
\lim_{\eps\to 0}\, \iint_{\bbS^{n-1}\times\bbS^{n-1}} \frac{\abs{u(x)-u(y)}^2}{\abs{x-y}^2}
\eta_\eps(\abs{x-y}) \dH^{n-1}_x\dH^{n-1}_y
= \bfK_{2,n-1} \int_{\bbS^{n-1}} \abs{\grad_\tau\,u}^2 \dH^{n-1}.
\]
\end{prp}}

\edit{Observe that this convergence is not enough to show that the unit ball is \textit{stable} when
$\gamma<1$, provided that $\lambda$ is large enough. In order to show this, we also need an upper
bound similar to \cref{prp:boundK1n}, which is obtained in \cref{app:mollifsphere}.}

\begin{thm}\label{thm:unstable}
Assume that $G$ satisfies assumptions \crefrange{Krad}{KC1}.
If $\gamma>1$, there exists $\lambda_u=\lambda_u(n,G\edit{,\gamma})$ such that for any
$\lambda>\lambda_u$ the unit ball $B_1$ is an \textbf{unstable} critical point of the functional
$\calF_{\gamma,G_\lambda}$.
\end{thm}

\begin{proof}
Consider a vector field $X\in C^\infty_c(\Ren;\Ren)$ inducing a volume-preserving flow on $B_1$ such
that $\delta^2 P(B_1)[X] > 0$.
Note that \edit{$G\in C^1(\Ren\setminus\{0\})$, $G(x)=o(\abs{x}^{\alpha-n})$ at the origin and
$G(x)=o(\abs{x}^{-(n+\beta)})$ by \cref{KC1}, \ref{Korig} and \cref{Ksuperdec} respectively}, thus the second
variation of $\calF_{\gamma,G_\lambda}$ at $B_1$ is given by \cref{prp:var2fk}, and by
\cref{lem:etaapproxid}, up to multiplying my $\bfK_{2,n-1}$, the family $(\eta_{G,\eps})_{\eps>0}$
is a $(n-1)$-dimensional approximation of identity.
In order to apply \cref{prp:limK1nsph}, we need to check that this family also satisfies assumption
\cref{eq:Tcond} in dimension $n=2$.
Let $r>0$. Given $\delta>0$, by \cref{Ksuperdec}, there exists $t_*$ such that, for any
$t>t_*$, we have $t^{3}g(t)\leq \delta r$. Thus there exists $\eps_*$ such that, for every
$\eps<\eps_*$ and every $s>r$, we have $\eps^{-1}s>t_*$, hence
\[
\left(\frac{s}{\eps}\right)^{3}g\left(\frac{s}{\eps}\right) \leq
\delta r,\qquad\forall s>r,
\]
which implies
\[
\eta_{G,\eps}(s)
=\eps^{-3}s^2g(\eps^{-1}s)
=\frac{1}{s} \left(\frac{s}{\eps}\right)^{3}g\left(\frac{s}{\eps}\right)
\leq \delta,\qquad\forall s>r,
\]
so that $\sup_{s>r}\, \eta_{G,\eps}(s)$ vanishes as $\eps$ goes to $0$, i.e., the family
$(\eta_{G,\eps})_{\eps>0}$ satisfies \cref{eq:Tcond} in dimension $n=2$, by the arbitrariness of
$r$.
Hence we can always apply \cref{prp:limK1nsph}, which gives
\begin{equation}\label{unstable:eq1}
\lim_{\lambda\to\infty}\, \iint_{\bbS^{n-1}\times\bbS^{n-1}}
\frac{\abs{\zeta_X(x)-\zeta_X(y)}^2}{\abs{x-y}^2} \eta_{G,1/\lambda}(x-y) \dH_x^{n-1}\dH_y^{n-1}= 
\int_{\bbS^{n-1}} \abs{\grad_\tau\, \zeta_X}^2 \dH^{n-1}.
\end{equation}
Similarly, we compute
\begin{equation}\label{unstable:eq2}
\begin{aligned}
c_{G,\lambda,\bbS^{n-1}}^2
&=\frac{1}{\abs{\bbS^{n-1}}}\iint_{\bbS^{n-1}\times \bbS^{n-1}}
\eta_{G,1/\lambda}(x-y)\dH_x^{n-1}\dH_y^{n-1}\\
&=\frac{1}{\abs{\bbS^{n-1}}}\iint_{\bbS^{n-1}\times \bbS^{n-1}}
\frac{\abs{x-y}^2}{\abs{x-y}^2}\eta_{G,1/\lambda}(x-y)\dH_x^{n-1}\dH_y^{n-1}\\
&\xrightarrow{\lambda\to \infty}
\frac{1}{\abs{\bbS^{n-1}}}\int_{\bbS^{n-1}} \abs{\grad\tau\, \id}^2
\dH^{n-1}=n-1=c^2_{\bbS^{n-1}}.
\end{aligned}
\end{equation}
Combining \cref{unstable:eq1,unstable:eq2} with \cref{prp:var2fk}, we find
\[
\begin{aligned}
\lim_{\lambda\to\infty} \delta^2\calF_{\gamma,G_\lambda}(B_1)[X]
&= \left(1-\gamma\right)
\left(\int_{\bbS^{n-1}} \abs{\grad_\tau\, \zeta_X}^2-c_{\bbS^{n-1}}^2\,\zeta_X^2\dH^{n-1}\right)\\
&= \left(1-\gamma\right)\delta^2 P(B_1)[X]<0,
\end{aligned}
\]
since $\gamma>1$ and $\delta^2 P(B_1)[X]>0$. This shows that there exists $\lambda_u>0$ such that
for any $\lambda>\lambda_u$, the unit ball is unstable for $\calF_{\gamma,G_\lambda}$, which
	concludes the proof.
\end{proof}

\subsubsection{Decomposition in spherical harmonics}

We have just seen that if $\gamma>1$, then the unit ball is unstable for large values of $\lambda$.
Conversely, we prove that if $\gamma<1$, then the unit is stable for large values of $\lambda$,
which shows that this stability threshold is sharp. In order to prove this, we are going to use the
decomposition of the Jacobi operator associated with the second variation of
$\calF_{\gamma,G_\lambda}$ in spherical harmonics.

Looking at the quadratic form of the second variation given by \cref{prp:var2fk}, we see that the
three different terms on the right-hand side of \cref{eq:var2fkexp} derive, first, from the
Laplace-Beltrami operator on the sphere, second, the identity map on the sphere, and third, the
integral operator
\begin{equation}\label{eq:defscrR}
\scrR_{G,\lambda}(u)(x) \coloneqq 2\,\int_{\bbS^{n-1}}
\frac{u(x)-u(y)}{\abs{x-y}^2}\eta_{G,1/\lambda}(\abs{x-y})\dH_y^{n-1},
\end{equation}
(we will justify later that it is well-defined on $C^1(\bbS^{n-1})$).
Indeed, formally, we have
\[
\int_{\bbS^{n-1}} \abs{\grad_\tau\, u}^2\dH^{n-1} =
\int_{\bbS^{n-1}} (-\Delta_{\bbS^{n-1}}\, u) u\dH^{n-1}
\]
and, by symmetry,
\[
\int_{\bbS^{n-1}} u\,\scrR_{G,\lambda}(u)\dH^{n-1}= \iint_{\bbS^{n-1}\times\bbS^{n-1}}
\frac{\abs{u(x)-u(y)}2}{\abs{x-y}^2} \eta_{G,1/\lambda}(\abs{x-y})\dH^{n-1}_x\dH^{n-1}_y.
\]
\edit{It is well-known that spherical harmonics of a given degree are eigenfunctions of the
Laplace-Beltrami operator on the sphere. We show below that they are also eigenfunctions for the
integral operators~$\scrR_{G,\lambda}$}. In the end, to establish the stability of the ball, one
needs only study the eigenvalues of the three operators, and deduce that the eigenvalues of the
operator associated with the second variation of $\calF_{\gamma,G_\lambda}$ are all positive for
large values of $\lambda$.

Let us now recall a few basic facts on spherical harmonics and justify that the $\scrR_{G,\lambda}$
are well-defined. For $k\in\N$, we denote by $\calS_k$ the finite dimensional subspace of
$L^2(\bbS^{n-1})$ made of the spherical harmonics of degree $k$, and $(Y_k^i)_{1\leq i\leq d(k)}$ an
orthonormal basis of $\calS_k$ in $L^2(\bbS^{n-1})$. When there can be no confusion, we denote by
$Y_k$ a generic element of $\calS_k$. Let us recall that $Y_k\in\calS_k$ and $Y_l\in\calS_l$ are
orthogonal in $L^2(\bbS^{n-1})$ whenever $k\neq l$, and that the family
\[
(Y_k^i)_{\substack{0\leq k<\infty,\\1\leq i\leq d(k)}}
\]
is total, i.e., it is an orthonormal basis of $L^2(\bbS^{n-1})$. For any $u\in
L^2(\bbS^{n-1})$, we denote by $a^i_k(u)$ its $(i,k)$-coordinate in the basis, that is,
\begin{equation}\label{eq:decompharmf}
u = \sum_{k=0}^\infty \sum_{i=1}^{d(k)} a^i_k(u)Y^i_k.
\end{equation}

In the following, we will often use the following to integrate $\bbS^{n-1}$ by
slices (see e.g. \autocite[Section D.2]{Gra2014} or \autocite[Corollary A.6]{ABW}):
a $\calH^{n-1}$-measurable function $f$ is integrable on $\bbS^{n-1}$ if and only
$(x,t)\mapsto(1-t^2)^{\frac{n-3}{2}}f(\sqrt{1-t^2}x,t)$ is integrable on $\bbS^{n-2}\times(-1,1)$,
and in that case we have
\begin{equation}\label{eq:intsphere}
\int_{\bbS^{n-1}} f\dH^{n-1} = \int_{-1}^1 (1-t^2)^{\frac{n-3}{2}} \int_{\bbS^{n-2}}
f(\sqrt{1-t^2}x,t) \dH^{n-2}_x\dt.
\end{equation}

\begin{lem}\label{lem:intsphere}
Let $f$ be a nonnegative measurable function on $(0,2)$.
Then for every $x\in \bbS^{n-1}$, the map $F:y\mapsto f(\abs{x-y})$ belongs to $L^1(\bbS^{n-1})$ if
and only if
\[
\int_{\bbS^{n-1}} f(\abs{x-y})\dH^{n-1}_y = \abs{\bbS^{n-2}} \int_0^2
\left(1-\frac{s^2}{4}\right)^{\frac{n-3}{2}} s^{n-2} f(s) \dd s < \infty.
\]
If $n\geq 3$ and
\begin{equation}\label{eq:intspherehyp}
\int_0^2 f(r)r^{n-2}\dr<+\infty,
\end{equation}
then $F\in L^1(\bbS^{n-1})$ and we have,
\begin{equation}\label{eq:intsphereineq1}
\int_{\bbS^{n-1}} f(\abs{x-y})\dH^{n-1}_y \leq \abs{\bbS^{n-2}} \int_0^{2}
s^{n-2} f(s) \dd s.
\end{equation}
For $n=2$, if \cref{eq:intspherehyp} holds and we have $f(r)\leq M$ on $(1,2)$, for some
$M>0$, then $F\in L^1(\bbS^{n-1})$ and
\begin{equation}\label{eq:intsphereineq2}
\int_{\bbS^{n-1}} f(\abs{x-y})\dH^{n-1}_y \leq \abs{\bbS^{n-2}} \left(\frac{2}{\sqrt{3}}
\int_0^{1} f(s)s^{n-2}\dd
s+\frac{M\pi}{3}\right).
\end{equation}
\end{lem}

\begin{proof}[Proof of {\cref{lem:intsphere}}]
Up to a change of variables, we can assume that $x=N=(0,\dotsc,0,1)$ is the ``north pole''. Our
computations will show that $y\mapsto f(\abs{N-y})\in L^1(\bbS^{n-1})$.
Applying \cref{eq:intsphere} to $y\mapsto f(\abs{N-y})$ we find
\[
\begin{aligned}
\int_{\bbS^{n-1}} f(\abs{N-y})\dH^{n-1}_y
&= \int_{-1}^1 (1-t^2)^{\frac{n-3}{2}} \int_{\bbS^{n-2}}
f\big(\big|N-(\sqrt{1-t^2}y,t)\big|\big)\dH^{n-2}_y\dt\\
&= \int_{-1}^1 (1-t^2)^{\frac{n-3}{2}} \int_{\bbS^{n-2}}
f\big(\sqrt{2(1-t)}\big)\dH^{n-2}_y\dt\\
&= \abs{\bbS^{n-2}}\int_{-1}^1 (1-t^2)^{\frac{n-3}{2}} f\big(\sqrt{2(1-t)}\big)\dt.
\end{aligned}
\]
Changing variables with $s=\sqrt{2(1-t)}$, it follows
\[
\begin{aligned}
\int_{\bbS^{n-1}} f(\abs{N-y})\dH^{n-1}_y
&= \abs{\bbS^{n-2}} \int_0^{2}
\left[\Big(s^2\Big(1-\frac{s^2}{4}\Big)\right]^{\frac{n-3}{2}} f(s) s\dd s\\
&= \abs{\bbS^{n-2}} \int_0^{2} \left(1-\frac{s^2}{4}\right)^{\frac{n-3}{2}} s^{n-2} f(s) \dd s.
\end{aligned}
\]
There remains to show that the integral on the right-hand side is finite, with the desired upper
bounds.
When $n\geq 3$, we have
\[
\left(1-\frac{s^2}{\edit{4}}\right)^{\frac{n-3}{2}} \leq 1,\quad\forall s\in (0,2),
\]
which gives the required estimate \cref{eq:intsphereineq1} and shows that the integral is finite.
When $n=2$, assuming that $f(r)\leq M$ on $\edit{(1,2)}$ and splitting the integral into two parts
gives
\begin{equation}\label{intsphere:eq1}
\begin{aligned}
\int_0^{2} \left(1-\frac{s^2}{4}\right)^{-\frac{1}{2}} f(s)s^{n-2}\dd s
&= \int_0^{1} \left(1-\frac{s^2}{4}\right)^{-\frac{1}{2}} f(s)s^{n-2}\dd s
+ \int_1^2 \left(1-\frac{s^2}{4}\right)^{-\frac{1}{2}} f(s)\dd s\\
&\leq \frac{2}{\sqrt{3}} \int_0^1 f(s)s^{n-2}\dd s+\frac{M\pi}{3},
\end{aligned}
\end{equation}
hence the required estimate \cref{eq:intsphereineq2}.
\end{proof}

\begin{prp}\label{prp:scrRwelldef}
If $n=2$, assume that $G(x)\leq M$ on $B_{R}^\compl$, for some positive constants $M$ and $R$, and
if $n\geq 3$, let us just set $R:=0$ and assume no upper bound on $G$.
Then for every $\lambda>R$, the operator $\scrR_{G,\lambda}$ defined by \cref{eq:defscrR} is an
operator from $C^1(\bbS^{n-1})$ to $L^2(\bbS^{n-1})$.
If in addition $G(x)=o(\abs{x}^{\alpha-n})$ near the origin for some $\alpha>1$, then for every
$\lambda>R$, $\scrR_{G,\lambda}$ is actually an operator from $L^2(\bbS^{n-1})$ to
$L^2(\bbS^{n-1})$. In any case, we have
\begin{equation}\label{eq:seminormKop}
\int_{\bbS^{n-1}} u\,\scrR_{G,\lambda}(u)\dH^{n-1}= \iint_{\bbS^{n-1}\times\bbS^{n-1}}
\frac{\abs{u(x)-u(y)}2}{\abs{x-y}^2} \eta_{G,1/\lambda}(\abs{x-y})\dH^{n-1}_x\dH^{n-1}_y
\end{equation}
whenever $u\in C^1(\bbS^{n-1})$.
\end{prp}

\begin{rmk}\label{rmk:scrRwelldef}
The extra assumption that $G(x)\leq M$ on some $B_{R}^\compl$ in dimension $n=2$ is justified by
\cref{lem:intsphere}. In any case, it holds if $G$ satisfies \cref{Ksuperdec}.
\end{rmk}

\begin{proof}[Proof of {\cref{prp:scrRwelldef}}]
Assume that $\lambda>R$.
We first focus on the case where $G$ satisfies the stronger assumption $G(x)=o(\abs{x}^{\alpha-n})$
near the origin.\\
\textit{Case 1.} Assume $G=o(\abs{x}^{\alpha-n})$ near the origin, for some $\alpha>1$.
In view of \cref{lem:intsphere}, we have
\begin{equation}\label{sphereharm:eq1}
\begin{aligned}
\int_{\bbS^{n-1}} \frac{\eta_{G,1/\lambda}(\abs{x-y})}{\abs{x-y}^2} \dH_y^{n-1}
&=\lambda^{n+1}\int_{\bbS^{n-1}} g(\lambda\abs{x-y})\dH_y^{n-1}\\
&=\lambda^{n+1}\abs{\bbS^{n-2}} \int_0^2 \left(1-\frac{s^2}{4}\right)^{\frac{n-3}{2}}
s^{n-2}g(\lambda s)\dd s,
\end{aligned}
\end{equation}
whose last integral is finite in both cases $n=2$ and $n\geq 3$, since $r^{n-2}g(\lambda
r)=o(r^{\alpha-2})$ in a neighborhood of $0$, and by the extra assumption $G(x)\leq M$ on
$B_R^\compl$ for $n=2$.
Thus, using Cauchy-Schwarz inequality, for any $u\in
L^2(\bbS^{n-1})$ we find
\begin{equation}\label{sphereharm:eq2}
\begin{aligned}
&\int_{\bbS^{n-1}} \left(\int_{\bbS^{n-1}}
\abs{u(y)}\frac{\eta_{G,1/\lambda}(\abs{x-y})}{\abs{x-y}^2}\dH_y^{n-1}
\right)^2 \dH_x^{n-1}\\
&\hphantom{\int_{\bbS^{n-1}}}
\leq \int_{\bbS^{n-1}} \left(\int_{\bbS^{n-1}} \abs{u(y)}^2
\frac{\eta_{G,1/\lambda}(\abs{x-y})}{\abs{x-y}^2}\dH_y^{n-1}\right)
\left(\int_{\bbS^{n-1}} \frac{\eta_{G,1/\lambda}(\abs{x-y})}{\abs{x-y}^2} \dH_y^{n-1}\right)
\dH^{n-1}_x\\
&\hphantom{\int_{\bbS^{n-1}}}
\leq C\int_{\bbS^{n-1}} \int_{\bbS^{n-1}} \abs{u(y)}^2
\frac{\eta_{G,1/\lambda}(\abs{x-y})}{\abs{x-y}^2}\dH_y^{n-1}\dH_x^{n-1}\\
&\hphantom{\int_{\bbS^{n-1}}}
\leq C\int_{\bbS^{n-1}} \abs{u(y)}^2 \left(\int_{\bbS^{n-1}}
\frac{\eta_{G,1/\lambda}(\abs{x-y})}{\abs{x-y}^2}\dH_x^{n-1}\right) \dH_y^{n-1}
\leq C\int_{\bbS^{n-1}} \abs{u(y)}^2 \dH_y^{n-1},
\end{aligned}
\end{equation}
for some $C=C(n,G,\lambda)$ and for every $u\in L^2(\bbS^{n-1})$.
As a consequence, by Fubini's theorem, the integral in \cref{eq:defscrR} converges for
$\calH^{n-1}$-almost every $x\in\bbS^{n-1}$, and $\scrR_{G,\lambda}(u)\in L^2(\bbS^{n-1})$
whenever $u\in L^2(\bbS^{n-1})$.\\
\textit{Case 2.} In the general case, $r^{n-2}g(\lambda r)$ is not necessarily integrable at the
origin, so that the corresponding kernel is \textsl{hypersingular}, and $\scrR_{G,\lambda} (u)$ may
not be well-defined for every \mbox{$u\in L^2(\bbS^{n-1})$}.
Nonetheless, when $u\in C^1(\bbS^{n-1})$, we have
\[
\frac{\abs{u(x)-u(y)}}{\abs{x-y}^2} \eta_{G,1/\lambda}(\abs{x-y})
\leq C\norm{u}_{C^1(\bbS^{n-1})} \frac{\eta_{G,1/\lambda}(\abs{x-y})}{\abs{x-y}},
\]
for every $x,y\in\bbS^{n-1}$, so that $\scrR_{G,\lambda}(u)(x)$ is still well-defined for every
$x\in\bbS^{n-1}$, and $\scrR_{G,\lambda}(u)\in L^2(\bbS^{n-1})$. Indeed, by
\cref{lem:intsphere} the integral
\[
\begin{aligned}
\int_{\bbS^{n-1}} \frac{\eta_{G,1/\lambda}(\abs{x-y})}{\abs{x-y}} \dH_y^{n-1}
&= \lambda^{n+1}\int_{\bbS^{n-1}} \edit{\abs{x-y}}g(\lambda \abs{x-y})\dH_y^{n-1}\\
&= \lambda^{n+1}\abs{\bbS^{n-2}}\int_0^2 \left(1-\frac{s^2}{4}\right)^{\frac{n-3}{2}}
s^{n-1}g(\lambda s)\dd s
\end{aligned}
\]
is finite, by the integrability of $G$ given by \cref{Kint} and the extra assumption on $G$ in
dimension $n=2$.
Then \cref{eq:seminormKop} is trivial by using the symmetry of $\bbS^{n-1}$ and of the kernel.
\end{proof}

Let us recall the so-called Funk-Hekke formula, from which we deduce that spherical harmonics are
actually eigenfunctions of the operator $\scrR_{G,\lambda}$ in both cases $\alpha\in(0,1]$ and
$\alpha\in(1,+\infty)$, and obtain an expression of their associated eigenvalues.

\begin{thm}[Funk-Hekke formula {\autocite[Theorem 1.7]{Samko}}]\label{thm:funk}
Let $f:(-1,1)\to \R$ be such that
\[
\int_{-1}^1 \left(1-t^2\right)^{\frac{n-3}{2}} \abs{f(t)}\dt <\infty.
\]
Then for every $Y_k\in \calS_k$ and every $x\in\bbS^{n-1}$, we have
\begin{equation}\label{funk:eq1}
\int_{\bbS^{n-1}} f(x\cdot y)Y_k(y)\dH_y^{n-1} = \mu_k Y_k(x),
\end{equation}
where $\mu_k$ is given by
\[
\mu_k=\abs{\bbS^{n-1}} \int_{-1}^1 P_{n,k}(t)f(t) (1-t^2)^{\frac{n-3}{2}} \dt,
\]
and $P_{n,k}$ denotes the Legendre polynomial of degree $d$ in dimension $n$, that is,
\[
P_{n,k}(t)=(-1)^{k}\frac{\Gamma\left(\frac{n-1}{2}\right)}{2^k\Gamma\left(k+\frac{n-1}{2}\right)}
(1-t^2)^{-\frac{n-3}{2}}\left(\frac{\dd}{\dt}\right)^k (1-t^2)^{k+\frac{n-3}{2}}.
\]
\end{thm}

\begin{prp}\label{prp:sphereharmop}
If $n=2$, assume that $G(x)\leq M$ on $B_{R}^\compl$, for some positive constants $M$ and $R$, and if
$n\geq 3$, let us just set $R:=0$ and assume no upper bound on $G$.
Then spherical harmonics are eigenfunctions of the operators $\scrR_{G,\lambda}$, for any
$\lambda>R$. More precisely, for every $k\in\N$, there exists $\mu_{G,\lambda,k}\in\R$ such that
\[
\scrR_{G,\lambda}(Y_k)=\mu_{G,\lambda,k} Y_k,\quad\forall Y_k\in\calS_k.
\]
As a consequence, for every $u\in C^1(\bbS^{n-1})$, we have
\begin{equation}\label{eq:seminormKspec}
\iint_{\bbS^{n-1}\times\bbS^{n-1}}
\frac{\abs{u(x)-u(y)}2}{\abs{x-y}^2} \eta_{G,1/\lambda}(\abs{x-y})\dH^{n-1}_x\dH^{n-1}_y
= \sum_{k=0}^\infty \sum_{i=1}^{d(k)} \mu_{G,\lambda,k}\, a^i_k(u)^2,
\end{equation}
and $\mu_{G,\lambda,k}$ is given by
\begin{equation}\label{eq:eigenvecK}
\mu_{G,\lambda,k} = \iint_{\bbS^{n-1}\times\bbS^{n-1}}
\frac{\abs{Y_k(x)-Y_k(y)}^2}{\abs{x-y}^2}\eta_{G,1/\lambda}(\abs{x-y})\dH_x^{n-1}\dH_y^{n-1},
\end{equation}
for any $Y_k\in\calS_k$ \edit{with $\norm{Y_k}_{L^2(\bbS^{n-1})}=1$}.
\end{prp}

\begin{proof}
Let $\lambda>R$.
The situation is actually easier if we add the extra assumption that $G(x)=o(\abs{x}^{\alpha-n})$ at
the origin, for some $\alpha>1$, so for the clarity of the exposition we first prove that spherical
harmonics are eigenfunctions for the operator $\scrR_{G,\lambda}$ under this extra assumption, and
we drop it completely in a second step.
In a third step we easily deduce \cref{eq:seminormKspec,eq:eigenvecK}.\\
\textit{Step 1.}
Here we assume $g(r)=o(r^{\alpha-n})$ near $0$, for some $\alpha>1$. Then we can consider the
operator
\[
\calR_{G,\lambda}(u)(x) \coloneqq \int_{\bbS^{n-1}}
\frac{u(y)}{\abs{x-y}^2}\eta_{G,1/\lambda}(\abs{x-y})\dH^{n-1}_y,
\]
and proceeding as above in \cref{sphereharm:eq1,sphereharm:eq2}, we see that it is a bounded linear
operator from $L^2(\bbS^{n-1})$ into $L^2(\bbS^{n-1})$.
Now let us see that we can just use the Funk-Hekke formula.  Recalling that
$\abs{x-y}=\sqrt{2(1-x\cdot y)}$ for every $x,y\in\bbS^{n-1}$, let us define
\begin{equation}\label{sphereharmop:eq0}
f_{G,\lambda}(t)\coloneqq\frac{\eta_{G,1/\lambda}(\sqrt{2(1-t)})}{2(1-t)},
\end{equation}
so that for any $Y_k\in\calS_k$, we have
\[
\calR_{G,\lambda}(Y_k)(x)= \int_{\bbS^{n-1}} f_{G,\lambda}(x\cdot y) Y_k(y)\dH_y^{n-1}.
\]
Changing variables, we compute
\begin{equation}\label{sphereharmop:eq1}
\begin{aligned}
\int_{-1}^1 (1-t^2)^{\frac{n-3}{2}} \abs{f_{G,\lambda}(t)}\dt
&=\lambda^{n+1}\int_0^2\left(1-\frac{s^2}{4}\right)^{\frac{n-3}{2}}g(\lambda s)s^{n-2} \dd s,
\end{aligned}
\end{equation}
which is finite in view of the fact that $g(r)=o(r^{\alpha-n})$ with $\alpha>1$ at the origin, and
by the extra assumption $G(x)\leq M$ in $B_R^\compl$ when $n=2$ (cutting the integral into two parts
as in \cref{lem:intsphere}).
Thus \cref{thm:funk} applies and shows that spherical harmonics are eigenfunctions of the operator
$\calR_{G,\lambda}$. In addition, for every $k\in\N$, all the $Y_k\in\calS_k$ share the same
eigenvalue, which we denote by $\mu_{G,\lambda,k}^*$.
Taking $Y_k=1$ in \cref{funk:eq1}, one gets in particular
\[
\mu_{G,\lambda,0}^* = \int_{\bbS^{n-1}} \frac{\eta_{G,1/\lambda}(\abs{x-y})}{\abs{x-y}^2}
\dH_y^{n-1},\quad\forall x\in\bbS^{n-1},
\]
so that $\scrR_{G,\lambda}$ is simply given by
\[
\scrR_{G,\lambda} = 2\left(\mu_{G,\lambda,0}^* \id-\calR_{G,\lambda}\right).
\]
Whence we deduce that spherical harmonics are also eigenfunctions for the operator
$\scrR_{G,\lambda}$, and the eigenvalue associated with any $Y_k\in\calS_k$ is
\[
\mu_{G,\lambda,k}\coloneqq 2(\mu_{G,\lambda,0}^*-\mu_{G,\lambda,k}^*).
\]
\textit{Step 2.}
When we do not know that $G(x)=o(\abs{x}^{\alpha-n})$ at the origin, the last integral on the
right-hand side of \cref{sphereharmop:eq1} is not necessarily finite, so we cannot apply the
Funk-Hekke formula directly on $f_{G,\lambda}$ to prove that spherical harmonic are eigenfunctions
of the operator $\calR_{G,\lambda}$.
In fact, $\calR_{G,\lambda}(u)$ is not even well-defined for every $u\in C^\infty(\bbS^{n-1})$,
since
\[
\int_{\bbS^{n-1}} \frac{\eta_{G,1/\lambda}(\abs{x-y})}{\abs{x-y}^2}\dH_y^{n-1}
\]
may not converge. Instead, we proceed similarly to \autocite[Lemmas 6.25 \& 6.26]{Samko}, defining
the operators
\begin{equation}\label{sphereharmop:eq2}
\begin{aligned}
\scrR^\eps_{G,\lambda}(u)(x)
\coloneqq & 2\int_{\bbS^{n-1}\,\cap\,\{\abs{x-y}>\eps\}}
\frac{u(x)-u(y)}{\abs{x-y}^2}\eta_{G,1/\lambda}(\abs{x-y})\dH_y^{n-1}\\
=& 2\int_{\bbS^{n-1}} (u(x)-u(y))\,f_{G,\lambda}^\eps(x\cdot y)\dH_y^{n-1},
\end{aligned}
\end{equation}
for every $\eps\in(0,1)$, where
\[
f_{G,\lambda}^\eps(t)\coloneqq f_{G,\lambda}(t)\ind_{\big(0,1-\frac{\eps^2}{2}\big)}(t),
\]
and $f_{G,\lambda}$ is given by \cref{sphereharmop:eq0}.
Now by introducing this cutoff we have removed the singularity making $f_{G,\lambda}$
nonintegrable on $(-1,1)$ and the operator hypersingular, so that we may now use the Funk-Hekke
formula with $f_{G,\lambda}^\eps$, which gives
\begin{equation}\label{sphereharmop:eq3}
\scrR^\eps_{G,\lambda}(Y_k)(x)= 2\mu_{G,\lambda,k}^\eps Y_k(x),\quad\forall Y_k\in\calS_k,~\forall
x\in\bbS^{n-1},
\end{equation}
where
\begin{equation}\label{sphereharmop:eq4}
\mu_{G,\lambda,k}^\eps = \abs{\bbS^{n-1}}\int_{-1}^{1-\frac{\eps^2}{2}} (1-P_{n,k}(t))
f_{G,\lambda}(t)(1-t^2)^{\frac{n-3}{2}}\dt,
\end{equation}
and $P_{n,k}$ is given by \cref{thm:funk}. It is known (see references in the proof of
\autocite[Lemma 6.25]{Samko}) that
\[
\abs{P_{n,k}(t)-1} \leq \frac{k(k+n-2)}{n-1}(1-t),\quad\forall t\in(-1,1),
\]
thus
\begin{equation}\label{sphereharmop:eq5}
\abs*{(1-P_{n,k}(t))f_{G,\lambda}(t)(1-t^2)^{\frac{n-3}{2}}}
\leq \frac{k(k+n-2)}{n-1} f_{G,\lambda}(t)(1+t)^{\frac{n-3}{2}}(1-t)^{\frac{n-1}{2}},
\end{equation}
for all $t\in(-1,1)$.
A change of variable gives
\[
\begin{aligned}
\int_{-1}^1 \abs{f_{G,\lambda}(t)}(1+t)^{\frac{n-3}{2}}(1-t)^{\frac{n-1}{2}}\dt
&=\frac{1}{2}\int_{0}^2\left(1-\frac{s^2}{4}\right)^{\frac{n-3}{2}} s^{n-2}\,\eta_{G,1/\lambda}(s)\dd
s,
\end{aligned}
\]
which is finite in view of the integrability of $r\mapsto r^n g(r)$ given by \cref{Kint}, and the
extra assumption on $G$ when $n=2$ (still cutting the integral as in \cref{lem:intsphere}).
Thus the integrand of the integral on the right-hand side of \cref{sphereharmop:eq4} belongs to
$L^1(-1,1)$.
Hence, letting $\eps$ go to $0$ in \cref{sphereharmop:eq4} gives
\[
\lim_{\eps\to 0}~\mu_{G,\lambda,k}^\eps=\abs{\bbS^{n-1}}\int_{-1}^{1}
(1-P_{n,k}(t))f_{G,\lambda}(t)(1-t^2)^{\frac{n-3}{2}}\dt\eqqcolon \mu_{G,\lambda,k}.
\]
Combining this with \cref{sphereharmop:eq3}, and passing to the limit in \cref{sphereharmop:eq2}
(with $u=Y_k$) yields
\[
\scrR_{G,\lambda}(Y_k) = \mu_{G,\lambda,k} Y_k,
\]
which shows that spherical harmonics are eigenfunctions of the operator $\scrR_{G,\lambda}$ in that
case too.\\
\textit{Step 3.}
Using the decomposition of $u\in L^2(\bbS^{n-1})$ in the orthonormal basis of spherical
harmonics given by \cref{eq:decompharmf}, and the spectral representation of $\scrR_{G,\lambda}$
given the previous step, from \cref{eq:seminormKop} we deduce
\[
\iint_{\bbS^{n-1}\times\bbS^{n-1}}
\frac{\abs{u(x)-u(y)}2}{\abs{x-y}^2} \eta_{G,1/\lambda}(\abs{x-y})\dH^{n-1}_x\dH^{n-1}_y
= \sum_{k=0}^\infty \sum_{i=1}^{d(k)} \mu_{G,\lambda,k}\, a^i_k(u)^2
\]
for every $u\in C^1(\bbS^{n-1})$, which is precisely \cref{eq:seminormKspec}. Taking $u=Y_k$
\edit{with $\norm{Y_k}_{L^2(\bbS^{n-1})}=1$} in \cref{eq:seminormKspec} gives~\cref{eq:eigenvecK}.
\end{proof}

\subsubsection[Stability when \texorpdfstring{$\gamma<1$}{gamma smaller than 1}]{Stability when \texorpdfstring{\forcebold{$\gamma<1$}}{gamma smaller than 1}} 

\begin{thm}\label{thm:stable}
Assume that $G$ satisfies \edit{assumptions \crefrange{Krad}{KC1}}.
If $\gamma<1$, there exists $\lambda_s=\lambda_s(n,G\edit{,\gamma})$ such that for any
$\lambda>\lambda_s$ the unit ball $B_1$ is a \textbf{stable} critical point of the functional
$\calF_{\gamma,G_\lambda}$.
\end{thm}

\begin{proof}
Consider any vector field $X\in C^\infty_c(\Ren;\Ren)$ inducing a volume-preserving flow on $B_1$.
We can use \cref{prp:var2fk}, \edit{since $G$ satisfies \cref{Korig,Ksuperdec,KC1}}, so that the
second variation of $\calF_{\gamma,G_\lambda}$ at $B_1$ is given by
\begin{equation}\label{stable:eq1}
\begin{aligned}
&\delta^2 \calF_{\gamma,G_\lambda}(B_1)[X]\\
&= \int_{\bbS^{n-1}} \abs{\grad_\tau\,\zeta_X}^2\dH^{n-1}
-\gamma\iint_{\bbS^{n-1}\times\bbS^{n-1}}
\frac{\abs{\zeta_X(x)-\zeta_X(y)}^2}{\abs{x-y}^2}\eta_{G,1/\lambda}(\abs{x-y})\dH_x^{n-1}\dH_y^{n-1}\\
&\phantom{=\int_{\bbS^{n-1}}}
-\left(c^2_{\bbS^{n-1}}-\gamma\, c_{G,\lambda,\bbS^{n-1}}^2 \right)\int_{\bbS^{n-1}}
\zeta_X^2\dH^{n-1}.
\end{aligned}
\end{equation}
It is well-known (see e.g. \autocite{Muller}) that for any $u\in H^1(\bbS^{n-1})$, we have
\begin{equation}\label{stable:eq2}
\int_{\bbS^{n-1}} \abs{\grad_\tau\, u}^2 \dH^{n-1}
= \sum_{k=1}^\infty \sum_{i=1}^{d(k)} l_k\, a_k^i(u)^2,
\end{equation}
where
\[
l_k = k(k+n-2),\quad\forall k\in\N.
\]
Since $G$ satisfies \cref{Ksuperdec}, by \cref{prp:scrRwelldef,rmk:scrRwelldef}, the operator
$\scrR_{G,\lambda}$ is well-defined on $C^1(\bbS^{n-1})$ whenever $\lambda>\lambda_s$, for some
$\lambda_s$ depending only on $G$, and the conclusion of \cref{prp:sphereharmop} holds as well.
Recall that the coordinate functions $x\mapsto x_i$ are spherical harmonics of degree $1$, thus
inserting these functions into \cref{eq:eigenvecK} and summing over $i$ gives
\begin{equation}\label{stable:eq3}
\begin{aligned}
\mu_{G,\lambda,1}
&= \frac{1}{\abs{\bbS^{n-1}}} \iint_{\bbS^{n-1}\times\bbS^{n-1}}
\eta_{G,1/\lambda}(\abs{x-y})\dH_x^{n-1}\dH_y^{n-1}\\
&= \int_{\bbS^{n-1}} \eta_{G,1/\lambda}(\abs{x-y})\dH_y^{n-1}=c_{G,\lambda,\bbS^{n-1}}^2.
\end{aligned}
\end{equation}
Furthermore, since $X$ induces a volume-preserving flow on $B_1$, we have
\[
\int_{\bbS^{n-1}} \zeta_X \dH^{n-1}=0.
\]
Writing
\[
\zeta_X = \sum_{k=0}^\infty \sum_{i=1}^{d(k)} a_k^i(\zeta_X) Y_k^i
\]
and using the well-known fact that
\[
\int_{\bbS^{n-1}} Y_k \dH^{n-1}=0
\]
for all spherical harmonic $Y_k$ of degree $k>0$, it follows that $a_0^1(\zeta_X)=0$. Combining this
with \cref{stable:eq1,stable:eq2,stable:eq3}, using \cref{prp:sphereharmop} and recalling
$c_{\bbS^{n-1}}^2=n-1=l_1$, we find
\begin{equation}\label{stable:eq3b}
\begin{aligned}
\delta^2 \calF_{\gamma,G_\lambda}(B_1)[X]
&= \sum_{k=\forcebold{1}}^\infty \sum_{i=1}^{d(k)}
\Big(l_k-l_1-\gamma\left(\mu_{G,\lambda,k} -\mu_{G,\lambda,1}\right)\Big) a_k^i(\zeta_X)^2\\
&= \sum_{k=\forcebold{2}}^\infty \sum_{i=1}^{d(k)}
\Big(l_k-\gamma\,\mu_{G,\lambda,k}+ \gamma\,\mu_{G,\lambda,1}-l_1\Big) a_k^i(\zeta_X)^2.
\end{aligned}
\end{equation}
We have seen in the proof of \cref{thm:unstable} that the family $(\eta_{G,\eps})_{\eps>0}$ is a
$(n-1)$-dimensional approximation of identity (up to multiplication by $\bfK_{2,n-1}$) and that in
dimension $n=2$, it satisfies the assumption \cref{eq:Tcond} (which is ensured by \cref{Ksuperdec}).
In particular for every $\lambda>0$, $\sup_{r>1}\, \eta_{G,1/\lambda}(r)<\infty$, so that we can
apply \cref{prp:boundK1nsph} to \cref{eq:eigenvecK} (even in dimension $n=2$), which yields
\[
\mu_{G,\lambda,k}
\leq Q(\eta_{G,1/\lambda})\int_{\bbS^{n-1}} \abs{\grad_\tau\, Y_k}^2\dH^{n-1}
= Q(\eta_{G,1/\lambda})l_k,
\]
where $Q(\eta_{G,1/\lambda})$ is given by \cref{eq:defQeta}. By \cref{lem:limQeps},
\edit{$Q(\eta_{G,1/\lambda})$ tends to $1$ as $\lambda$ goes to infinity}, thus for $\eps>0$ to be
fixed later, there exists $\lambda_s$ such that, for any $\lambda>\lambda_s$, we have
\[
\mu_{G,\lambda,k} \leq \left(1+\eps\right)l_k,
\]
whence
\begin{equation}\label{stable:eq13}
\begin{aligned}
l_k-\gamma\,\mu_{G,\lambda,k} \geq \big(1-\gamma(1+\eps)\big)l_k.
\end{aligned}
\end{equation}
As we have seen in the proof of \cref{thm:unstable} (see \cref{unstable:eq2}) we have
\[
\mu_{G,\lambda,1} = c_{G,\lambda,\bbS^{n-1}}^2=\int_{\bbS^{n-1}} \eta_{G,1/\lambda}(x-y)\dH_y^{n-1}
\xrightarrow{\lambda\to \infty}
c^2_{\bbS^{n-1}}=l_1,
\]
thus up to choosing $\lambda_s$ even larger if necessary, we have
\begin{equation}\label{stable:eq14}
\mu_{G,\lambda,1} \geq l_1(1-\eps),\quad\quad\forall \lambda>\lambda_s.
\end{equation}
Combining \cref{stable:eq13,stable:eq14}, it follows
\begin{equation}\label{stable:eq14b}
l_k-\gamma\,\mu_{G,\lambda,k} +\gamma\,\mu_{G,\lambda,1}-l_1
\geq \big(1-\gamma(1+\eps)\big)l_k +\big(\gamma(1-\eps)-1\big)l_1
\end{equation}
for all $\lambda>\lambda_s$, where $\lambda_s=(n,\gamma,G,\eps)$.
Since $\gamma<1$, we may choose $\eps$ small enough depending only on $n$ and $\gamma$ so that
\[
1-\gamma(1+\eps) > 0,
\]
thus, noting that $l_k$ is bounded from below by $l_2>0$ for all $k\geq 2$, because $(l_k)_{k\in\N}$
is an increasing sequence, \cref{stable:eq14b} yields
\[
\begin{aligned}
l_k-\gamma\,\mu_{G,\lambda,k} +\gamma\,\mu_{G,\lambda,1}-l_1
&\geq \big(1-\gamma(1+\eps)\big)l_2 +\big(\gamma(1-\eps)-1\big)l_1\\
&= \left(1-\gamma\right)(l_2-l_1) -\gamma(l_1+l_2)\eps,
\end{aligned}
\]
whenever $\lambda>\lambda_s(n,\gamma,G,\eps)$. Hence, since $\gamma<1$, choosing $\eps$ even smaller
if needed, depending only on $n$ and $\gamma$, we have
\[
l_k-\gamma\,\mu_{G,\lambda,k} +\gamma\mu_{G,\lambda,1}-l_1 > 0,
\]
for all $k\geq 2$ and every $\lambda>\lambda_s$, where $\lambda_s=\lambda_s(n,\gamma,G)$,
hence $\delta^2\calF_{\gamma,G_\lambda}(B_1)[X]\geq 0$ in view of \cref{stable:eq3b}, and the unit
ball $B_1$ is a volume-constrained stable set for $\calF_{\gamma,G_\lambda}$. 
\end{proof}
\begin{rmk}
In the proof of \cref{thm:stable}, we actually show that for $\lambda>\lambda_s$, $\delta^2
\calF_{\gamma,G_\lambda}(B_1)[X]>0$ for every vector field $X\in C^\infty_c(\Ren;\Ren)$ inducing a
volume-preserving flow on $B_1$ such that the projection of $\zeta_X$ on the vector space of
spherical harmonics of degree $d\geq 2$ is not trivial. This can actually be rephrased as the
strict positivity of a quadratic form.

Indeed, in view of the expression of $\delta^2 \calF_{\gamma,G_\lambda}(B_1)[X]$, it would be
natural to define the quadratic form $\calQ\calF_{\gamma,G_\lambda}(u)\coloneqq \calQ
P(u)-\gamma\calQ\Per_{G_\lambda}(u)$, where
\[
\calQ P(u) \coloneqq \int_{\bbS^{n-1}} \abs{\grad_\tau\,u}^2\dH^{n-1}-(n-1)\int_{\bbS^{n-1}}
u^2\dH^{n-1},
\]
and
\[
\begin{aligned}
\calQ \Per_{G_\lambda}(u)
\coloneqq & \iint_{\bbS^{n-1}\times\bbS^{n-1}}
G_\lambda(x-y)\abs{u(x)-u(y)}^2\dH_x^{n-1}\dH_y^{n-1}\\
&\hphantom{\iint}
-\iint_{\bbS^{n-1}\times\bbS^{n-1}} G_\lambda(x-y)\abs{x-y}^2
u(x)^2\dH_x^{n-1}\dH_y^{n-1},
\end{aligned}
\]
on the vector space
\[
\wt{H}^1(\bbS^{n-1})
= \left\{u\in H^1(\bbS^{n-1}) ~:~\int_{\bbS^{n-1}} u\dH^{n-1}=0\right\},
\]
(recall that $\zeta_X=X\cdot \nu_B$ is null-averaged on $\bbS^{n-1}$ whenever $X$ induces a
volume-preserving flow on $B_1$, since $B_1$ is a volume-constrained stationary set). Note that
$\wt{H}^1(\bbS^{n-1})$ is the quotient space of $H^1(\bbS^{n-1})$ by $\calS_0$, denoted by
$H^1(\bbS^{n-1})/\calS_0$, since $\calS_0$ is the space of constant of functions on $\bbS^{n-1}$.
It is worth pointing out that, by the proof of \autocite[Theorem 7.1]{FFMMM},
$\delta^2\calF_{\gamma,G_\lambda}(B_1)[X]\geq 0$ for every vector field $X$ inducing a
volume-preserving flow on $B_1$ if and only if $\calQ\calF_{\gamma,G_\lambda}$ is nonnegative on
$\wt{H}^1(\bbS^{n-1})$. The proof of \cref{thm:stable} actually shows that for $\lambda>\lambda_s$,
$\calQ\calF_{\gamma,G_\lambda}$ is \textsl{strictly positive} on $\wt{H}^1(\bbS^{n-1})/\calS_1$.
It is natural to quotient by $\calS_1$, which corresponds to infinitesimal translations of the ball,
in view of the translation invariance of the problem.

The strict stability could prove useful to study the \textsl{local minimality} of the unit ball.
However, going from stability to local minimality often requires uniform regularity results for
minimizers, as in \autocite{FFMMM,BC}, where the fact that minimizers are \textsl{uniform}
quasi-minimizers for the perimeter is necessary.
\end{rmk}

%
%
%


\appendix

\section{Asymptotics for nonlocal seminorms on the sphere}
\label{app:mollifsphere}

In this appendix, we prove Bourgain-Brezis-Mironescu-type results similar to the ones in
\cref{prp:boundK1n,prp:limK1n}, in the case where $\Ren$ is replaced by the $(n-1)$-dimensional
sphere $\partial B_1=\bbS^{n-1}$.
In \autocite{KreMor} the case of a general Riemannian manifold is considered, yet the monotonicity
of the radial kernels is required, which is too strong to be applicable in our case.
In addition, it does not provide a satisfying upper bound on the quantity
\begin{equation}\label{mollifsphere:eq1}
\iint_{\bbS^{n-1}\times\bbS^{n-1}}
\frac{\abs{u(x)-u(y)}^2}{\abs{x-y}^2}\eta_{\eps}(x-y)\dH_x^{n-1}\dH_y^{n-1},
\end{equation}
in the sense that we want a bound which is uniform in $u$, and asymptotically sharp as
$\eps$ goes to $0$.

The following is the counterpart of \cref{prp:boundK1n}.

\begin{prp}\label{prp:boundK1nsph}
Let $\eta:(0,2)\to [0,+\infty)$ be a nonnegative measurable function such that
\begin{equation}\label{boundK1nsph:hyp1}
\int_0^2 \eta(t)t^{n-2}\dd t < \infty.
\end{equation}
When $n=2$ we assume in addition that $\abs{\eta(r)}\leq M$ for all $r\in (R,2)$, for some
constants $M>0$ and $R\in (0,2)$.
Then for any $u\in H^1(\bbS^{n-1})$, we have
\begin{equation}\label{eq:boundK1nsph}
\iint_{\bbS^{n-1}\times\bbS^{n-1}} \frac{\abs{u(x)-u(y)}^2}{\abs{x-y}^2}\eta(\abs{x-y})
\dH_x^{n-1}\dH_y^{n-1}
\leq \bfK_{2,n-1}\,Q(\eta)\int_{\bbS^{n-1}} \abs{\grad_\tau\,u}^2\dH^{n-1},
\end{equation}
where
\begin{equation}\label{eq:defQeta}
Q(\eta) := \abs{\bbS^{n-2}}\int_0^\pi
\left(\frac{\theta}{2}\right)^2\frac{\left(\sin\theta\right)^{n-2}}
{\left(\sin\left(\frac{\theta}{2}\right)\right)^2}\,
\eta\left(2\sin\left(\frac{\theta}{2}\right)\right)\dd\theta
\end{equation}
is finite.
\end{prp}

\begin{proof}
By Fatou's Lemma and density of $C^1(\bbS^{n-1})$ in $H^1(\bbS^{n-1})$, we may assume that $u\in
C^1(\bbS^{n-1})$.
For $x,y\in\bbS^{n-1}$, we denote by $\overarc{xy}$ the shortest arc in $\bbS^{n-1}$ from
$x$ to $y$ ($x$ and $y$ excluded), which is well-defined and nonempty whenever $x\neq \pm y$. Note
that the set $\{(x,y)\in \bbS^{n-1}\times\bbS^{n-1}~:~ x\neq \pm y\}$ is in fact
$\calH^{2(n-1)}$-negligible. For $x\neq\pm y$, we parametrize the arc $\overarc{xy}$ by
$\gamma_{x,y}:(0,1)\to \bbS^{n-1}$ with constant speed $\abs{\gamma_{x,y}'(t)}=\calH^1(\overarc{xy})$.
Integrating, we find
\[
u(y)-u(x)=\calH^1(\overarc{xy})\int_0^1 \grad_{\tau}\, u(\gamma_{x,y}(t))\cdot e_{x,y}(t)\dt,
\]
where $e_{x,y}(t):=\frac{\gamma_{x,y}'(t)}{\abs{\gamma_{x,y}'(t)}}$ is the unit tangent vector of the
arc $\overarc{xy}$ at $\gamma_{x,y}(t)$.
By Cauchy-Schwarz inequality,
\[
\abs{u(y)-u(x)}^2 \leq \calH^1(\overarc{xy})^2\int_0^1 \abs{\grad_\tau\,
u(\gamma_{x,y}(t))\cdot e_{x,y}(t)}^2\dH^1_\xi.
\]
Hence, by Fubini's theorem we have
\begin{equation}\label{boundK1nsph:eq1}
\begin{aligned}
I&:=\iint_{\bbS^{n-1}\times\bbS^{n-1}} \frac{\abs{u(x)-u(y)}^2}{\abs{x-y}^2}\eta(\abs{x-y})
\dH_x^{n-1}\dH_y^{n-1}\\
&\leq \iint_{\bbS^{n-1}\times\bbS^{n-1}}\int_0^1
\abs{\grad_\tau\,u(\gamma_{x,y}(t))\cdot
e_{x,y}(t)}^2\frac{\calH^1(\overarc{xy})^2}{\abs{x-y}^2}\eta(\abs{x-y})\dt\dH_x^{n-1}\dH_y^{n-1}\\
&= \int_{\calN}~\abs{\grad_\tau\,u(\gamma_{x,y}(t))\cdot e_{x,y}(t)}^2
\frac{\calH^1(\overarc{xy})^2}{\abs{x-y}^2}\eta(\abs{x-y})\dH_{x,y,t}^{2n-1}
=: II,
\end{aligned}
\end{equation}
where $\calN$ is the smooth $(2n-1)$-dimensional submanifold of $\R^{2n+1}$ defined by
\[
\calN := \left\{ (x,y,t) \in \bbS^{n-1}\times\bbS^{n-1}\times(0,1)~:~ x\neq\pm y\right\}.
\]
Now we want to make the proper change of variables. Let us consider the smooth $(2n-1)$-dimensional
submanifold
\[
\calM := \left\{ (e,f,\theta,\vphi) \in \bbS^{n-1}\times\bbS^{n-1}\times(0,\pi)\times(0,\pi)
~:~ e\cdot f=0~\text{ and }~\theta+\vphi<\pi\right\}\subsq \R^{2n+2}.
\]
Now we define the smooth map $\Phi:\calM\to\calN$ by
\begin{equation}\label{boundK1nsph:defchangevar}
\Phi(e,f,\theta,\vphi) = (x,y,t) =
\left((\cos\theta)e-(\sin\theta)f,(\cos\vphi)e+(\sin\vphi)f,\frac{\theta}{\theta+\vphi}\right).
\end{equation}
Note that the map $\Phi$ is one-to-one: given $(x,y,t)\in\calN$, defining $e$ by $\gamma_{x,y}(t)$,
$f$ by the rotation of $e$ in the plane $\mathrm{Vect}(x,y)$ of angle $\pi/2$ (where the orientation
is given by the arc from $x$ to $y$), $\theta$ by the angle between $x$ and $e$, and $\vphi$ by the
angle between $y$ and $e$, one gets $(e,f,\theta,\vphi)\in\calM$ with
$\Phi(e,f,\theta,\vphi)=(x,y,t)$: see \cref{fig:changevar} to picture the situation.

\begin{figure}
\centering
\includestandalone[width=.4\textwidth]{figure1}
\caption{The situation of \cref{boundK1nsph:defchangevar} in the proof of \cref{prp:boundK1nsph}.
Here $t=\frac{\theta}{\theta+\vphi}$.}
\label{fig:changevar}
\end{figure}

Then the change of variables $(x,y,t)=\Phi(e,f,\theta,\vphi)$ gives
\begin{equation}\label{boundK1nsph:eq2}
\begin{aligned}
II
&=\int_{\Phi(\calM)} \abs{\grad_\tau\,u(\gamma_{x,y}(t))\cdot e_{x,y}(t)}^2\dH_{x,y,t}^{2n-1}\\
&=\int_{\calM} \abs{\grad_\tau\,u(e)\cdot f}^2 \abs{\det J_\Phi}\dH_{e,f,\theta,\vphi}^{2n-1}\\
&=
\iint\limits_{\substack{\theta,\vphi\in(0,\pi)\\\theta+\vphi<\pi}}~
\iint\limits_{\substack{\bbS^{n-1}\times\bbS^{n-1}\\ e\cdot f=0}}
\abs{\grad_\tau\,u(e)\cdot f}^2
\frac{(\theta+\vphi)^2}{\Big(2\sin\left(\frac{\theta+\vphi}{2}\right)\Big)^2}
\eta\Big(2\sin\left(\frac{\theta+\vphi}{2}\right)\Big)
\abs{\det J_\Phi}\dH^{2n-3}_{(e,f)}\dd\theta\dd\vphi,
\end{aligned}
\end{equation}
where we have used the facts that $e_{x,y}(t)=f$, $\calH^1(\overarc{xy})=\theta+\vphi$ and
\[
\abs{x-y}=\sqrt{2(1-\cos(\theta+\vphi))}=2\sin\left(\frac{\theta+\vphi}{2}\right).
\]
Here $J_\Phi$ denotes the matrix of the differential $D\Phi: \Tan_{\calM}\, (e,f,\theta,\vphi) \to
\Tan_{\calN}\, \Phi(e,f,\theta,\vphi)$ in orthonormal bases of the tangent spaces. After
computations that we detail further below, we get
\begin{equation}\label{boundK1nsph:eqjacPhi}
\abs{\det J_\Phi} = \frac{\left(\sin(\theta+\vphi)\right)^{n-2}}{\sqrt{2}(\theta+\vphi)}.
\end{equation}
Hence, \cref{boundK1nsph:eq2} yields
\begin{equation}\label{boundK1nsph:eq3}
\begin{multlined}
II
=
\iint\limits_{\substack{\theta,\vphi\in(0,\pi)\\\theta+\vphi<\pi}}
\frac{(\theta+\vphi)(\sin(\theta+\vphi))^{n-2}}
{\Big(2\sin\left(\frac{\theta+\vphi}{2}\right)\Big)^2}
\eta\Big(2\sin\left(\frac{\theta+\vphi}{2}\right)\Big)
\dd\theta\dd\vphi\\
\times\frac{1}{\sqrt{2}}\int\limits_{\substack{\bbS^{n-1}\times\bbS^{n-1}\\ e\cdot f=0}}
\abs{\grad_\tau\,u(e)\cdot f}^2 \dH^{2n-3}_{e,f}.
\end{multlined}
\end{equation}
Now, on one hand, the Jacobian determinant of the map $F:\bbS^{n-1}\times\bbS^{n-1}\to\R$ given by
$F(e,f)=e\cdot f$ is $\abs{\det J_F(e,f)} = \sqrt{\abs{P_e(f)}^2 +
\abs{P_f(e)}^2}=\sqrt{2}{\abs{1-e\cdot f}}$, where $P_e(f)$ and $P_f(e)$ denote respectively the
projection of $f$ on $\Tan_{\bbS^{n-1}}(e)$ and the projection of $e$ on $\Tan_{\bbS^{n-1}}(f)$.
On the other hand, for every $e\in \bbS^{n-1}$, the Jacobian determinant of the map
$G_e:\bbS^{n-1}\to\R$ given by $G_e(f)=e\cdot f$ is $\abs{\det
J_{G_e}(f)}=\abs{P_e(f)}=\abs{1-e\cdot f}$, hence
\begin{equation}\label{boundK1nsph:eq4}
\int\limits_{\substack{\bbS^{n-1}\times\bbS^{n-1}\\e \cdot f=0}} \abs{\grad_\tau\, u(e)\cdot
f}^2\dH_{e,f}^{2n-3}
=\sqrt{2}\int_{\bbS^{n-1}}\left(\int_{\bbS^{n-1}\,\cap\,\{e\}^\perp} \abs{\grad_\tau\,u(e)\cdot
f}^2\dH_f^{n-2}\right)\dH_e^{n-1}.
\end{equation}
In addition, by definition of $\bfK_{2,n-1}$, we have
\begin{equation}\label{boundK1nsph:eq5}
\int_{\bbS^{n-1}\,\cap\,\{e\}^\perp} \abs{\grad_\tau\, u(e)\cdot f}^2\dH_f^{n-2}
=\bfK_{2,n-1}\abs{\bbS^{n-2}}\abs{\grad_{e^\perp}\, u(e)}^2
\leq \bfK_{2,n-1}\abs{\bbS^{n-2}}\abs{\grad_\tau\, u(e)}^2 ,
\end{equation}
where $\grad_{e^\perp}\, u(e)$ denotes the projection of $\grad_\tau\, u(e)$ on $\{e\}^\perp$.
Finally, changing variables, we compute
\begin{equation}\label{boundK1nsph:eq6}
\iint\limits_{\substack{\theta,\vphi\in(0,\pi)\\\theta+\vphi<\pi}}
\frac{(\theta+\vphi)(\sin(\theta+\vphi))^{\frac{n-2}{2}}}
{\Big(2\sin\left(\frac{\theta+\vphi}{2}\right)\Big)^2}
\eta\Big(2\sin\left(\frac{\theta+\vphi}{2}\right)\Big)
\dd\theta\dd\vphi
= \int_0^\pi \frac{t^2\big(\sin(t)\big)^{n-2}}{\big(2\sin(t/2)\big)^2}\eta\big(2\sin(t/2)\big)\dt.
\end{equation}
Whence, plugging \cref{boundK1nsph:eq4,boundK1nsph:eq5,boundK1nsph:eq6} into \cref{boundK1nsph:eq3}
gives
\begin{equation}\label{boundK1nsph:eq7}
II = \bfK_{2,n-1} Q(\eta)\int_{\bbS^{n-1}} \abs{\grad_\tau\, u}^2\dH^{n-1},
\end{equation}
which implies \cref{eq:boundK1nsph}, in view of \cref{boundK1nsph:eq1}.
Now let us justify that $Q(\eta)$ is finite. A further change of variables gives
\begin{equation}\label{eq:Qetaaltexp}
Q(\eta)
=\int_0^2
\left(\frac{2}{t}\arcsin\left(\frac{t}{2}\right)\right)^2\left(1-\frac{t^2}{4}
\right)^{\frac{n-3}{2}}t^{n-2}\eta(t)\dd t.
\end{equation}
When $n\geq 3$, for any $t\in(0,2)$, we have $\left(1-\frac{t^2}{4}\right)^{\frac{n-3}{2}}\leq 1$,
and the function $t\mapsto \frac{1}{t}\arcsin t$ is continuous, so by \cref{boundK1nsph:hyp1} we see
that $Q(\eta)$ is finite. When $n=2$, we use the fact that $\eta(r)\leq M$ for every
$t\in(R,+\infty)$, and we cut the integral into two parts, so that
\[
Q(\eta) \leq C_R\int_0^R t^{n-2}\eta(t)\dt + C_R M<\infty,
\]
where $C_R$ is a constant depending only on $R$.
There only remains to show \cref{boundK1nsph:eqjacPhi} to conclude the proof, which is sketched
separately just below, for the sake of completeness.
\end{proof}

\begin{proof}[Proof of {\cref{prp:boundK1nsph}}, continued: Computation of $\abs{\det J_\Phi}$.]
For simplicity, we assume that $n=3$. It will be made clear what happens when $n=2$ and in higher
dimension.
Let us first give an orthormal basis of $\Tan_\calM (e,f,\theta,\vphi)$. Considering the curve
$\gamma(t)=(e(t),f(t),\theta,\vphi)\in \calM$ such that $\gamma(0)=(e,f,\theta,\vphi)$ and
$\gamma'(0)=(\sigma,\tau,0,0)$, the condition $e(t)\cdot f(t)=0$ implies that $\sigma\cdot
f+\tau\cdot e=0$, hence $(f,-e,0,0)$ belongs to $\Tan_\calM (e,f,\theta,\vphi)$. One sees then
easily that an orthonormal basis of $\Tan_\calM (e,f,\theta,\vphi)$ is given by
\begin{equation}\label{baseM}
\begin{multlined}
u_0:=\frac{1}{\sqrt{2}}(f,-e,0,0),\quad u_1:=(0,0,1,0),\quad u_2:=(0,0,0,1)\\
u_3:=(e\wedge f,0,0,0),\quad u_4:=(0,e\wedge f,0,0),
\end{multlined}
\end{equation}
where $e\wedge f$ is the cross product of the vectors $e$ and $f$.
As for $\Tan_{\calN} (x,y,t)$, let us denote by $x^\perp$ and $y^\perp$ respectively the rotations
of $x$ and $y$ of angle $\pi/2$ in the plane $\mathrm{Vect}(x,y)$, where the orientation is given by
$(x,y,x\wedge y)$. Then it is easy to see that $(x^\perp,0,0)$ and $(0,y^\perp,0)$ belong to
$\Tan_\calN (x,y,t)$, and so do $(x\wedge y,0,0)$, $(0,x\wedge y,0)$ and $(0,0,1)$, hence an
orthonormal basis is
\begin{equation}\label{baseN}
\begin{multlined}
v_0:=(x^\perp,0,0),\quad v_1:=(0,y^\perp,0),\quad v_2:=(0,0,1),\\
\quad v_3:=\frac{1}{\abs{x\wedge y}}(x\wedge y,0,0)=(e\wedge f,0,0),
\quad u_4:=\frac{1}{\abs{x\wedge y}}(0,x\wedge y,0)=(0,e\wedge f,0),
\end{multlined}
\end{equation}
where we used the fact that $\frac{x\wedge y}{\abs{x\wedge y}}=e\wedge f$.
Here, since we are in dimension $n=3$, there is only one orthogonal direction to the plane
$\mathrm{Vect}(e,f)$. In higher dimension, the cross products appearing in the bases can just be
replaced by an orthonormal basis of $\mathrm{Vect}(e,f)^\perp$, which is a $(n-2)$-dimensional
vector space. Straightforward computations give that the matrix of $D\Phi$ at $(e,f,\theta,\vphi)$
in theses bases is
\[
J_\Phi = \begin{pmatrix}
\frac{1}{\sqrt{2}} & -1 & 0 & 0 & 0\\
\frac{1}{\sqrt{2}} & 0 & 1 & 0 & 0\\
0 & \frac{\vphi}{(\theta+\vphi)^2} & \frac{-\theta}{(\theta+\vphi)^2} & 0 & 0\\
0 & 0 & 0 & \cos\theta & -\sin\theta\\
0 & 0 & 0 & \cos\vphi &\sin\vphi,
\end{pmatrix}
\]
whose determinant is $\frac{\sin(\theta+\vphi)}{\sqrt{2}(\theta+\vphi)}$. In dimension $n=2$, there
is no orthogonal direction to $\mathrm{Vect}(e,f)$, so there is no block with sine and cosine
functions, while for $n>3$, there are $(n-2)$ orthogonal directions, hence $(n-2)$ such blocks,
which gives \cref{boundK1nsph:eqjacPhi}.
\end{proof}

\begin{lem}\label{lem:limQeps}
Let $(\eta_\eps)_{\eps>0}$ be a $(n-1)$-dimensional approximation of identity.
In dimension $n=2$, assume in addition that it satisfies
\begin{equation}\label{eq:Tcond}
\sup_{r>R}~ \eta_{\eps}(r) \,\xrightarrow{\eps\to 0}\, 0,\quad\forall R\in(0,2).
\end{equation}
Then
\[
Q(\eta_\eps)\,\xrightarrow{\eps\to 0}\, 1,
\]
where $Q(\eta_\eps)$ is defined by \cref{eq:defQeta}.
\end{lem}

\begin{proof}
Let $\delta>0$, and then let $\theta_0>0$ small to be chosen later.
Let us split the outer integral defining $Q(\eta_\eps)$ into the two parts $\theta<\theta_0$ and
$\theta>\theta_0$.
The integrand is always nonnegative, and for the integral on $(0,\theta_0)$, by changing variables
as in \cref{eq:Qetaaltexp}, we have
\begin{equation}\label{limQeps:eq1}
\int_0^{\theta_0}
\frac{\theta^2\left(\sin\theta\right)^n}{\left(2\sin\left(\frac{\theta}{2}\right)\right)^2}
\eta_\eps\left(2\sin\left(\frac{\theta}{2}\right)\right)\dd\theta
=\int_0^{2\sin\left(\frac{\theta_0}{2}\right)}
\left(\frac{2}{t}\arcsin\left(\frac{t}{2}\right)\right)^2\left(1-\frac{t^2}{2}
\right)^{\frac{n-3}{2}}t^{n-2}\eta_\eps(t)\dd t.
\end{equation}
Since
\[
t\mapsto \left(\frac{2}{t}\arcsin\left(\frac{t}{2}\right)\right)^2\left(1-\frac{t^2}{2}
\right)^{\frac{n-3}{2}} 
\]
converges to $1$ as $t$ vanishes, we may choose $\theta_0$ small enough such that
\begin{equation}\label{limQeps:eq2}
\abs*{\int_0^{2\sin\left(\frac{\theta_0}{2}\right)}
\left(\frac{2}{t}\arcsin\left(\frac{t}{2}\right)\right)^2\left(1-\frac{t^2}{2}
\right)^{\frac{n-3}{2}}t^{n-2}\eta_\eps(t)\dd t
-\int_0^{2\sin\left(\frac{\theta_0}{2}\right)}
t^{n-2}\eta_\eps(t)\dd t}\leq \delta,
\end{equation}
for every $\eps>0$.
Now since $(\delta_\eps)_{\eps>0}$ is a $(n-1)$-dimensional approximation of identity, we also have
\begin{equation}\label{limQeps:eq3}
\int_0^\infty t^{n-2}\eta_\eps(t)\dd t=1,\quad\text{ and }\quad
\int_{2\sin\left(\frac{\theta_0}{2}\right)}^\infty
t^{n-2}\eta_\eps(t)\dd t\,\xrightarrow{\eps\to 0}\, 0.
\end{equation}
Thus, by \cref{limQeps:eq1,limQeps:eq2,limQeps:eq3}, we have
\begin{equation}\label{limQeps:eq4}
\abs*{\int_0^{\theta_0}
\frac{\theta^2\left(\sin\theta\right)^n}{\left(2\sin\left(\frac{\theta}{2}\right)\right)^2}
\eta_\eps\left(2\sin\left(\frac{\theta}{2}\right)\right)\dd\theta-1} \leq 2\delta,
\end{equation}
for our choice of $\theta_0$, and any $\eps$ large enough.
Now let us focus on the integral on $(\theta_0,\pi)$, and distinguish the cases $n=2$ and $n\geq 3$.
When $n=2$, we have
\begin{equation}\label{limQeps:eq5}
\begin{aligned}
\int_{\theta_0}^\pi
\frac{\theta^2\left(\sin\theta\right)^n}{\left(2\sin\left(\frac{\theta}{2}\right)\right)^2}
\eta_\eps\left(2\sin\left(\frac{\theta}{2}\right)\right)\dd\theta
&\leq C \left[\sup_{\theta\,\in\, (\theta_0,\pi)}~
\eta_\eps\left(2\sin\left(\frac{\theta}{2}\right)\right)\right]\\
&\leq C \left(\sup_{r>\delta_0}~ \eta_\eps(r)\right) \,\xrightarrow{\eps\to 0}\, 0,
\end{aligned}
\end{equation}
for some $C=C(n,\theta_0)$, where $\delta_0:=2\arcsin\left(\frac{\theta_0}{2}\right)$,
and where we used the fact that $\sin(\theta/2)$ is bounded from below and decreasing on $\theta\in
(\theta_0,\pi)$, as well as assumption \cref{eq:Tcond}.
When $n\geq 3$, a change of variable gives
\begin{equation}\label{limQeps:eq6}
\begin{aligned}
\int_{\theta_0}^\pi
\frac{\theta^2\left(\sin\theta\right)^n}{\left(2\sin\left(\frac{\theta}{2}\right)\right)^2}
\eta_\eps\left(2\sin\left(\frac{\theta}{2}\right)\right)\dd\theta
&=\int_{2\sin\left(\frac{\theta_0}{2}\right)}^2 \left(\frac{2}{t}\arcsin\left(\frac{t}{2}\right)
\right)^2\left(1-\frac{t^2}{4}\right)^{\frac{n-3}{2}}t^{n-2}\eta_\eps(t)\dt\\
&\leq C\int_{2\sin\left(\frac{\theta_0}{2}\right)}^2 t^{n-2}\eta(t)\dt
\,\xrightarrow{\eps\to 0}\, 0,
\end{aligned}
\end{equation}
since $\left(1-\frac{t^2}{4}\right)^{\frac{n-3}{2}}\leq 1$ in $(0,2)$ and $(\eta_\eps)_{\eps>0}$ is
a $(n-1)$-dimensional approximation of identity.
The result follows by gathering both cases \cref{limQeps:eq5,limQeps:eq6}, in view of
\cref{limQeps:eq4} and the arbitrariness of $\delta$.
\end{proof}

\begin{cor}\label{cor:equicont}
Let $(\eta_\eps)_{\eps>0}$ be a $(n-1)$-dimensional approximation of identity, and in dimension
$n=2$, assume that it satisfies in addition assumption \cref{eq:Tcond}.
Then there exists $\eps_0>0$ such that the maps $N_{\eps} : H^1(\bbS^{n-1}) \to [0,+\infty)$ defined by
\[
N_{\eps}(u) \coloneqq \left(\iint_{\bbS^{n-1}\times\bbS^{n-1}} \frac{\abs{u(x)-u(y)}^2}{\abs{x-y}^2}
\eta_{\eps}(\abs{x-y})\dx\dy\right)^{\frac{1}{2}},
\]
for any $0<\eps<\eps_0$, are uniformly equicontinuous.
\end{cor}

\begin{proof}
In dimension $2$, since the family $(\eta_\eps)_{\eps>0}$ satisfies \cref{eq:Tcond}, we have in
particular
\[
\sup_{r>1}~ \eta_\eps(r)<\infty.
\]
Thus, by \cref{prp:boundK1nsph} and the reverse triangle inequality, we have
\[
\abs{N_{\eps}(g)-N_{\eps}(h)} \leq N_{\eps}(g-h)\leq
\big(\bfK_{2,n-1}Q(\eta_{\eps})\big)^{\frac{1}{2}} \norm{g-h}_{H^1(\bbS^{n-1})},
\]
for all $g,h\in H^1(\bbS^{n-1})$. Then by \cref{lem:limQeps}, $(Q_{\eps})_{k\in\N}$ vanishes as
$\eps$ goes to $0$, whence the uniform equicontinuity of the maps $N_{\eps}$ when $\eps<\eps_0$,
for some positive $\eps_0$.
\end{proof}

We can now compute the limit of \cref{mollifsphere:eq1} and prove \cref{prp:limK1nsph}, which is a
counterpart of \cref{prp:limK1n} for the sphere.

\begin{proof}[Proof of {\cref{prp:limK1nsph}}]
Let $u\in H^1(\bbS^{n-1})$. Since $C^1(\bbS^{n-1})$ is dense in $H^1(\bbS^{n-1})$ and, by
\cref{cor:equicont}, the functions $N_\eps$ are equicontinuous whenever $\eps<\eps_0$ for some
positive $\eps_0$, we shall assume without loss of generality that $u\in C^1(\bbS^{n-1})$.
Let us define the quantities
\[
I(\eta) 
:= \iint_{\bbS^{n-1}\times\bbS^{n-1}} \frac{\abs{u(x)-u(y)}^2}{\abs{x-y}^2}
\eta(\abs{x-y}) \dH^{n-1}_x\dH^{n-1}_y,
\]
and
\[
II(\eta) := \bfK_{2,n-1}Q(\eta)\int_{\bbS^{n-1}} \abs{\grad_\tau\, u}^2\dH^{n-1},
\]
whenever $\eta$ satisfies the assumptions of \cref{prp:boundK1nsph}.
The aim is then to show
\begin{equation}\label{limK1nsph:eq1}
E(\eta_\eps) := I(\eta_\eps)-II(\eta_\eps)\,\xrightarrow{\eps\to 0}\, 0.
\end{equation}
Let us cut $\eta_\eps$ into two parts: one for the contribution close to the origin $\eta_\eps^c :=
\eta_\eps \ind_{(0,r)}$, and one for the contribution far from the origin $\eta_\eps^f := \eta_\eps
\ind_{(r,\infty)}$, for some small $r\in(0,2)$ to be chosen later. Note that
\[
I(\eta_\eps)=I(\eta_\eps^c)+I(\eta_{\eps}^f),\quad\text{ and }\quad
II(\eta_\eps)=II(\eta_\eps^c)+II(\eta_{\eps}^f).
\]
One the one hand, for any $x\in \bbS^{n-1}$, we have
\[
I(\eta_{\eps}^f)
\leq C_r \norm{u}_{\infty,\bbS^{n-1}}^2
\int_{\bbS^{n-1}} \eta_{\eps}^f(\abs{x-y})\dH_y^{n-1}
= C_r \norm{u}_{\infty,\bbS^{n-1}}^2
\int_r^2 \left(1-\frac{t^2}{4}\right)^{\frac{n-3}{2}}t^{n-2}\eta_{\eps}(t)\dt,
\]
where we used \cref{lem:intsphere} for the last equality.
When $n\geq 3$, this immediately implies
\[
I(\eta_{\eps}^f)
\leq C_r \norm{u}_{\infty,\bbS^{n-1}}^2
\int_r^2 \eta_{\eps}(t)\dt \,\xrightarrow{\eps\to 0}\, 0,
\]
since $(\eta_\eps)_{\eps>0}$ is a $(n-1)$-dimensional approximation of identity.
When $n=2$, we split the integral into two parts and use \cref{eq:Tcond} to obtain
\[
I(\eta_{\eps}^f)
\leq C_{r,u} \left[\int_r^1 \eta_{\eps}(t)\dt+ \left(\sup_{r>1}\, \eta_{\eps}(r)\right)\right]
\,\xrightarrow{\eps\to 0}\, 0,
\]
Thus, in both cases, one has
\[
I(\eta_{\eps}^f) \,\xrightarrow{\eps\to 0}\, 0.
\]
On the other hand, by \cref{lem:limQeps}, we have
\[
II(\eta_{\eps}^f) \,\xrightarrow{\eps\to 0}\, 0.
\]
Hence \cref{limK1nsph:eq1} amounts to showing that for any $\delta_0>0$, we may find $r$ small enough
so that
\begin{equation}\label{limK1nsph:eq2}
E(\eta_\eps^c) = I(\eta_\eps^c)-II(\eta_\eps^c)\leq \delta_0,\quad\forall \eps>0.
\end{equation}
Now, by the proof of \cref{prp:boundK1nsph} (in particular \cref{boundK1nsph:eq1,boundK1nsph:eq7},
with $\eta=\eta_\eps^c$), we have
\[
II(\eta_\eps^c)
= \iint_{\bbS^{n-1}\times\bbS^{n-1}} \left(\int_0^1 \abs{\grad_\tau\,
u(\gamma_{x,y}(t)}^2\dt\right) \frac{\calH^1(\overarc{xy})^2}{\abs{x-y}^2}
\eta_\eps^c(\abs{x-y})\dH_x^{n-1}\dH_y^{n-1},
\]
so that, using the identity
\[
\left(\int_0^1 f(t)\dt\right)^2-\int_0^1 f(t)^2\dt
= -\frac{1}{2} \int_0^1 (f(t)-f(s))^2\dt\dd s
\]
with $f(t):= \grad_\tau\, u(\gamma_{x,y}(t))\cdot e_{x,y}(t)$, we find
\begin{equation}\label{limK1nsph:eq3}
\begin{multlined}
I(\eta_\eps^c)-II(\eta_\eps^c)\\
=-\frac{1}{2}
\iint_{\bbS^{n-1}\times\bbS^{n-1}} \left(\int_0^1\int_0^1
\abs{\grad_\tau\,u(\gamma_{x,y}(t))\cdot{e_{x,y}(t)}-\grad_\tau\,
u(\gamma_{x,y}(s))\cdot{e_{x,y}(s)}}^2 \dt\dd s\right)\\
\frac{\calH^1(\overarc{xy})^2}{\abs{x-y}^2}\eta_{\eps}^c(\abs{x-y}) \dH_x^{n-1}\dH_y^{n-1}.
\end{multlined}
\end{equation}
Notice that the integrand in \cref{limK1nsph:eq3} vanishes when $\abs{x-y}>r$. Recall that $u\in
C^1(\bbS^{n-1})$ and that $e_{x,y}(t)$ is the unit tangent vector of $\bbS^{n-1}$ at
$\gamma_{x,y}(t)$. Thus, given any small $\delta>0$, we may choose $r$ small enough so that
\[
\abs{e_{x,y}(t)-e_{x,y}(s)}\norm{\grad_\tau\, u}_{\infty,\bbS^{n-1}}\leq \delta,
\]
and also, by uniform continuity of $\grad_\tau\, u$ on $\bbS^{n-1}$, so that
\[
\abs{\grad_\tau\,u(\gamma_{x,y}(t))-\grad_\tau\,u(\gamma_{x,y}(s))}\leq \delta,
\]
for every $x,y\in \bbS^{n-1}$ s.t. $\abs{x-y}<r$ and every $s,t\in (0,1)$.
Hence, using the squared triangle inequality, \cref{limK1nsph:eq3} yields
\[
\abs{I(\eta_\eps^c)-II(\eta_\eps^c)}
\leq \delta^2\iint_{\bbS^{n-1}\times\bbS^{n-1}} \frac{\calH^1(\overarc{xy})^2}{\abs{x-y}^2}
\eta_\eps^c(\abs{x-y})\dH_x^{n-1}\dH_y^{n-1}.
\]
Then up to choosing $r$ even smaller if needed, we may assume
\[
\calH^1(\overarc{xy})\leq 2\abs{x-y},
\]
which gives
\begin{equation}\label{limK1nsph:eq4}
\abs{I(\eta_\eps^c)-II(\eta_\eps^c)}
\leq 4\delta^2\abs{\bbS^{n-1}}\int_{\bbS^{n-1}} \eta_{\eps}^c(\abs{x-e})\dH_x^{n-1},
\end{equation}
for any $e\in\bbS^{n-1}$. Finally, a change of variables shows that the integral of the right-hand
side of \cref{limK1nsph:eq4} is always finite and uniformly bounded (even in dimension $n=2$, since
$r<2$). Hence, by arbitrariness of $\delta$, \cref{limK1nsph:eq2} holds, which concludes the proof.
\end{proof}

Note that \cref{lem:limQeps} and \cref{prp:limK1nsph} justify that the upper bound given in
\cref{prp:boundK1nsph} is in some sense asymptotically sharp.

\begin{rmk}
Let us point out that as a consequence of \cref{prp:limK1nsph}, we recover the well-known
convergence of the Gagliardo-Sobolev $H^s$ seminorm to the $H^1$ seminorm on the
sphere as $s$ goes to $1$ (see e.g. \autocite[(8.4)]{FFMMM}), that is,
\[
(1-s)\seminorm{u}_{H^s(\bbS^{n-1})}^2 ~\xrightarrow{s\,\uparrow\,1}~ \frac{\om_{n-1}}{2}
\norm{\grad_\tau\, u}_{L^2(\bbS^{n-1})}^2,\quad\forall u\in H^1(\bbS^{n-1}),
\]
where
\[
\seminorm{u}^2_{H^s(\bbS^{n-1})}=\iint_{\bbS^{n-1}\times\bbS^{n-1}}
\frac{\abs{u(x)-u(y)}^2}{\abs{x-y}^{n-1+2s}} \dH_x^{n-1}\dH_y^{n-1}.
\]
Indeed, we have
\[
\begin{aligned}
&\lim_{s\,\uparrow\, 1^-}~ (1-s)\seminorm{u}_{H^s(\bbS^{n-1})}^2\\
&\phantom{\lim_{s\,\uparrow\, 1^-}~}
=\lim_{\eps\,\downarrow\,0}~\frac{\eps}{2}\seminorm{u}^2_{H^{1-\frac{\eps}{2}}(\bbS^{n-1})}\\
&\phantom{\lim_{s\,\uparrow\, 1^-}~}
=\lim_{\eps\,\downarrow\,0}~\frac{\eps}{2}\iint_{\bbS^{n-1}\times\bbS^{n-1}}
\frac{\abs{u(x)-u(y)}^2}{\abs{x-y}^{n+1-\eps}} \dH_x^{n-1}\dH_y^{n-1}\\
&\phantom{\lim_{s\,\uparrow\, 1^-}~}
=\lim_{\eps\,\downarrow\,0}~2^{\eps-1}\abs{\bbS^{n-2}}\iint_{\bbS^{n-1}\times\bbS^{n-1}}
\frac{\abs{u(x)-u(y)}^2}{\abs{x-y}^2} \eta_\eps(\abs{x-y})\dH_x^{n-1}\dH_y^{n-1},
\end{aligned}
\]
where $\eta_\eps:[0,+\infty)\to [0,+\infty)$ is defined by
\[
\eta_\eps(r) \coloneqq \frac{\eps}{2^\eps r^{n-1-\eps}\abs{\bbS^{n-2}}}\ind_{[0,2]}(r).
\]
We readily check that $(\eta_\eps)_{\eps>0}$ is a $(n-1)$-dimensional approximation of identity
which satisfies \cref{eq:Tcond}, thus
\[
\begin{multlined}
\lim_{\eps\,\downarrow 0}~
2^{\eps-1}\abs{\bbS^{n-2}}\iint_{\bbS^{n-1}\times\bbS^{n-1}}
\frac{\abs{u(x)-u(y)}^2}{\abs{x-y}^2} \eta_\eps(\abs{x-y})\dH_x^{n-1}\dH_y^{n-1}\\
= \frac{\bfK_{2,n-1}\abs{\bbS^{n-2}}}{2}\norm{\grad_\tau\, u}^2_{L^2(\bbS^{n-1})}
= \frac{\om_{n-1}}{2}\norm{\grad_\tau\,u}^2_{L^2(\bbS^{n-1})},
\end{multlined}
\]
where we used \cref{lem:expKpn} for the last equality.
\end{rmk}

\addtocontents{toc}{\protect\setcounter{tocdepth}{1}}

\printbibliography 

\vspace{10pt}

{\small\noindent\textsc{Univ. Lille, CNRS, Inria, UMR 8524 - Laboratoire Paul Painlevé, F-59000 Lille, France}\\
\href{mailto:marc.pegon@univ-lille.fr}{marc.pegon@univ-lille.fr}}

\listoftodos

\end{document}

%% file: figure1.tex
\usetikzlibrary{math}
\usetikzlibrary{arrows.meta}

\tikzmath{
\scale=2;
\textscale=0.8;
\Rx=\scale*1;
\Ry=\scale*0.4;
\th=0.7*pi;
\ph=pi/6;
\axi=deg(1.2*pi/4);
\ax=\axi-deg(\th);
\ay=\axi+deg(\ph);
\xii1=\Rx*cos(\axi);
\xii2=\Ry*sin(\axi);
\x1=\Rx*cos(\ax);
\x2=\Ry*sin(\ax);
\y1=\Rx*cos(\ay);
\y2=\Ry*sin(\ay);
\z2=\Rx*1.3;
\rex=0.27;
\arcexRx=\rex*\Rx;
\arcexRy=\rex*\Ry;
\arcexx1=\rex*\x1;
\arcexx2=\rex*\x2;
\arcexxi1=\rex*\xii1;
\arcexxi2=\rex*\xii2;
\rey=0.5;
\arceyRx=\rey*\Rx;
\arceyRy=\rey*\Ry;
\arceyy1=\rey*\y1;
\arceyy2=\rey*\y2;
\arceyxi1=\rey*\xii1;
\arceyxi2=\rey*\xii2;
\midey1=0.75*(\y1+\xii1)/2;
\midey2=0.75*(\y2+\xii2)/2;
\midex1=0.8*(\x1+\xii1)/2;
\midex2=0.8*(\x2+\xii2)/2;
\af=\axi+90;
\f1=\Rx*cos(\af);
\f2=\Ry*sin(\af);
\sqscale=0.2;
\sqe1=\sqscale*\xii1;
\sqe2=\sqscale*\xii2;
\sqf1=\sqscale*\f1;
\sqf2=\sqscale*\f2;
\sqmid1=(\sqe1+\sqf1);
\sqmid2=(\sqe2+\sqf2);
}

\tikzset{
    partial ellipse/.style args={#1:#2:#3}{
        insert path={+ (#1:#3) arc (#1:#2:#3)}
    }
}

\begin{tikzpicture}[>=Latex]

\draw[dash pattern=on 2pt off 2pt] (0,0) ellipse ({\Rx} and {\Ry});
\draw (0,0) circle ({\Rx});
\filldraw[color=black!40, fill=black!5] (0,0) [partial ellipse={\ax}:{\axi}:{\arcexRx} and {\arcexRy}];
\fill[black!5] (0,0) -- ({\arcexx1},{\arcexx2}) -- ({\arcexxi1},{\arcexxi2}) -- cycle;
\draw ({\midex1},{\midex2}) node[scale=\textscale] {$\theta$};
\filldraw[color=black!40, fill=black!5] (0,0) [partial ellipse={\axi}:{\ay}:{\arceyRx} and {\arceyRy}];
\fill[black!5] (0,0) -- ({\arceyy1},{\arceyy2}) -- ({\arceyxi1},{\arceyxi2}) -- cycle;
\draw ({\midey1},{\midey2}) node[scale=\textscale] {$\varphi$};
\draw [->] (0,0) -- (\xii1,\xii2);
\draw (\xii1,\xii2) node[scale=\textscale,outer xsep=-11pt,outer ysep=0pt,anchor=south west]
{$\gamma_{x,y}(t)=e$};
\draw [->] (0,0) -- (\f1,\f2);
\draw (\f1,\f2) node[scale=\textscale,outer xsep=-12pt,anchor=south east] {$e_{x,y}(t)=f$};
\draw[red] ({\sqe1},{\sqe2}) -- ({\sqmid1},{\sqmid2}) -- ({\sqf1,\sqf2});
\draw [->] (0,0) -- (\x1,\x2);
\draw (\x1,\x2) node[scale=\textscale,anchor=north] {$x$};
\draw [->] (0,0) -- (\y1,\y2);
\draw (\y1,\y2) node[scale=\textscale,anchor=south] {$y$};
\draw [->,semithick] (0,0) -- (0,\z2);
\draw (0,\z2) node[scale=\textscale,anchor=south] {$\mathrm{Vect}(x,y)^\perp$};
\end{tikzpicture}